\documentclass[10pt]{article}
\usepackage{amsthm,amsfonts,amssymb,epsfig,graphics}
\usepackage[latin1]{inputenc}
\usepackage[francais,english]{babel}
\usepackage{makeidx}
\makeindex

\newcommand{\A}{{\mathbb{A}}}
\def\t{\tilde }

\def\ie{\emph{i.e.} }
\def\cf{\emph{cf.} }

\def\homeo{hom\'eo\-mor\-phisme}

\def\inpa{\mathrm{IP}}

\def\com{\mathrm{Comm}}
\def\indi{\mathrm{Indice}}
\def\adhe{\mathrm{Adh}}
\def\inte{\mathrm{Int}}

\newcommand{\R}{{\mathbb{R}}}

\def\fixe{\mathrm{Fixe}}

\def\wfs{W_{\rightarrow S}}
\def\wfn{W_{\rightarrow N}}

\def\fs {(\rightarrow S)}
\def\fn {(\rightarrow N)}
\def\sf {(S \rightarrow)}
\def\nf {(N \rightarrow)}
\def\res{\textit}

\newcommand{\ii}[1]{{#1}}

\def\S{\mathbb{S}}
\newcommand{\Z}{\mathbb{Z}}

\newtheorem{theo}{Th\'eor\`eme}[section]
\newtheorem{ques}{Question}
\newtheorem*{ques*}{Question}
\newtheorem*{theo*}{Th\'eor\`eme}
\newtheorem{affi}[theo]{Affirmation}
\newtheorem{coro}[theo]{Corollaire}

\newtheorem{prop}[theo]{Proposition}
\newtheorem{lemm}[theo]{Lemme}

\def\?{$^{***}$\marginpar{?}}

\newcommand{\incn}[1][]{}
\newenvironment{demo}[1][]{\addvspace{6mm} \emph{Preuve #1
    ---~~}}{~~~$\square$\bigskip}

\newlength{\espaceavantspecialthm}
\newlength{\espaceapresspecialthm}
\newenvironment{defi}[1][]{\refstepcounter{theo} 
\vskip \espaceavantspecialthm \noindent \textbf{D\'efinition~\thetheo
#1.} }%
{\vskip \espaceapresspecialthm}

\newenvironment{rema}[1][]{\refstepcounter{theo} 
\vskip \espaceavantspecialthm \noindent \textbf{Remarque~\thetheo
#1.} }%
{\vskip \espaceapresspecialthm}

\author{Fr\'ed\'eric Le Roux \\ Universit\'e Paris Sud \\ Bat. 430 \\ 91405 Orsay Cedex   FRANCE \\ e-mail : frederic.le-roux@math.u-psud.fr}
\date{}
 \title{Dynamique des hom\'eo\-mor\-phismes de surfaces\\ Versions
 topologiques des th\'eor\`emes de la fleur de Leau-Fatou et de la
 vari\'et\'e stable}
 
\begin{document}
\sloppy
\selectlanguage{french}
\thispagestyle{empty}
\maketitle

\begin{abstract}
On étudie la dynamique d'un
homéomorphisme de surface au voisinage d'un point fixe isolé.
 Si l'indice du point fixe est strictement
plus grand que $1$, on construit une famille de pétales autour du
point fixe, alternativement attractifs et
répulsifs, ce qui généralise un énoncé de dynamique
holomorphe. Si l'indice est strictement plus petit que $1$, on
obtient une famille de branches alternativement stables et instables, ce qui
généralise un énoncé de dynamique hyperbolique.
\end{abstract}

\selectlanguage{english}
\begin{abstract}
The study of the dynamics of a surface homeomorphism in the
neighbourhood of an isolated fixed point leads us to the following results. If the
fixed point index is greater than $1$, a family of
attractive and repulsive petals is constructed, generalizing the
Leau-Fatou flower theorem
in complex dynamics. If the index is less than $1$,
we get a family of stable and unstable branches, generalizing the
stable manifold theorem in hyperbolic dynamics.

\end{abstract}
\selectlanguage{french}
\paragraph{Mots-clés}Homéomorphisme, surface, point fixe, Brouwer.
\paragraph{Classification AMS (2000)}37E30, 37C25.

\pagebreak
\footnotesize
\thispagestyle{empty}
\setcounter{tocdepth}{1}
\tableofcontents

\pagebreak
\thispagestyle{empty}
\normalsize
\addcontentsline{toc}{section}{Index}
{\small
\printindex}

\section{Introduction}
\subsection{Contexte}
Au tout d\'ebut  de l'apparition des syst\`emes dynamiques comme branche
des math\'e\-ma\-tiques, Poincar\'e avait soulign\'e l'importance des {orbites
p\'eriodiques} comme l'un des biais par lequel il est possible d'attaquer l'\'etude d'un
syst\`eme dynamique ; depuis, une bonne partie des probl\`emes de
dynamique consiste d'une part \`a chercher des orbites p\'eriodiques,
d'autre part \`a tenter de comprendre la dynamique {au voisinage} d'une
orbite p\'eriodique. Les  r\'esultats de dynamique locale sont souvent un
pr\'eliminaire incontournable \`a l'\'etude de la dynamique globale : c'est
le cas par exemple en \emph{dynamique hyperbolique} (th\'eo\-r\`eme de Hartman-Grobman, de la
vari\'et\'e stable, \'etude des bifurcations locales) ou en \emph{dynamique
holomorphe} \`a une variable (th\'eo\-r\`emes de lin\'earisation locale, th\'eo\-r\`eme de la fleur de
Leau-Fatou).

\subsection{Probl\'ematique}
Dans ce texte, nous nous int\'eressons \`a la \emph{dynamique topologique}, en dimension 2.
Que peut-on dire de la dynamique locale d'un hom\'eo\-morphisme de surface
au voisinage d'un point fixe isol\'e ?  Cette question a \'et\'e
abord\'ee il y a longtemps par G. D. Birkhoff, et plus r\'ecemment par
M. Brown,
E. Slaminka, S. Pelikan, S. Baldwin, P. Le Calvez, J.-C. Yoccoz,
S. Matsumoto (\cite{brow5,bald2,peli1,leca1,leca2,leca5}).
Contrairement aux cat\'egories avec plus de structures, il n'y a
aucun espoir de th\'eo\-r\`eme de lin\'earisation locale, ce qui explique
qu'il n'y ait pas eu de r\'eponse d\'efinitive.

Un des ingr\'edients principal est l'indice de Lefschetz du point fixe,
il permet de pr\'eciser la question initiale : quel lien y a-t-il entre
l'indice et la dynamique ?

Le lien le plus simple est cons\'equence de la th\'eorie des
hom\'eo\-morphismes de \hbox{Brouwer} : si l'indice est diff\'erent de $1$, alors
de nombreuses formes de r\'ecurrence (points p\'eriodiques, et m\^eme
points non errants) sont exclues dans un voisinage du point
fixe~; autrement dit, le \emph{comportement individuel} de chaque
orbite est tr\`es simple. Pr\'eciser le comportement dynamique 
consiste alors \`a essayer de dire
quelque chose du \emph{comportement collectif} des orbites.

\subsection{R\'esultats}
Dans ce texte, on  obtient des \'enonc\'es qui s'inspirent de th\'eo\-r\`emes
classiques des cat\'egories plus structur\'ees, le \emph{th\'eo\-r\`eme de la
vari\'et\'e stable} en dynamique hyperbolique, et le \emph{th\'eo\-r\`eme de la
fleur de Leau-Fatou} en dynamique holomorphe. 

La version topologique
du th\'eo\-r\`eme de la vari\'et\'e stable concerne le cas o\`u l'indice de Lefschetz est
strictement plus petit que $1$~: on prouve alors l'existence de $p$
\emph{branches stables} et $p$ \emph{branches instables locales},
cycliquement altern\'ees autour du point fixe (o\`u le nombre $p$ est la
diff\'erence entre $1$ et l'indice).

La version topologique du th\'eo\-r\`eme de la fleur de Leau-Fatou  concerne le cas dual,
 o\`u l'indice de Lefschetz est strictement plus grand que $1$~: on
 montre cette fois-ci l'existence de $p$ \emph{p\'etales attractifs} et
 $p$ \emph{p\'etales r\'epulsifs}, ici encore cycliquement altern\'es
 autour du point fixe (le nombre $p$ est la 
diff\'erence entre l'indice et $1$).

On obtient \'egalement une nouvelle preuve d'un r\'esultat de M. Brown~:
l'indice du point fixe pour toute puissance non nulle de
l'hom\'eo\-morphisme est \'egal \`a l'indice pour $h$.

\subsection{Strat\'egie et outils}
La strat\'egie des preuves est la suivante. Un r\'esultat d'extension
(th\'eo\-r\`eme~\ref{the.exte}) permet de se ramener \`a un cadre de dynamique
globale, celui d'un hom\'eo\-morphisme de la sph\`ere n'ayant que deux points
fixes, not\'es $N$ et $S$~: ainsi, la quasi-totalit\'e du texte portera
en r\'ealit\'e sur des questions de dynamique globale. En particulier, les
trois r\'esultats de dynamique locale se d\'eduiront rapidement d'un
m\^eme th\'eor\`eme, que nous expliquons maintenant.

Dans le cadre global, on peut utiliser la th\'eorie de \hbox{Brouwer}. Celle-ci permet
la construction, en chaque point (non fixe), d'une \emph{droite de
\hbox{Brouwer}}, c'est-\`a-dire d'un arc  dont les
extr\'emit\'es sont les points fixes (ou d'une courbe
ferm\'ee simple passant par l'un des points fixes), et qui ne rencontre
son image qu'aux points fixes. Le principal
r\'esultat de dynamique globale (th\'eo\-r\`eme~\ref{the.prin}) annonce l'existence
 d'une famille finie de
droites de \hbox{Brouwer} formant
le bord de structures 
dynamiques baptis\'ees \emph{croissants} et \emph{p\'etales}, dont le nombre et
la disposition  sont reli\'es aux
indices des points fixes $N$ et $S$.

La construction des croissants et des p\'etales
 utilise principalement deux outils.  D'une part, une notion
d'indice entre deux droites de \hbox{Brouwer} disjointes, appel\'ee \emph{indice
partiel}. D'autre part, les
\emph{d\'ecom\-po\-sitions en briques}, \'elabor\'ees et d\'ej\`a utilis\'ees par P. Le Calvez et
A. Sauzet~: il s'agit d'un certain type de triangulation du
compl\'ementaire des points fixes, adapt\'ee \`a la dynamique.
 Cette triangulation permet notamment d'\'evacuer une grande  partie des
probl\`emes d\'elicats de topologie plane, et d'obtenir des preuves
simples des r\'esultats de \hbox{Brouwer}, en construisant par une proc\'edure
 ``automatique'' des droites de
\hbox{Brouwer} \emph{simpliciales} (qui sont r\'eunion d'ar\^etes de la triangulation).

La plupart des preuves ont pour d\'ecor un autre cadre de
dynamique globale, celui des \emph{hom\'eo\-morphismes de \hbox{Brouwer}}~: ce
sont des hom\'eo\-morphismes du plan sans point fixe. Le passage du
cadre ``sph\`ere avec deux points fixes'' au cadre ``plan sans point
fixe'' s'effectue de la mani\`ere suivante~: en enlevant les deux points
fixes $N$ et $S$, puis en passant au rev\^etement universel, on
obtient une surface hom\'eomorphe au plan~; on montre alors  qu'il
existe une seule mani\`ere ``non triviale'' de relever la dynamique en un
hom\'eo\-morphisme du plan (qui est bien s\^ur sans point fixe). Ceci
est l'objet de la proposition~\ref{pro.rele}, et n\'ecessite la
construction pr\'ealable d'une droite de \hbox{Brouwer} reliant les deux points
fixes $N$ et $S$. Celle-ci est obtenue par des arguments portant
principalement sur la combinatoire des droites de \hbox{Brouwer}
simpliciales. 

On prouve notamment  une version du th\'eo\-r\`eme principal, concernant
l'existence  de
croissants et de p\'etales, dans le cadre des
hom\'eo\-morphismes de \hbox{Brouwer} (th\'eor\`eme~\ref{the.prinbrou}).
 Puis on en d\'eduit la version initiale du th\'eo\-r\`eme.

\subsection{Plan du texte}
Dans la premi\`ere partie, on \'enonce pr\'ecis\'ement les r\'esultats locaux,
et on illustre ces r\'esultats et leurs limites sur des exemples.
La deuxi\`eme partie, tr\`es courte,  contient le th\'eo\-r\`eme d'extension qui rattache le
cadre local au cadre global. La troisi\`eme \'enonce le r\'esultat principal
de dynamique globale, et construit les outils n\'ecessaires \`a sa
d\'emons\-tra\-tion (th\'eorie de \hbox{Brouwer} ``classique'',
 indice partiel, d\'ecom\-po\-sition en briques). La quatri\`eme
contient le c\oe ur de la preuve. Enfin, la derni\`ere partie applique
th\'eo\-r\`eme d'extension et r\'esultat global pour r\'ecolter les
r\'esultats de dynamique locale.

Les parties concernant le cadre global (\ref{par.enon} et
\ref{par.preu}) sont largement ind\'ependantes des autres parties.

\subsection{Remerciements}
Je tiens \`a remercier les professeurs du Coll\`ege de France, et tout
particuli\`erement Jean-Christophe Yoccoz, de m'avoir invit\'e, en janvier
2001, \`a faire
la s\'erie d'expos\'es qui m'a incit\'e \`a \'ecrire ce texte.
 Cette \'etude a pour origine
le premier sujet de th\`ese que  m'avait donn\'e Lucien Guillou,  \`a l'\'epoque
sans beaucoup de succ\`es... Le sujet a red\'emarr\'e apr\`es plusieurs
discussions avec Patrice Le Calvez, en particulier gr\^ace \`a sa suggestion d'utiliser
les ``d\'ecom\-po\-sitions en briques'', qui jouent ici un r\^ole central. 
Je remercie Lucien et Patrice pour
leur int\'er\^et constant pour ce travail. Ce texte a \'egalement b\'en\'efici\'e des
conversations avec François B\'eguin, Sylvain Crovisier, Vincent Guirardel et Duncan Sands.


\clearpage
\part{Pr\'esentation des r\'esultats de dynamique locale}\label{par.pres}
Dans cette partie, on \'enonce et on illustre les r\'esultats de dynamique
topologique locale de dimension 2. Ceux-ci sont affili\'es \`a deux th\'eo\-r\`emes c\'el\`ebres,
l'un en dynamique diff\'erentiable, l'autre en dynamique holomorphe ;
nous commençons par rappeler ces \'enonc\'es. Apr\`es avoir pr\'esent\'e les
r\'esultats nouveaux et certains de leurs anc\^etres, on d\'ecrit quelques
exemples, de complexit\'e croissante, pour la plupart d\'ej\`a
connus ; les deux derniers exemples sont nouveaux (figures~\ref{difigbis}
et \ref{ifig10}). Ces exemples ne
sont pas utilis\'es dans la suite du texte (mais en facilite sans doute
la compr\'ehension).

\section{Rappels de quelques  r\'esultats classiques}

On note $0$ le point $(0,0)$ du plan $\R^2$. On notera $\inte(E)$,
$\adhe(E)$ et  $\partial E$ l'int\'erieur, l'adh\'erence et la fronti\`ere
d'un ensemble $E$.\index{$\inte(E)$, $\adhe(E)$, $\partial E$}

Rappelons que $x$ est un
\emph{point fixe} d'une application $f$ si $f(x)=x$. L'ensemble des
points fixes de $f$ sera not\'e $\fixe(f)$.
\begin{defi}\label{def.holo}
\index{hom\'eo\-morphisme local}
Un \res{hom\'eo\-morphisme local} $h:U \rightarrow V$  est un
   hom\'eo\-morphisme $h$ entre deux voisinages connexes $U$ et $V$ de
   $0$ dans $\R^2$, pr\'eservant l'orientation, dont $0$ est un point fixe.
\end{defi}
\begin{defi}
\index{germe}
Un \res{germe d'hom\'eo\-morphisme} est une classe
d'\'equivalence d'hom\'eo\-morphismes locaux pour la relation : 
$$
(h_1,U_1,V_1) \sim (h_2,U_2,V_2) \Leftrightarrow \mbox{ il existe un voisinage $U'$ de
$0$ tel que } h_{1\mid U'}=h_{2 \mid U'}.
$$
\end{defi}
\begin{defi}
Deux hom\'eo\-morphismes locaux $(h_1,U_1,V_1)$ et $(h_2,U_2,V_2)$
   sont \res{topologiquement conjugu\'es} si
il existe un hom\'eo\-morphisme local $g: U_1 \cup V_1 \rightarrow U_2
   \cup V_2$, envoyant $U_1$ sur $U_2$ et $V_1$ sur $V_2$,
   tel que, sur $U_2$, $g \circ h_1 \circ
g^{-1}=h_2$. 
\end{defi}
Cette d\'efinition induit une
d\'efinition sur les germes. 
On d\'efinit de mani\`ere similaire les diff\'eo\-mor\-phismes et les
diff\'eo\-mor\-phismes holomorphes locaux, ainsi que leurs germes.

\`{A} partir de maintenant, tous les hom\'eo\-morphismes consid\'er\'es auront
$0$ comme unique point fixe.

\subsection{Dynamique diff\'erentiable} 

\begin{figure}[htp]
$$
\begin{array}{ccccccc}
\left( 
\begin{array}{cc} 
1/2&0\\
0&1/2
\end{array}
\right)
&\;\;\;\;\;\;\;\; &
\left( 
\begin{array}{cc} 
2&0\\
0&2
\end{array}
\right)
&\;\;\;\;\;\;\;\;&
\left( 
\begin{array}{cc} 
2&0\\
0&1/2
\end{array}
\right)
&\;\;\;\;\;\;\;\;&
\left( 
\begin{array}{cc} 
-2&0\\
0&-1/2
\end{array}
\right)
\end{array}
$$
\centerline{\hbox{\input{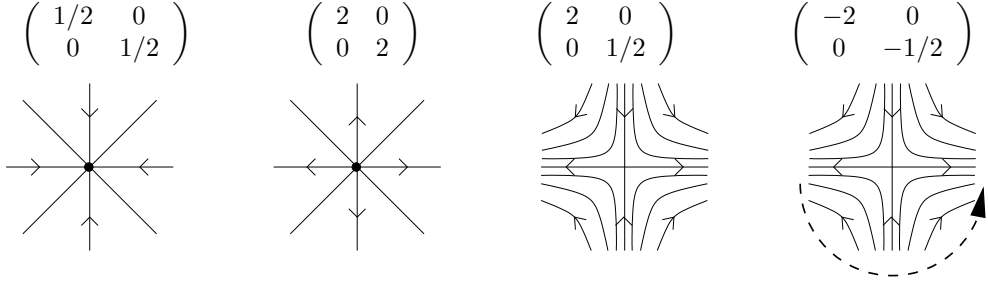}}}
\caption{\label{ifig2}Les quatre mod\`eles de diff\'eo\-mor\-phismes hyperboliques
locaux}
\end{figure}

   Soit $f$ un germe de diff\'eo\-mor\-phisme en $0$. On dit que $f$ est
   \emph{topologiquement lin\'earisable} s'il est topologiquement conjugu\'e \`a sa
   diff\'erentielle en $0$.
\begin{theo}[Hartman-Grobman, 1962]
 Supposons que la diff\'erentielle en $0$ soit
   hyperbolique (\emph{i.e.} les modules des valeurs propres sont
   diff\'erents de $1$). Alors $f$ est topologiquement lin\'earisable.  
\end{theo}
   On en d\'eduit qu'\`a conjugaison topologique
   pr\`es, il n'existe  que quatre germes de
   diff\'eo\-mor\-phismes hyperboliques (pr\'eservant
   l'orientation) : homoth\'etie contractante ou dilatante,
   point selle avec ou sans rotation (figure \ref{ifig2}).

Soit $f : U \rightarrow V$ un diff\'eo\-mor\-phisme local hyperbolique de
type selle, c'est-\`a-dire que l'une
des valeurs propres a sa valeur absolue  plus grande que $1$, et
l'autre plus petite que $1$.

Si $U'$ est un voisinage de $0$ inclus dans $U$, on note
 $W_{\rightarrow 0}^{U'}$
 et $W_{0 \rightarrow}^{U'}$ les \emph{ensembles
   stables et instables locaux}, \emph{i.e.} $W_{\rightarrow 0}^{U'}=\{x \in U' \mid \forall n
   \geq 0, f^n(x) \in U'\}$ et  $W_{0 \rightarrow}^{U'}=\{x \in U' \mid \forall n
   \leq 0, f^n(x) \in U'\}$.
Le th\'eo\-r\`eme d'Hartman-Grobman permet bien s\^ur de d\'ecrire le
comportement des orbites ; tout se passe
comme sur la figure \ref{ifig2} :
\begin{coro} \label{cor.vast}
Supposons que la diff\'erentielle est de type selle.
Alors il existe un voisinage $U'$ de $0$ arbitrairement petit, 
tel que~: 
\begin{itemize}
\item les ensembles $W_{\rightarrow 0}^{U'}$ et $W_{0 \rightarrow}^{U'}$ sont deux arcs
d'intersection $\{0\}$, d'extr\'emit\'es sur la fronti\`ere $\partial U'$~;
\item tous les points de  $W_{\rightarrow 0}^{U'}\setminus \{0\}$ ont un it\'er\'e
n\'egatif hors de $U'$ et leurs it\'er\'es positifs qui tendent vers $0$.
\item tous les points de  $W_{0 \rightarrow}^{U'}\setminus \{0\}$ ont un it\'er\'e
positif hors de $U'$ et leurs it\'er\'es n\'egatifs qui tendent vers $0$.
\end{itemize}
\end{coro}
Le th\'eo\-r\`eme de la vari\'et\'e stable ajoute que ces arcs sont
diff\'erentiables et tangents aux directions propres de la
diff\'erentielle (mais ici seuls les aspects
topologiques de la dynamique nous int\'eressent).

\subsection{Dynamique holomorphe} 
On peut trouver des pr\'ecisions sur les r\'esultats ci-dessous dans \cite{miln1,carl1}.
Soit $f$ un germe de diff\'eo\-mor\-phisme holomorphe.

Si $\mid f'(0) \mid \neq 1$ (point hyperbolique),
   alors $f$ est lin\'earisable, \emph{i.e.} la dynamique est
   localement (holomorphiquement) conjugu\'ee \`a sa
   diff\'erentielle $z \mapsto f'(0) \times z$.  Remarquons
   qu'il n' y a que deux classes de conjugaison
   topologique, correspondant aux homoth\'eties contractantes ou dilatantes ;

 Si $f'(0) = e^{i2\pi\theta}$, avec $\theta$
   irrationnel ``mal approch\'e'' par des rationnels (condition diophantienne), alors
   $f$ est encore lin\'earisable, et donc (holomorphiquement) conjugu\'e \`a
une rotation irrationnelle (Siegel, Bryuno).  

 Si $f'(0) =
   e^{i2\pi\theta}$, avec $\theta$ irrationnel ``bien
   approch\'e'' par des rationnels, alors il existe des exemples o\`u $f$ n'est pas
   lin\'earisable (Cremer, Yoccoz).
   
Le dernier cas  va nous int\'eresser plus sp\'ecialement :
   
 Si $f'(0) = e^{i2\pi\theta}$, avec $\theta$
   rationnel (point parabolique), on a encore conjugaison
\`a un mod\`ele simple (pour simplifier, on se
   place dans le cas o\`u $\theta=0$) :
\begin{theo}[Camacho, \cite{cama1}]
\label{the.cama}
 Si $f(z)=z+\alpha z^{(p+1)}+...$ ($p>0$ et $\alpha \neq 0$),
   alors  $f$ est
   topologiquement conjugu\'e \`a $z \mapsto
   z+z^{(p+1)}$.
\end{theo}
On peut  alors d\'ecrire la
dynamique locale \`a l'aide des d\'efinitions suivantes :
\begin{defi}\label{def.dito}
\index{disque topologique!ferm\'e}
 On appelle
\res{disque topologique ferm\'e} un ensemble hom\'eomorphe \`a un disque
euclidien ferm\'e du plan.
\end{defi}
\begin{defi}[ (figure~\ref{ifig1})]
\index{p\'etale}
\index{p\'etale!r\'egulier}
Soit $h:U \rightarrow V$ un hom\'eo\-morphisme local.
  Un \res{p\'etale attractif} est un  disque topologique ferm\'e $P$
  inclus dans $U$,
  dont le bord contient $0$, v\'erifiant 
$h(P) \subset \inte(P) \cup \{ 0\}$.
 Le p\'etale $P$ est dit
   \res{r\'egulier} si de plus $\bigcap_{n>0}h^n(P)=\{0\}$.
 Un \res{p\'etale r\'epulsif} est un p\'etale attractif pour $h^{-1}$.
\end{defi}
\begin{coro}[th\'eo\-r\`eme de la fleur de Leau-Fatou, 1897-1917]\label{cor.leau}
 Soit $f$ v\'eri\-fiant les hypoth\`eses du th\'eo\-r\`eme
   \ref{the.cama}. Alors, pour tout hom\'eo\-morphisme local $h:U
   \rightarrow V$ de germe $f$,
   il existe $p$ p\'etales 
   attractifs  r\'eguliers  et $p$ p\'etales r\'epulsifs r\'eguliers cycliquement
   altern\'es autour de $0$, dont la r\'eunion forme
   un voisinage de $0$.  
\end{coro}
Remarquons que le th\'eo\-r\`eme \ref{the.cama} se d\'eduit facilement du
corollaire \ref{cor.leau}.

\begin{figure}[htp]
\par
\centerline{\hbox{\input{ifig1.pstex_t}}}
\par
\caption{\label{ifig1}Dynamique de $z \mapsto z+z^3$ au voisinage de $0$}
\end{figure}

\section{R\'esultats de dynamique topologique}
On \'enonce ici les r\'esultats de dynamique topologique locale.

Certain des r\'esultats classiques ci-dessus sont  du type :

\vspace{0.2cm}

\centerline{\large Hypoth\`ese sur la diff\'erentielle $\Longrightarrow$
 conjugaison \`a un mod\`ele simple.} 

\vspace{0.2cm}

 Si on essaie d'obtenir des r\'esultats
 analogues en dynamique topologique, on se heurte \`a deux types de
 probl\`emes~:   
\begin{enumerate} 
\item au niveau des hypoth\`eses, il
   faut trouver une g\'en\'eralisation topologique
   d'hypoth\`eses portant sur la diff\'erentielle ; 
\item au
   niveau des conclusions, il n'y a aucun espoir d'obtenir
   une conjugaison \`a un petit nombre de mod\`eles (ceci sera clair sur les
   exemples d\'ecrits ci-dessous, section \ref{sec.exem}).  
\end{enumerate} 
La deuxi\`eme objection indique simplement qu'on ne pourra pas obtenir
   de th\'eo\-r\`emes de ``lin\'earisation'' (comme les th\'eo\-r\`emes de
   Camacho ou d'Hartman-Grobmann) , et qu'il faudra se
   contenter de th\'eo\-r\`emes de description partielle de la
   dynamique (comme les corollaires \ref{cor.vast} ou \ref{cor.leau}).
La premi\`ere objection est plus s\'erieuse, mais par
   miracle il existe une notion qui va permettre d'y
   rem\'edier : l'indice d'un point fixe. La d\'efinition est
   expliqu\'ee plus bas (section \ref{ss.indi}) ; pour l'instant, remarquons
   simplement deux choses : 
\begin{itemize} 
\item si
   $f(z)=z+\alpha z^{(1+p)}+\cdots$, l'indice du point
   fixe est $1+p$ (en particulier, il est strictement plus
   grand que $1$) ; 
\item pour un point selle sans
   rotation, l'indice est $-1$.
\end{itemize}
   
Si $h$ est un hom\'eo\-morphisme local fixant uniquement $0$, on notera
$\indi(h,0)$ l'indice du point fixe (qui ne d\'epend que du germe de $h$). 

\subsection{En indice $>1$, p\'etales attractifs et r\'epulsifs}
\begin{defi}
Soient $U$ un voisinage de $0$ dans le plan, $F$ et $F'$ deux
familles finies de parties compactes et connexes de $U$ contenant $0$,
 et supposons que les \'el\'ements de $F
\cup F'$ sont
deux \`a deux d'intersection \'egale \`a $\{0\}$.
 On dira que $F$ et $F'$ sont \res{cycliquement altern\'ees autour de
$0$} si il existe une courbe ferm\'ee simple $\gamma$ incluse dans $U$,
entourant $0$,
rencontrant chaque \'el\'ement de $F \cup F'$ en un unique point, telle
qu'en parcourant $\gamma$ on rencontre alternativement les \'el\'ements de
$F$ et ceux de $F'$.
\end{defi}
Remarquons qu'on pourrait \^etre plus pr\'ecis en d\'efinissant
l'ordre cyclique autour de $0$ de n'importe quelle famille finie de
parties connexes de $U$ deux \`a deux d'intersection \'egale \`a $\{0\}$.

\begin{theo}[version topologique du ``th\'eor\`eme de la fleur'' de Leau-Fatou]~
\label{the.fleu}

\nopagebreak
Soit $h:U \rightarrow V$ un hom\'eo\-mor\-phisme local fixant
  uniquement $0$, d'indice strictement sup\'erieur \`a $1$ ; on \'ecrit 
  $\indi(h,0)=1+p$ avec $p>0$.
 
\nopagebreak
Alors  il existe $p$ p\'etales attractifs et $p$ p\'etales r\'epulsifs
dans $U$, deux \`a deux d'intersection \'egale \`a $\{0\}$,
les p\'etales attractifs et r\'epulsifs \'etant cycliquement altern\'es autour de $0$.
\end{theo}

\begin{figure}[htp]
  \par
\centerline{\hbox{\input{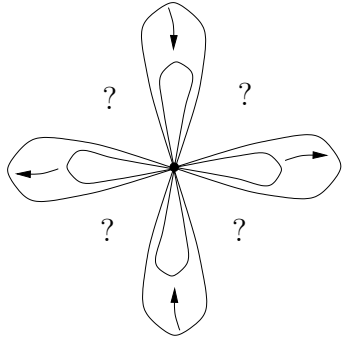}}}
\par
  \caption{\label{iifig8}Dynamique autour d'un point fixe d'indice plus
  grand que $1$ (ici d'indice $3$)}
\end{figure}

Ce th\'eo\-r\`eme est illustr\'e par la figure \ref{iifig8}.

\subsection{En indice $<1$, branches stables et instables
   locales}

Soit $h:U \rightarrow V$ un  hom\'eo\-mor\-phisme local, et $U'$ un disque
topologique ferm\'e inclus dans $\inte(U \cap V)$.
\begin{defi}
Un ensemble $E \subset U'$ est \res{plein} si $U \setminus E$ est
connexe.\footnote{Ceci revient \`a dire que toute composante connexe de $U' \setminus E$
rencontre $\partial U'$.}
\end{defi}

\begin{figure}[htp]
  \par
\centerline{\hbox{\input{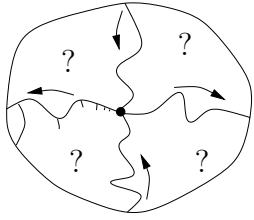}}}
\par
  \caption{\label{iifig7}Dynamique autour d'un point fixe d'indice plus
  petit que $1$ (ici d'indice $-1$)}
\end{figure}

\begin{defi}
Un ensemble $E \subset U'$
 est \res{positivement invariant} si $h(E) \subset E$. 
\end{defi}
\begin{defi}\label{def.bran}
\index{branche stable}
\index{plein}
  Une \res{branche stable $U'$-locale} est un ensemble $k$ v\'erifiant
les propri\'et\'es suivantes~:
\begin{enumerate}
\item $k$ est  inclus dans $U'$, il contient $0$ et rencontre la
fronti\`ere de $U'$~;
\item $k$  est compact, connexe, plein~;
\item $k$ est positivement invariant, et la suite des it\'er\'es positifs
de tout point de $k$ tend vers $0$~;
\item  il existe un ensemble $W_k$, voisinage de $k \setminus \{0\}$
dans $U'$, tel que tout  point de $W_k$ a un it\'er\'e n\'egatif hors de $U'$.
\end{enumerate} 
Une \res{branche instable $U'$-locale} est une branche stable
 $U'$-locale pour $h^{-1}$. 
\end{defi}
\begin{theo}[Version topologique du th\'eor\`eme de la vari\'et\'e stable]
\label{the.stlo}
Soit $h:U \rightarrow V$ un hom\'eo\-mor\-phisme local fixant
uniquement $0$, d'indice strictement inf\'erieur \`a $1$ ; on \'ecrit 
  $\indi(h,0)=1-p$ avec $p>0$. 

Alors il existe un disque topologique ferm\'e  $U'$, voisinage de $0$ 
inclus dans $\inte(U \cap V)$, $p$ branches $U'$-locales stables et $p$ branches
$U'$-locales instables, deux \`a deux d'intersection \'egale \`a $\{0\}$,
les branches stables et instables \'etant cycliquement altern\'ees autour
de $0$.
\end{theo}
Ce th\'eo\-r\`eme est illustr\'e par la figure \ref{iifig7}.

Remarquons que les p\'etales du th\'eo\-r\`eme~\ref{the.fleu} et le voisinage
$U'$ du th\'eo\-r\`eme~\ref{the.stlo} peuvent \^etre choisis arbitrairement
petits (en appliquant chaque th\'eo\-r\`eme \`a une restriction de $h$ \`a un
voisinage arbitrairement petit de $0$).

\subsection{Indices des it\'er\'es d'un hom\'eo\-morphisme
local}
On obtiendra \'egalement une nouvelle
preuve d'un th\'eo\-r\`eme de M. Brown (\cite{brow5}) :
\begin{theo}[M. Brown]
\label{the.brow}
  Soit $h:U\rightarrow V$ un  \homeo\ local fixant uniquement $0$,
dont le point fixe est d'indice diff\'erent de $1$.
 Alors il existe un voisinage $U'$ de $0$ tel que pour tout
  $n \neq 0$, 
\begin{enumerate}
\item $0$ est le seul point fixe de $h^n$ dans $U'$~;
\item $\indi(h^{n},0)=\indi(h,0)$.
\end{enumerate}
\end{theo}

   
\section{Ant\'ec\'edents}\label{ss.ante}
\subsection{La th\'eorie de \hbox{Brouwer}}
La th\'eorie des hom\'eo\-morphismes de  \hbox{Brouwer} explique notamment qu'il
  ne peut pas y avoir de r\'ecurrence
au voisinage d'un point fixe d'indice diff\'erent de $1$. Plus pr\'ecis\'ement, 
il existe un voisinage
$U$ de $0$ tel que la suite des it\'er\'es d'un point de $U$ n'a que deux possibilit\'es~:
tendre vers le point fixe ou sortir de $U$.
Ce r\'esultat est un pr\'eliminaire essentiel pour l'\'etude des dynamiques locales.
La th\'eorie de \hbox{Brouwer} est d\'etaill\'ee dans la partie \ref{sec.thbr}.
   
\subsection{Anc\^etres du th\'eo\-r\`eme de
   la fleur} En \'etudiant les points fixes d'indice strictement plus
   grand que $1$, C.~P.~Simon et N.~A.~Nikitin ont
   d'abord montr\'e que les points fixes des
   diff\'eo\-mor\-phismes pr\'eservant l'aire avaient un indice
   inf\'erieur ou \'egal \`a $1$   (\cite{simo1}, \cite{niki1}).
 S.~Pelikan et E.~E.~Slaminka (\cite{peli1})
   ont ensuite \'etendu ce r\'esultat aux hom\'eo\-morphismes.
   P.~Le Calvez a \'et\'e le premier \`a obtenir une
   propri\'et\'e de dynamique topologique qui interdit de
   pr\'eserver l'aire : en indice strictement sup\'erieur \`a $1$, 
tout voisinage du point fixe
   contient tous les it\'er\'es positifs d'un ouvert errant
   (\emph{i. e.} disjoint de tous ses it\'er\'es,
 voir \cite{leca1,leca5}). Le th\'eo\-r\`eme de la fleur
   topologique renforce ces r\'esultats, en d\'ecrivant une
   propri\'et\'e dynamique li\'ee \`a la valeur exacte de l'indice.
   
   Mais ce th\'eo\-r\`eme ne contient pas encore toutes les
   propri\'et\'es impliqu\'ees par la condition d'indice : par
   exemple, une autre propri\'et\'e dit que tout voisinage du
   point fixe contient une orbite enti\`ere (diff\'erente du point fixe).
  Ceci est une
   cons\'equence du th\'eo\-r\`eme d'indice de P. Le Calvez et J.-C.Yoccoz
   (\cite{leca2}) ;  auparavant, S. A. Andr\'ea
    l'avait prouv\'e, en indice $2$, dans un contexte global
   (\cite{andr1}).
   
\subsection{Anc\^etres du th\'eo\-r\`eme de la
   ``vari\'et\'e'' stable} Ce th\'eo\-r\`eme est extr\^emement proche
   d'un \'enonc\'e de S. Baldwin et E. E. Slaminka (\cite{bald2}) ; celui-ci
   comporte une hypoth\`ese suppl\'ementaire (l'aire est
   pr\'eserv\'ee) et une propri\'et\'e en moins dans la conclusion (les orbites
   n\'egatives des points des branches stables ne sortent
   pas n\'ecessairement du voisinage). D'autre part, la preuve semble
   utiliser implicitement une hypoth\`ese additionnelle assez forte~: l'existence
   d'un voisinage qui ne contient aucune orbite enti\`ere.
   \footnote{Voir~\cite{bald2}, page~827, ligne~7~: ``Since $x \in S_0
   - p$, there exists an $N >0$ such that $\hat h^{-N}(x) \not \in D$.''
Or sans l'hypoth\`ese additionnelle donn\'ee ci-dessus, il n'est pas
   clair que l'ensemble $S_0$ aie cette propri\'et\'e. 
D'autre part, la construction est d'abord effectu\'ee pour un
   hom\'eo\-mor\-phisme $\hat h$ qui est une perturbation de $h$, et
la mani\`ere d'obtenir les branches stables de $h$ \`a partir de celles de
   $\hat h$, qui n'est pas d\'etaill\'ee, semble \'egalement poser probl\`eme.}

   L'un des  buts de ce travail \'etait
   d'obtenir une preuve compl\`ete de cet \'enonc\'e,
   en \'evitant la  th\'eorie de la ``mise en
   position canonique'' utilis\'ee par S. Baldwin et E. E. Slaminka.
 Cette th\'eorie est due aux efforts successifs de B. S. Schmidt (\cite{schm0}),  E. E. Slaminka
   (\cite{slam1,slam2}), et M. Bonino (\cite{boni1,boni3,boni4}) ; 
 bien que ne faisant appel qu'\`a des notions
   \'el\'ementaires, les preuves sont assez d\'elicates.
La mise en position canonique peut  aussi \^etre utilis\'ee 
   pour \'etendre les germes en \'evitant l'apparition de points
   fixes, et pour calculer les
   indices des it\'er\'es $h^n$ (\cite{brow5}). On donne \'egalement
    des preuves alternatives de ces deux
   r\'esultats (th\'eo\-r\`eme \ref{the.brow} ci-dessus, et th\'eo\-r\`eme
   d'extension \ref{the.exte} ci-dessous).

K. Hiraide et J. Lewowicz ont \'egalement prouv\'e et utilis\'e un r\'esultat de
``vari\'et\'e'' stable topologique, dans un contexte diff\'erent, celui des hom\'eo\-morphismes
expansifs de surfaces (voir \cite{hira2}, \cite{lewo1}).

   Pour ce qui est des anc\^etres plus lointain, rappelons
   que dans son \'etude des hom\'eo\-morphismes de l'anneau
   d\'eviant la verticale, G. D. Birkhoff montrait et
   utilisait une propri\'et\'e de l'ensemble stable local :
   sous une hypoth\`ese tr\`es g\'en\'erale (n'\^etre ni puit ni
   source, ce qui est le cas si l'aire est pr\'eserv\'ee ou si
   l'indice est diff\'erent de $1$), la composante connexe
   du point fixe est non triviale (la preuve est tr\`es courte, voir par
   exemple \cite{lero9}, lemme~5).

\subsection{Autres r\'esultats de dynamique topologique locale} 
P. Le Calvez et J.-C. Yoccoz, puis P.~Le Calvez
 ont reli\'e la suite des indices des it\'er\'es $(h^n)$ d'un
hom\'eo\-morphisme local aux propri\'et\'es dynamiques de $h$
 (\cite{leca2,leca5}). Notamment, P.~Le Calvez retrouve le th\'eor\`eme de
 M.~Brown \'enonc\'e ci-dessus.
En dynamique globale, un  corollaire
important de cette \'etude  est la preuve de la non-existence
d'hom\'eo\-morphismes minimaux de l'anneau ouvert $\R \times \S^1$ ;
P. Le Calvez obtient \'egalement des th\'eo\-r\`emes d'existence
d'orbites p\'eriodiques.

\section{Limites des r\'esultats, questions}
\label{sub.liqu}
\subsection{En indice $<1$}
Par d\'efinition, l'ensemble stable d'un point fixe est l'ensemble des
 points dont la demi-orbite positive tend
vers le point fixe. Contrairement au cadre diff\'erentiable hyperbolique,  
pour un hom\'eo\-morphisme, l'ensemble stable  n'est pas
n\'ecessairement connexe
 (voir l'exemple du paragraphe \ref{ss.moco} ci-dessous ; cette remarque r\'epond \`a la
 question pos\'ee par  S. Baldwin et E. E. Slaminka \`a la fin de~\cite{bald2}). Le th\'eo\-r\`eme
topologique de la ``vari\'et\'e'' stable ne peut donc pas prendre en compte
l'int\'egralit\'e des ensembles stables et instables. D'autre part, il ne
 dit rien sur la dynamique entre les branches stables et instables ;
 dans quelle mesure ressemble-t-elle \`a celle des diff\'eo\-mor\-phismes
 selles ? Sur la figure \ref{iifig7},  cette absence d'information
est symbolis\'ee par les points d'interrogation dans les
secteurs entre les branches.

\subsection{En indice $>1$}
Par rapport au cadre
holomorphe du th\'eo\-r\`eme de Leau-Fatou,
 on a deux restrictions essentielles : d'une part les p\'etales ne sont pas
n\'ecessairement r\'eguliers, c'est-\`a-dire qu'on doit autoriser les p\'etales \`a contenir enti\`erement
l'orbite d'un point (voir  l'exemple du paragraphe \ref{ss.flpc}
 ci-dessous, qui ne poss\`ede aucun p\'etale r\'egulier) ; et d'autre
part on ne peut pas esp\'erer que la r\'eunion des p\'etales forme un
voisinage du point fixe (voir par exemple le deuxi\`eme dessin de la
figure \ref{ifig4}) : il resterait donc ici aussi \`a d\'ecrire la dynamique entre
les p\'etales.

\subsection{Les constructions ne sont pas canoniques}
Les constructions effectu\'ees font appara\^itre beaucoup 
de choix arbitraires (en particulier dans la construction d'une ``d\'ecom\-po\-sition
en briques''), ce qui
emp\^eche d'obtenir des objets canoniques. En essayant de rem\'edier \`a
cet inconv\'enient, on tombe notamment sur les questions suivantes :
\begin{ques}
En indice $<1$, peut-on trouver un exemple avec une infinit\'e de
branches stables ``maximales'' disjointes ?
\end{ques}
\begin{ques}
En indice $1+p>1$, disons que deux p\'etales attractifs $P_1$ et $P_2$ sont \'equivalents
si l'union des it\'er\'es n\'egatifs de $P_1$ co\"\i ncide avec l'union des
it\'er\'es n\'egatifs de $P_2$.
Peut-on montrer qu'il n'existe qu'un nombre fini
(sup\'erieur \`a $p$) de classes d'\'equivalence maximales de p\'etales attractifs ?
\end{ques}

Une des difficult\'es  est qu'en
g\'en\'eral, dans un m\^eme hom\'eo\-morphisme, on trouve simultan\'ement 
des ``zones hyperboliques'' (contribuant \`a de l'indice
n\'egatif) et des ``zones elliptiques'' (contribuant \`a de l'indice positif)
: ceci est clair sur les exemples de la section \ref{ss.simp}, figures
\ref{ifig3} et \ref{ifig4}. 
En particulier, le th\'eo\-r\`eme de la fleur dit juste que si l'indice est
$1+p>1$, on aura forc\'ement au moins $2p$ zones elliptiques ; mais il peut aussi y
avoir des zones hyperboliques, et dans ce cas, il devra y avoir 
plus de $2p$ zones elliptiques.
\begin{ques}
Peut-on donner un sens pr\'ecis aux termes ``secteurs hyperboliques'' et ``secteurs
elliptiques'' ? 
\end{ques}
Une approche de cette question sera propos\'ee par l'auteur dans
\cite{lero10}, en partant de l'id\'ee suivante : m\^eme sans savoir les d\'efinir
pr\'ecis\'ement, on peut compter les secteurs hyperboliques et elliptiques !

\section{Exemples, d\'efinition informelle de l'indice}\label{sec.exem} 
On d\'ecrit ici un certain nombre d'exemples utiles \`a la
compr\'ehension du texte ; en particulier,  le corollaire
\ref{cor.pere} et la proposition  \ref{pro.moco} montrent chacun une
limite des th\'eo\-r\`emes \ref{the.fleu} et \ref{the.stlo}.
\emph{Sur le plan logique, la suite du texte ne fait pas appel \`a cette
section.}

\subsection{Flots simples}
\label{ss.simp}
   
Un flot (diff\'erentiable) est une famille de
   diff\'eo\-mor\-phismes obtenue en int\'egrant un champ de
   vecteurs. On peut d\'efinir un \textit{flot topologique} en
   ne gardant que la condition de composition : c'est une
   famille $(h_t)_{t \in \R}$ d'hom\'eo\-morphismes v\'erifiant
   $h_t \circ h_s=h_{t+s}$ (et $h_0$ est l'identit\'e). On montre, comme dans le cas
   diff\'erentiable, que les orbites $\{h_t(x) \mid t \in
   \R\}$ forment un feuilletage du plan (\cite{whit1}).
   Les exemples les plus faciles de germes
   d'hom\'eo\-morphismes se construisent en prenant le \textit{temps
   un}, $h_1$, d'un flot $(h_t)_{t \in \R}$.

   Les flots les plus simples se construisent \`a l'aide de
   six mod\`eles de secteurs angulaires (\emph{cf}
   figure~\ref{ifig3})\footnote{On
 a repr\'esent\'e le feuilletage en orbites : pour
$x$ fix\'e, quand $t$ varie, $h_t(x)$ parcourt la feuille passant par
$x$ ; en particulier, l'hom\'eo\-mor\-phisme $h_1$ pr\'eserve chaque
feuille du feuilletage et pousse les points dans le sens des fl\`eches.}~:
\begin{figure}[htp]
\par
\centerline{\hbox{\input{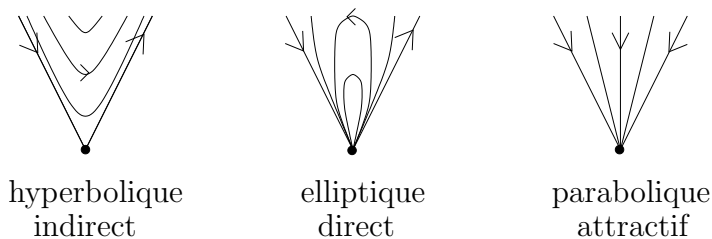}}}
\par
\caption{\label{ifig3}Mod\`eles de secteurs des flots
``simples''}
\end{figure} 
\begin{itemize}
 \item {\bf secteurs hyperboliques} direct et indirect avec lesquels on
   fabrique les points hyperboliques selles ; 
\item {\bf secteurs elliptiques} direct et indirect avec lesquels on
   fabrique les points paraboliques (!) de la
   dynamique holomorphe ; 
\item{ \bf secteurs
   paraboliques} attractif et r\'epulsif. 
\end{itemize}
On passe du mod\`ele direct au mod\`ele indirect et d'attractif \`a r\'epulsif
en changeant le sens du temps (des fl\`eches). La figure \ref{ifig4}
montre deux exemples.

\begin{figure}[htp]
\par
\centerline{\hbox{\input{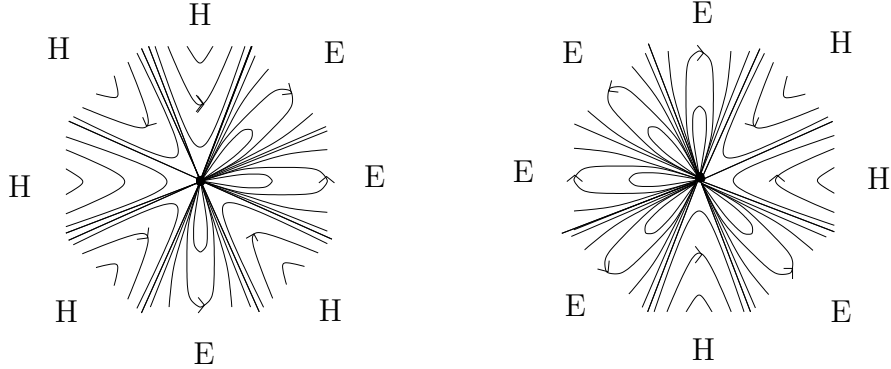}}}
\par
\caption{\label{ifig4}Deux exemples complets, duaux, de germes de flots ``simples''}
\end{figure}
   
\paragraph{}
F. Dumortier a montr\'e  que tout germe de champ de vecteurs
   $C^\infty$ du plan au voisinage d'une singularit\'e
   v\'erifiant une certaine condition (portant sur
   l'expansion de Taylor et l'existence d'une s\'eparatrice) est de
   cette forme (\cite{dumo1,dumo2}).

\subsection{L'indice}\label{ss.indi}\index{indice}
La d\'efinition sera pr\'ecis\'ee \`a la section \ref{ss.rain}.
\subsubsection*{$\bullet$ D\'efinition informelle}
\begin{figure}[hp]
\par
\centerline{\hbox{\input{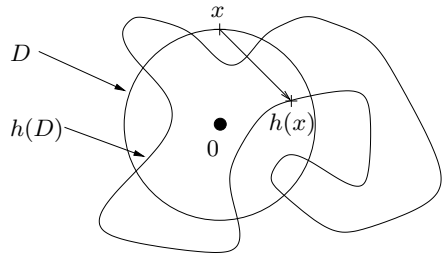}}}
\par
\caption{\label{ifig12}D\'efinition de l'indice}
\end{figure}

Soit $h: U \rightarrow V$ un hom\'eo\-mor\-phisme local fixant uniquement $0$.
 Pour obtenir
   l'indice de $h$ en $0$, on prend un disque
   $D$, inclus dans $U$, qui contient le point fixe dans son int\'erieur
 (figure \ref{ifig12}).
 En tout point
   $x$ de $\partial D$, on trace le vecteur $\overrightarrow{x
   h(x)}$, qui n'est pas nul puisque $\partial D$ ne
   contient pas de point fixe ; puis on calcule l'indice
   du champ de vecteurs ainsi d\'efini le long de $\partial
   D$ : autrement dit, le nombre (alg\'ebrique) de tours faits par
   $\overrightarrow{x h(x)}$ quand on parcourt une fois
   $\partial D$. On montre que cette d\'efinition ne d\'epend
   pas du disque $D$ choisi ; elle ne d\'epend donc que du germe de
   $h$. 

Par exemple, pour $f(z)=z+z^{p+1}$, on a $f(z)-z=z^{p+1}$ qui fait
$p+1$ tours quand $z$ parcourt le cercle unit\'e :
le point fixe $0$ est d'indice $p+1$.

\subsubsection*{$\bullet$ Formule pour les flots
   simples} Pour un flot compos\'e des quatre types de
   secteurs d\'ecrits pr\'ec\'edemment, l'indice est donn\'e par la formule :
   
$$
   \indi(h,0)=\frac{\# \mbox{secteurs elliptiques} -
   \# \mbox{secteurs hyperboliques}} {2} +1.
$$
   En effet, dans un secteur elliptique, le vecteur
   $\overrightarrow{x h(x)}$ fait un demi-tour dans le sens de
   parcours de $\partial D$, plus un ``chou\" \i a'' ; dans
   un secteur hyperbolique, il fait un demi-tour dans le
   sens contraire, moins un ``chou\" \i a'' ; dans les
   secteurs paraboliques il ne tourne presque pas (juste
   un ``chou\" \i a''). Les ``chou\" \i as'' sont en fait
   exactement les angles des secteurs, et ils
   s'additionnent pour donner le ``1'' de la formule.  
Le premier exemple  de la figure \ref{ifig4} est donc d'indice $0$,
   le deuxi\`eme est d'indice $2$.

L'int\'er\^et principal de l'indice r\'eside dans la {\bf formule
 de Lefschetz}, qui fait le lien entre le local et le global~:
 pour un \homeo\ d\'efini sur une surface
 compacte $S$, isotope \`a l'identit\'e  et
n'ayant qu'un nombre fini de points fixes, la somme des indices de ces
points fixes est \'egale \`a la caract\'eristique d'Euler-Poincar\'e de la
surface (par exemple, les deux dessins de la figure~\ref{ifig4}
 se recollent en un hom\'eo\-mor\-phisme de la sph\`ere $\S^2$, et on a bien
 $0+2=2$). 
Cette formule permet de d\'etecter des points fixes ou p\'eriodiques : par
exemple, une cons\'equence du th\'eo\-r\`eme de S. Pelikan et E. E. Slaminka (tout
germe pr\'eservant l'aire est d'indice strictement plus petit que $1$)
est que tout hom\'eo\-morphisme de la sph\`ere $\S^2$ (de caract\'eristique
 $2$), pr\'eservant l'aire, 
isotope \`a l'identit\'e, a
au moins deux points fixes\footnote{Ceci peut aussi se voir comme
 cons\'equence de la th\'eorie de \hbox{Brouwer}.}. Elle est \'egalement l'une des cl\'es de la preuve
du th\'eo\-r\`eme de P. Le Calvez et J.-C. Yoccoz sur la non-existence
d'hom\'eo\-morphismes minimaux de l'anneau.

 \subsubsection*{$\bullet$ Illustration des 
   th\'eo\-r\`emes} 
\begin{figure}[htp]
  \par
\centerline{\hbox{\input{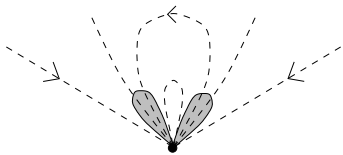}}}
\par
  \caption{\label{ifig15}P\'etales dans un secteur elliptique}
\end{figure}

Pour ces flots simples, on peut facilement
   montrer les th\'eo\-r\`emes \ref{the.stlo} et \ref{the.fleu} :

\begin{itemize}
\item {\bf ``vari\'et\'e'' stable :} soit $h$ un germe
   d'indice $1-p$ strictement plus petit que $1$, construit \`a
   l'aide des six mod\`eles de secteur d\'ecrits
   pr\'ec\'edemment. D'apr\`es la formule, il y a au moins $2p$ secteurs hyperboliques, donc au
   moins $p$ ayant la m\^eme orientation (directe
ou indirecte) ; \`a cause de l'orientation, les bords de ces secteurs
sont disjoints deux \`a deux ; ils constituent donc $p$ branches stables
et $p$ branches instables cycliquement altern\'ees.

\item {\bf fleur topologique :}
en indice $1+p$ strictement plus grand que $1$, il y a au
   moins $p$ secteurs elliptiques avec la m\^eme
   orientation ; dans chacun on trouve un p\'etale attractif
   et un r\'epulsif (voir la figure \ref{ifig15}).
\end{itemize}

Remarquons que  dans un secteur
elliptique ou parabolique, tout disque assez petit est disjoint
 de tous ses it\'er\'es positifs (ou n\'egatifs), qui sont inclus dans le
voisinage~; la dynamique ne peut donc pas pr\'eserver la mesure de
Lebesgue (ou n'importe quelle autre ``bonne'' mesure, \emph{i. e.}
non-atomique et chargeant les ouverts).
Par cons\'equent, un flot
``simple'' pr\'eservant une bonne mesure est constitu\'e uniquement de
secteurs hyperboliques.
   

\subsection{Flots compliqu\'es} Tous les germes de flots n'admettent pas
   une description aussi simple ; par exemple :
   
\begin{itemize} 
\item on peut compliquer un secteur
   hyperbolique en faisant pousser une infinit\'e de
   ``bulles'' elliptiques (\emph{cf} figure \ref{ifig5}) ; 
\begin{figure}[hp]
\par
\centerline{\hbox{\input{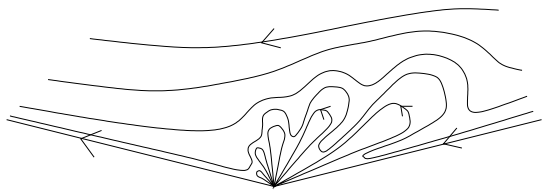}}}
\par
\caption{\label{ifig5}De l'elliptique dans l'hyperbolique}
\end{figure} 
remarquons que si on
   n'ajoute qu'un nombre fini de bulles, le germe obtenu
   est encore un flot simple (on a juste chang\'e un
   secteur hyperbolique en $k+1$ secteurs hyperboliques et
   $k$ secteurs elliptiques altern\'es).

\item on peut compliquer un secteur elliptique en y
   introduisant des ``composantes de Reeb'' (\emph{cf.} figure \ref{ifig6}) ;
\begin{figure}[htp]
\par
\centerline{\hbox{\input{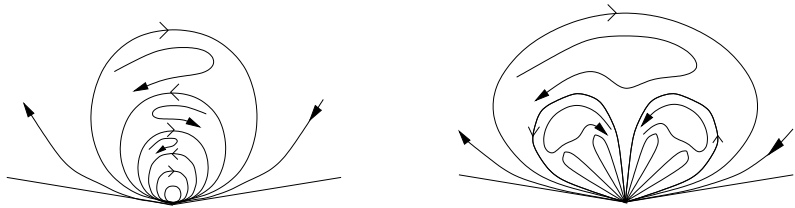}}}
\par
\caption{\label{ifig6}Des composantes de Reeb dans l'elliptique}
\end{figure}
l\`a encore, tout processus fini redonne un flot
  simple.

Les germes de feuilletages de surface au voisinage d'une singularit\'e
isol\'ee ont \'et\'e \'etudi\'es en toute g\'en\'eralit\'e par Kaplan (\cite{kapl4}).		   
\end{itemize}
   
\subsection{Flot le plus compliqu\'e possible}\label{ss.flpc}
\begin{figure}[htp]
\par
\centerline{\hbox{\input{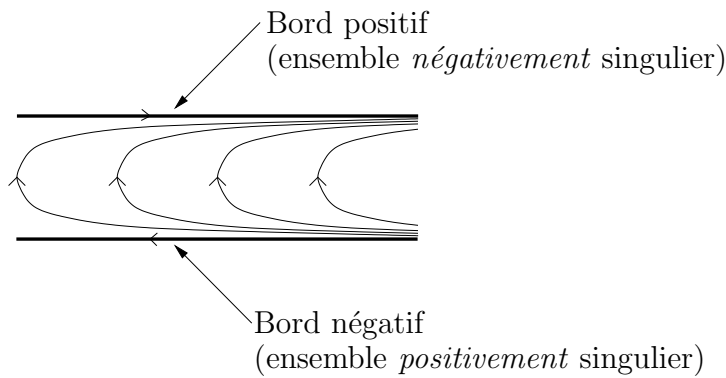}}}
\par
\caption{\label{ifig8}Une composante de Reeb}
\end{figure}

\index{singulier!ensemble}
   En premi\`ere approximation, on peut mesurer la complexit\'e d'un flot (du plan, par exemple) par son
   \textit{ensemble singulier} : celui-ci est d\'efini comme la r\'eunion
   des orbites qui ne sont pas s\'epar\'ees dans l'espace des
   feuilles du feuilletage par orbites ; ou comme
   l'ensemble des points de discontinuit\'e de l'application
   ``$x \mapsto \mbox{orbite de } x$'' (pour une topologie de
   Hausdorff) ; ou encore comme la r\'eunion des bords des
   composantes de Reeb, une composante de Reeb \'etant un
   sous-ensemble du plan, \emph{pas n\'ecessairement ferm\'e}, sur
   lequel le feuilletage est hom\'eomorphe au feuilletage de
   Reeb (repr\'esent\'e figure~\ref{ifig8}). 
L'ensemble singulier est alors la r\'eunion des ensembles 
\textit{positivement singulier} et \textit{n\'egativement singulier}, ceux-ci
   \'etant d\'efinis comme la r\'eunion des \textit{bords n\'egatifs}
   (respectivement des \textit{bords positifs}) des composantes de Reeb 
   (\emph{cf.} figure \ref{ifig8}). Remarquons encore que l'ensemble singulier
   est l'ensemble des points en lesquels la famille $\{h_1^n\}$ n'est
   pas \'equicontinue pour la m\'etrique sph\'erique, ce qui en fait
   l'analogue des ensembles de Julia de la dynamique
holomorphe. L'ensemble singulier est \'etudi\'e dans~\cite{kere2,homm1,naka1,naka2,begu1}.
  
\begin{prop}\label{pro.flde}
 Il existe un flot du plan, sans orbite p\'eriodique ni point singulier,
 dont les ensembles positivement et
   n\'egativement singuliers sont chacun denses dans le plan.
\end{prop}
   
\begin{demo} On peut trouver une construction
   compliqu\'ee due \`a T. Homma et H. Terasaka dans \cite{homm1}. Voici une construction
   \'el\'egante que m'a signal\'ee C. Bonatti. On consid\`ere
   le feuilletage du plan par les droites
   verticales. Enlevons
 les points de l'ensemble
$$
   F=\{(\frac{p}{q},q) \mid p,q \in \Z, q \neq 0 \} \  ;
$$
   il nous reste un ouvert $U$ muni d'un feuilletage $\cal F$, que
   l'on rel\`eve au rev\^etement universel $\tilde U \simeq \R^2$ pour
   obtenir un feuilletage $\tilde {\cal F}$ du plan. Soit $(\phi_t)_{t
   \in \R}$ 
n'importe quel flot de $U$ dont les orbites sont les
\begin{figure}[htp]
\par
\centerline{\hbox{\input{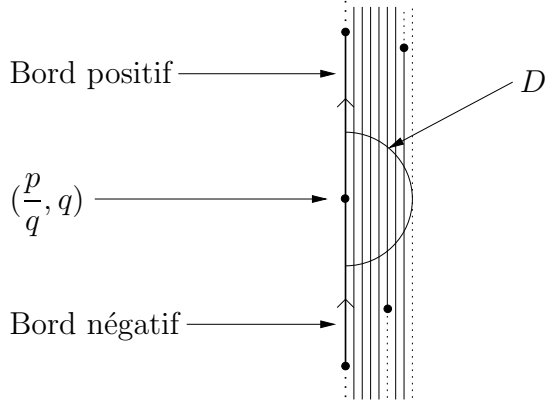}}}
\par
\caption{\label{ifig14}La r\'eunion des feuilles de $\cal F$ qui
rencontrent le demi-disque $D$ est une composante de Reeb}
\end{figure}
   feuilles de ${\cal F}$ ; alors les ensembles positivement et
   n\'egativement singuliers de $(\phi_t)_{t \in \R}$ 
   sont \'egaux, et constitu\'es de tous les
   points de $U$ d'abscisse rationnelle (figure \ref{ifig14} ;
   remarquons que la droite verticale d'abscisse $p/q$
   dans le plan contient une infinit\'e de points de l'ensemble $F$, tous les points $((kp)/(kq),kq)$).
 Les ensembles singuliers sont donc denses.
 Comme la pr\'eimage d'une composante de Reeb au
   rev\^etement universel  est une r\'eunion de
   composantes de Reeb, c'est aussi le cas pour n'importe quel flot le
   long de  $\tilde {\cal F}$.
\end{demo} 

Cet exemple montre une des limites du th\'eo\-r\`eme de la ``fleur
   topologique'' (\emph{cf} section \ref{sub.liqu}) :
\begin{coro}
\label{cor.pere}
Il   existe un hom\'eo\-morphisme local $h$ fixant uniquement $0$, d'indice $2$, pour
   lequel il n'y a ni p\'etale attractif  r\'egulier ni p\'etale r\'epulsif r\'egulier.
\end{coro}
\begin{demo}
Si $P$ est un p\'etale attractif r\'egulier et $D$ un disque inclus dans l'int\'erieur
de $P$, la suite des it\'er\'es $(h^n(D))_{n \geq 0}$ tend vers le point
fixe (par exemple au sens de la topologie de Hausdorff).

Consid\'erons alors  l'hom\'eo\-morphisme $h_1$  temps un de l'exemple de la 
proposition \ref{pro.flde}, au
voisinage du point fixe \`a l'infini. Ce point fixe est d'indice 2 (\`a
cause de la formule de Lefschetz). D'autre part, l'int\'erieur de
tout disque du plan rencontre l'ensemble positivement singulier de
$h_1$ (qui est le bord n\'egatif d'une composante de Reeb),
 donc la limite de la suite $(h_1^n(D))_{n \geq 0}$
 contient le bord positif d'une composante de Reeb : ceci montre
qu'aucun p\'etale attractif en l'infini n'est r\'egulier.
\end{demo}

Remarquons que le feuilletage $\tilde {\cal F}$
 ci-dessus poss\`ede un grand nombre de sym\'etries (tous les
automorphismes du rev\^etement utilis\'e dans la construction)~;
 par contre, on peut choisir un flot $(\phi_t)_{t \in \R}$
 le long de  $\tilde {\cal F}$ qui n'a aucune sym\'etrie, au sens o\`u les
 seuls hom\'eo\-mor\-phismes du plan commutant avec $\phi_1$ sont les
 $\phi_t$. On pourrait \'egalement construire un hom\'eo\-mor\-phisme $h$
 pr\'eservant chaque feuille de  $\tilde {\cal F}$, et ne commutant
 qu'avec ses it\'er\'es (en utilisant les id\'ees de \cite{lero5} et
 de~\cite{begu1}).

\subsection{Modification simple de flots}
\begin{figure}[htp]
\par
\centerline{\hbox{\input{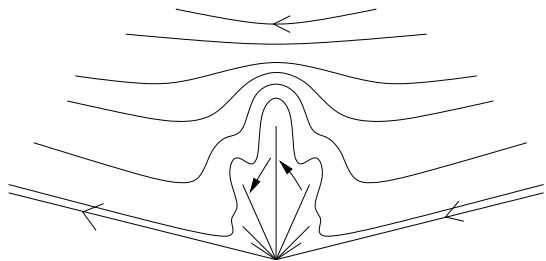}}}
\par
\caption{\label{ifig9}Cet hom\'eo\-morphisme n'est pas le temps un d'un flot}
\end{figure}
\begin{figure}[htp]
\par
\centerline{\hbox{\input{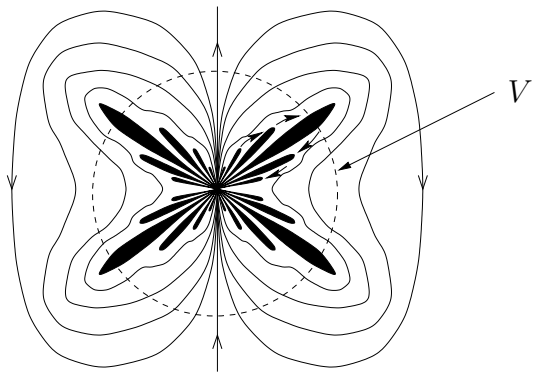}}}
\par
\caption{\label{difigbis} Il n'existe pas de ``petit'' ensemble connexe
totalement invariant}
\end{figure}

 Il existe des germes d'hom\'eo\-morphismes qui ne sont pas temps un d'un
   flot. C'est le cas de l'exemple de la figure \ref{ifig9}.
En effet, quand $h$ est temps un
   d'un flot, l'ensemble stable du point fixe (d\'efini
   comme $\{x \mid \lim_{n \rightarrow +\infty}h^n(x)=0\}$)
   est r\'eunion d'orbites du flot ; ça ne peut pas \^etre
   le cas pour ce dessin. Remarquons qu'on construit facilement une bonne
   mesure pr\'eserv\'ee par ce germe.
   
Dans le cas holomorphe parabolique, le voisinage de $0$ est feuillet\'e
par des courbes invariantes (figure \ref{ifig1}) ; en particulier,
on trouve des courbes totalement invariantes arbitrairement petites. 
Ceci est encore vrai pour des flots. Mais il n'y a pas d'\'equivalent
 de cette propri\'et\'e pour les hom\'eo\-mor\-phismes~: en effet, la figure~\ref{difigbis}
montre un germe d'hom\'eo\-morphisme, d'indice 2, pour lequel il existe un
voisinage $V$ du point fixe tel que tout
ensemble connexe,  ne contenant pas le point fixe, et totalement invariant
doit sortir de $V$ (\cite{lero1}).

\subsection{Modification compliqu\'ee de flots} 
\label{ss.moco}
Les exemples pr\'ec\'edents ressemblent encore beaucoup \`a
   des flots (notamment par le fait qu'ils se laissent
   facilement dessiner).  Les choses se compliquent
pour l'intuition (comme pour le dessinateur) 
quand des orbites ``se croisent''. Dans la proposition suivante,
 on note   $\wfn$ et $\wfs$ les ensembles
   stables de $N$ et $S$, \emph{i. e.} $\wfn=\{x \in \S^2 \mid \lim_{n
   \rightarrow +\infty}h^n(x)=N\}$.

\begin{prop}\label{pro.moco}
Pour tout entier $p$, il existe un hom\'eo\-morphisme $h$ de la sph\`ere $\S^2$, v\'erifiant :
\begin{itemize} 
\item $h$ a exactement deux points
   fixes $N$ et $S$ ; 
\item il existe un disque
   topologique $D$ tel que 
\begin{itemize} 
\item $\partial
   D \subset \wfn$ ; 
\item $D \cap \wfs \neq
   \emptyset$ ;
\end{itemize} 
\item les indices des points $N$ et $S$ sont respectivement $1+p$ et $1-p$.
\end{itemize} De plus, si $p>0$, $h$ pr\'eserve une bonne mesure au
voisinage de $S$.
   \end{prop} 

Comme pour beaucoup d'exemples, on va partir d'un flot tr\`es simple et
lui faire subir une s\'erie de \emph{modifications libres}~:
\begin{figure}[p]
\par
\centerline{\hbox{\input{ifig13.pstex_t}}}
\par
\caption{\label{ifig13}Modification libre}
\end{figure}
\begin{defi}[ (M. Brown, figure~\ref{ifig13})]
\index{modification libre}\index{support} 
\label{def.moli}
L'hom\'eo\-morphisme $h_2$ est une \res{modification libre de $h_1$ \`a
support dans $D$} si
$h_2=\phi \circ h_1$, o\`u $\phi$ est un hom\'eo\-mor\-phisme tel que~:
 \begin{enumerate}
\item le support de $\phi$ est inclus dans  $D$, autrement dit tout
point hors de $D$ est fixe par $\phi$~;
\item $D$ est un disque topologique ferm\'e d'int\'erieur \emph{libre}
pour $h_1$,
c'est-\`a-dire  v\'erifiant $h_1(\inte(D)) \cap D=\emptyset$.  
 \end{enumerate}
\end{defi}

L'int\'er\^et des modifications libres est qu'on a :
\begin{affi}\label{aff.moli}
Une modification libre ne change pas l'ensemble des points
fixes, ni les indices des points fixes isol\'es.
\end{affi}
La preuve est imm\'ediate (en particulier, on peut toujours calculer
les indices le long de courbes qui \'evitent le support de la
modification libre) .

\begin{demo}[de la proposition \ref{pro.moco}]
La figure \ref{ifig10} repr\'esente un flot sur une demi-sph\`ere, qui est
compos\'e d'un secteur hyperbolique en $S$, et en $N$ d'un secteur
elliptique avec une petite bulle elliptique.
Soit $h_1$ le temps un de ce flot, que l'on compl\`ete  sur la ``face
cach\'ee'' (demi-sph\`ere restante) par n'importe quelle dynamique sans
point fixe, de mani\`ere \`a obtenir les indices souhait\'es pour les points $N$
et $S$.
\begin{figure}[htp]
\par
\centerline{\hbox{\input{ifig10b.pstex_t}}}
\par
\caption{\label{ifig10}Flot \`a modifier}
\end{figure}
Soit $D$ un disque topologique ferm\'e comme sur la figure \ref{ifig10}, suffisamment
petit pour \^etre disjoint de son image. Remarquons que $D \subset
\wfn(h_1)$, mais que la suite $(h_1^n)_{n\geq 0}$ ne converge pas
uniform\'ement vers $N$ sur $D$ (le bord de la bulle est dans l'ensemble
positivement singulier de $h_1$).

La figure \ref{ifig11}
repr\'esente les it\'er\'es positifs de $D$ au voisinage du point fixe $S$.
Soit $\gamma$ un petit arc issu de $S$, transverse au flot, qui
rencontre  tous les it\'er\'es $h_1^n(\inte(D))$ pour ${n \geq
0}$ ; on choisit
pour chaque $n \geq 0$ un point $x_n$ sur $\inte(h_1^n(D)) \cap \gamma$.
\begin{figure}[htp]
\par
\centerline{\hbox{\input{ifig11b.pstex_t}}}
\par
\caption{\label{ifig11}Modifications libres successives}
\end{figure}

Le point $x_0$ tente de s'enfuir vers $N$ \`a chaque it\'eration par
$h_1$, le but du jeu va \^etre de le ramener vers $S$ par des
modifications libres.
On choisit pour chaque $n \geq 1$ un disque $D_n$ inclus
dans $h_1^n(D)$ et contenant $h_1(x_{n-1})$ et $x_n$ ; ainsi qu'un
hom\'eo\-morphisme $\phi_n$ qui est l'identit\'e hors de $D_n$ et qui envoie
$h_1(x_{n-1})$ sur $x_n$.
On pose alors $h_2=\phi_1 \circ h_1$, et par r\'ecurrence $h_{n+1}=\phi_n
\circ h_n$. D'apr\`es l'affirmation \ref{aff.moli}, $h_n$ n'a pas d'autres  points
fixes que $N$ et $S$, avec le m\^eme indice que $h_1$. D'autre part,
on a pu choisir les $D_n$ de mani\`ere \`a ce que la suite $(D_n)$ tende
vers $S$ (pour la topologie de Hausdorff), et dans ce cas la suite
$(h_n)$ converge uniform\'ement vers un hom\'eo\-morphisme $h_\infty$.

Cet hom\'eo\-morphisme r\'epond au probl\`eme ; en effet : 
\begin{itemize}
\item on n'a pas modifi\'e
la dynamique sur les it\'er\'es de $\partial D$ ; on a donc, comme pour
$h_1$, $\partial D \subset \wfn(h_\infty)$ ;

\item la suite $(x_n)_{n \geq 1}$ est devenue une demi-orbite positive de
$h_\infty$ ; donc $x_0$ est un point de $D$ dont les it\'er\'es positifs
tendent vers $S$.

\item si $p>0$, on a pu partir d'un flot $h_1$ qui, au point $S$, est un
germe de flot ``simple'' ne comportant que des secteurs hyperboliques
(voir section \ref{ss.simp}) ; ce flot pr\'eserve une bonne mesure au
voisinage de $S$. Pour
que $h_\infty$ pr\'eserve encore cette mesure, il suffit alors qu'elle
soit pr\'eserv\'ee par chacun des hom\'eo\-morphismes $\phi_n$ ; ceci est
possible (on montre tr\`es facilement que le groupe des hom\'eo\-morphismes du
disque pr\'eservant la mesure de Lebesgue agit transitivement sur
l'int\'erieur du disque ; et le
th\'eo\-r\`eme d'Oxtoby-Ulam \cite{oxto1} permet de g\'en\'eraliser ce r\'esultat
\`a n'importe quelle bonne mesure sur n'importe quel disque topologique ferm\'e).
\end{itemize}\nopagebreak
\end{demo} \pagebreak[1]

Cet exemple montre une des limites du
   th\'eo\-r\`eme de la ``vari\'et\'e'' stable topologique : ici,
   contrairement au cas des points hyperboliques de la
   dynamique diff\'erentiable, ni l'ensemble stable de $S$ ni son ensemble stable
   local, ni m\^eme leurs adh\'erences,
   ne sont connexes (cette remarque r\'epond \`a une question
   de S.~Baldwin et E.~E.~Slaminka, voir \cite{bald1}).

\subsection{Difficult\'es cach\'ees}
Finissons ce panorama d'exemples en mentionnant un des pi\`eges des
germes d'hom\'eo\-morphismes de surface (ou de diff\'eo\-mor\-phismes
non hyperboliques) : le simple dessin d'un germe
hyperbolique selle (troisi\`eme dessin de la figure \ref{ifig2}) est
trompeur,
 parce qu'on pourrait croire qu'il suffit \`a d\'efinir la dynamique ;
 en r\'ealit\'e, il existe une infinit\'e (non d\'enombrable) de
classes de conjugaison  de germes d'hom\'eo\-morphismes qui sont
toutes repr\'esent\'ees par ce dessin.
 Formellement, ces hom\'eo\-morphismes
pr\'eservent tous chaque feuille du feuilletage par hyperboles dessin\'e
sur la figure \ref{ifig2}, et d\'eplacent les points dans le sens des
fl\`eches ; mais la libert\'e dans le choix de la longueur du d\'eplacement
le long de chaque feuille permet cette coexistence d'une infinit\'e de
dynamiques distinctes (bien qu'extr\^emement semblables).  Certains de ces germes sont temps
un d'un flot, d'autres non, certains ne sont m\^eme conjugu\'es \`a aucun
de leurs it\'er\'es $h^n$. Le texte \cite{begu1}  donne une id\'ee de la
diversit\'e de cette faune, en d\'ecrivant un nouvel invariant de
conjugaison des hom\'eo\-morphismes de \hbox{Brouwer}, \emph{l'ensemble oscillant}.

\clearpage
\part{Interm\`ede : du local au global}\label{par.inte}
Les r\'esultats sur les hom\'eo\-morphismes locaux
vont se d\'eduire de l'\'etude de la dynamique d'une certaine
classe d'hom\'eo\-morphismes de la sph\`ere ; cette partie fait le
lien entre le contexte local et le contexte global.
\section{Le th\'eo\-r\`eme d'extension}
\begin{theo}[d'extension]
\index{extension!th\'eo\-r\`eme d'}
\label{the.exte}
Soit $h:U \rightarrow V$ un  hom\'eo\-morphisme local~\footnote{D\'efinition \ref{def.holo}.}
 fixant uniquement $0$. Alors il existe un
hom\'eo\-morphisme $H : \R^2 \rightarrow \R^2$, fixant uniquement $0$, ayant le m\^eme
germe que $h$ en $0$. 
\end{theo}

Ce th\'eo\-r\`eme est d\'ej\`a connu~: il appara\^{\i}t essentiellement dans
un article de O.~H.~Hamilton, qui l'attribue \`a M.~H.~A.~Newman (voir~\cite{hami1}).

 Nous donnons ici une nouvelle preuve,
 reposant uniquement sur la th\'eorie \'el\'ementaire des
rev\^etements, principalement sur les th\'eor\`emes de
classification des rev\^etements et de rel\`evement des applications
(voir \cite{span1}, chapitre 2, th\'eor\`emes 13 et 5 des sections 5 et
4). Plus pr\'ecis\'ement, nous allons d\'emontrer la proposition suivante~:
\begin{prop}
\label{pro.exte}
Soit $h:U\rightarrow V$ un ho\-m\'eo\-mor\-phisme
local fixant uniquement $0$~; on suppose que $U$
est un disque topologique ferm\'e~\footnote{D\'efinition \ref{def.dito}.}, 
et on note $W$ n'importe quel ouvert connexe,
contenant le point $0$, tel que $W$ et son image $h(W)$ sont inclus
dans $U$.

Alors la restriction de $h$ \`a $W$ s'\'etend en un hom\'eo\-mor\-phisme $H$ du plan
dont $0$ est le seul point fixe. 
\end{prop}

Voici l'id\'ee de la preuve~: il est facile d'\'etendre $h$ en un
hom\'eo\-mor\-phisme du plan si l'on autorise l'apparition de nouveaux
points fixes ; toute la difficult\'e consiste ensuite \`a \'eliminer ceux-ci. 
\`{A} l'aide de la th\'eorie des rev\^etements, on montre que l'application identit\'e de
$U$ peut s'\'etendre en une application $q$ du plan $\R^2$ sur le compl\'ementaire
de l'ensemble des points fixes \`a \'eliminer, qui est un rev\^etement
(sauf au dessus du point fixe $0$, qui n'a qu'un seul ant\'ec\'edent par $q$).
 On obtient alors le r\'esultat en
relevant l'hom\'eo\-mor\-phisme $h$ par cette application
$q$.\footnote{Cette application $q$ n'appara\^{\i}t pas explicitement
dans la preuve ci-dessous, mais correspond \`a l'application $p \circ
s$.} Soulignons l'un des int\'er\^ets de la m\'ethode par
rev\^etement~: cette partie de la construction est canonique
(voir~\cite{lero10} pour une application de cette remarque).

\begin{figure}[htp]
\par
\centerline{\hbox{\input{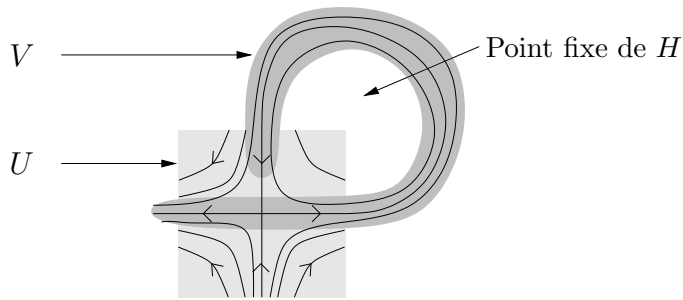}}}
\par
\caption{\label{fig49}On ne peut pas \'etendre $h$ sans rajouter de
point fixe}
\end{figure}

La figure~\ref{fig49} montre un hom\'eo\-mor\-phisme local $h : U \rightarrow
V$, fixant uniquement $0$, mais tel que tout hom\'eo\-mor\-phisme $H$ qui
\'etend $h$ poss\`ede au moins deux points fixes\footnote{On peut prouver
cette affirmation en utilisant le lemme~\ref{lem.libr} ci-dessous.}.
 Ceci montre qu'on ne
peut pas \'etendre $h$ sans restriction pr\'ealable, et explique le r\^ole
de $W$ dans l'\'enonc\'e .

\begin{demo}[de la proposition~\ref{pro.exte}]
Pour tout espace topologique $X$ localement simplement connexe, et
tout point $x_0$ de $X$, on notera $\pi_1(X,x_0)$ le groupe
fondamental de $X$ bas\'e en $x_0$. Quand cela a un sens, on note
$f_\#$ l'application induite par une application continue $f$ au
niveau des groupes fondamentaux.
On notera aussi, pour tout sous-ensemble $X$ du plan, $\dot{X}=X \setminus \{0\}$.

\paragraph{Extension avec points fixes}
Comme $U$ est un disque topologique ferm\'e, le th\'eor\`eme de
Schoenflies (voir appendice)
permet d'\'etendre $h$ en un hom\'eo\-mor\-phisme du plan, que l'on note
encore $h$ ; bien s\^ur, en g\'en\'eral, cet hom\'eo\-mor\-phisme a des points fixes hors de $U$.

Soit $O$ l'unique composante connexe de $\R^2 \setminus \fixe(h)$
contenant  $\dot{U}$ ; on a  $h(O)=O$.

\paragraph{Rev\^etement de $O$ par un anneau ouvert}
Soit $x_0$ un point de $\dot{W}$, et notons $i$  l'inclusion de
$\dot{U}$ dans $O$.
Le groupe $\pi_1(\dot{U},x_0)$ est isomorphe au groupe infini cyclique $\Z$, et
 s'identifie  \emph{via} $i_\#$ \`a un sous-groupe de $\pi_1(O,x_0)$ que
 l'on note $G$. Soit $p : (\t O,\t x_0) \rightarrow (O,x_0)$
 le rev\^etement tel que $p_\#(\pi_1(\t O,\t x_0))=G$ (th\'eo\-r\`eme~V.13
 de \cite{span1}).
L'espace topologique $\t O$ est une surface orientable, sans bord, de groupe fondamental
 isomorphe \`a $\Z$, il est par cons\'equent hom\'eomorphe \`a l'anneau ouvert
 $\A=\R/\Z \times \R$.

Soit $\t{\dot{U}}$ la composante connexe de $p^{-1}(\dot{U})$
qui contient $\t x_0$~; on voit facilement que la restriction de $p$
\`a $\t{\dot{U}}$ est un hom\'eo\-mor\-phisme  sur $\dot{U}$ (ceci vient du
fait que $p_\#(\pi_1(\t{\dot{U}}, \t x_0))=\pi_1(\dot{U}, x_0)$ et du
th\'eor\`eme de rel\`evement des applications). Dans la
compactification de l'anneau $\tilde O$ en deux bouts,  on note
$+\infty$  le bout qui est dans l'adh\'erence de  $\t{\dot{U}}$. 
On consid\`ere alors l'inverse de la restriction de $p$ \`a
$\t{\dot{U}}$~: elle s'\'etend 
 en $0$ en un hom\'eo\-morphisme $s$ entre $U$ et
$\t{\dot{U}}\cup \{+\infty\}$.

\paragraph{Relev\'e de $h$}
Puisque $x_0$ est dans $\dot W$, que $W$ et $h(W)$ sont inclus dans
$U$, et que $h$ pr\'eserve l'orientation,
on a $h_\#(G)=i_\# \pi_1(\dot{U},h(x_0))$.
\footnote{\label{foo:coin}C'est ici que la preuve ``coince'' si $h(W)$ n'est pas
inclus dans $U$~: en effet
 la premi\`ere extension de $h$ aurait alors pu  faire appara\^itre des points
fixes ``dans les trous'' entre $W$ et $h(W)$ (voir figure~\ref{fig49})~; dans ce cas, on
n'aurait pas cette \'egalit\'e.}
 Par cons\'equent, $p_\#(\pi_1(\t O, s\circ h(x_0)))=h_\#(G)=(h\circ
p)_\#(\pi_1(\t O,\t x_0))$. Le th\'eo\-r\`eme de rel\`evement des application
entra\^ine alors l'existence 
d'une (unique) application $\t h : \t O \rightarrow \t
O$ telle que $p \t h=h p$, et qui envoie $\t x_0=s(x_0)$ sur $s(h(x_0))$. 
Cette application est un hom\'eo\-mor\-phisme (car c'est un rev\^etement
 qui induit un isomorphisme au niveau des groupes fondamentaux). Cet
 hom\'eo\-morphisme se prolonge en un hom\'eo\-morphisme du plan  $\t O \cup
 \{+\infty\}$ qui fixe $+\infty$, et qui n'a pas d'autre point fixe
 (car $h$ n'a pas de point fixe autre que $0$ dans $O$).
L'application $s$ \'etablit une conjugaison entre le germe de $h$ en $0$
et celui de $\t h$ en $+\infty$.

Le th\'eor\`eme de Schoenflies permet de prolonger
$s$ en un hom\'eo\-mor\-phisme $S$ entre  $\R^2$ et $\t
O \cup \{+\infty\}$. On v\'erifie que l'hom\'eo\-mor\-phisme du plan
$H=S^{-1} \t h S$ fixe seulement $0$, et qu'il co\"\i ncide avec $h$
sur $W$ (car $\t h(s(W)) \subset s(U)$).
\end{demo}

\clearpage
\part{Dynamique globale : \'enonc\'e et r\'esultats pr\'eliminaires}
\label{par.enon}

Dans cette partie, on \'enonce un th\'eo\-r\`eme sur la dynamique de certains
hom\'eo\-morphismes de la sph\`ere (section \ref{sec.deen}), et on introduit les outils n\'ecessaire \`a
sa d\'emons\-tra\-tion : la th\'eorie des hom\'eo\-morphismes de  \hbox{Brouwer} 
(section~\ref{sec.thbr}), la notion 
d'indice partiel (section~\ref{sec.inpa}), et les d\'ecom\-po\-sitions en
briques (section~\ref{sec.deco}). \emph{Les parties \ref{par.enon}
 et \ref{par.preu} sont ind\'ependantes
des parties \ref{par.pres} et \ref{par.inte}. Dans la
partie~\ref{par.enon},  les sections~\ref{sec.inpa} et
\ref{sec.deco} sont ind\'ependantes, et peuvent \^etre lues dans un
ordre quelconque}.

\section{D\'efinitions, \'enonc\'e}
\label{sec.deen}
On note $\S^2$ la sph\`ere topologique orient\'ee de dimension $2$, $N$ et $S$ deux points
distincts de $S^2$.
\textbf{  On se donne un \homeo\ $h$ de $S^2$ qui
pr\'eserve l'orientation. On suppose que
$\fixe(h)=\{N,S\}$}, o\`u $\fixe(h)$ est l'ensemble des \textit{points fixes} de $h$,
$$\fixe(h)=\{x \in \S^2 \mid h(x)=x\}.$$
\index{$\fixe(h)$}
L'ensemble de ces donn\'ees sera appel\'e ``\textbf{hypoth\`ese (H2)}''.
\index{hypoth\`ese (H2)}
\subsection{D\'efinitions}
\begin{defi} 
\label{def.attr}
\index{attracteur}
\index{attracteur!strict}
  Un ensemble $E \subset \S^2$ est
\res{positivement invariant} si $h(E) \subset E$ ; on dit aussi que
$E$ est un \res{attracteur}.
 L'ensemble $E$ est un \res{attracteur
strict} si $h(E) \subset \inte(E) \cup \fixe(h)$. On d\'efinit aussi les
ensembles
\res{n\'egativement invariants}, \res{r\'epulseurs} et \res{r\'epulseurs stricts} en changeant $h$
en $h^{-1}$. Un ensemble $E$ est \res{totalement invariant} si il est \`a la
fois positivement et n\'egativement invariant ($h(E)=E$).
\end{defi}

\begin{defi}
\index{disque topologique}
Un \res{disque topologique ferm\'e} (respectivement \res{ouvert})
est un ensemble hom\'eo\-morphe au
disque unit\'e ferm\'e (respectivement ouvert) du plan.
\end{defi}

\begin{figure}[htp]
  \par
\centerline{\hbox{\input{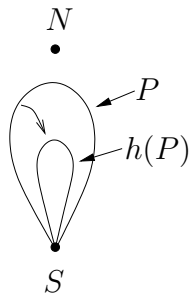}}}
\par
  \caption{\label{fig4}Un p\'etale attractif en $S$}
\end{figure}

\begin{defi}[ (figure \ref{fig4})]
\index{p\'etale}\label{def.peta}
 Un \res{p\'etale attractif bas\'e en $S$} est un disque topologique  ferm\'e
 $P$ dans $\S^2$, v\'erifiant :
\begin{itemize}
\item $S \in \partial P$ ;
\item $N \not \in P$;
\item $P$ est un attracteur strict.
\end{itemize}
 Un \res{p\'etale r\'epulsif} est un p\'etale attractif pour
$h^{-1}$. On d\'efinit de m\^eme les p\'etales attractifs et
r\'epulsifs bas\'es en $N$.
\end{defi}

\begin{figure}[htp]
  \par
\centerline{\hbox{\input{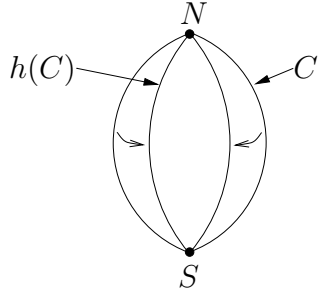}}}
\par
  \caption{\label{fig5}Un croissant attractif}
\end{figure}

\pagebreak
\begin{defi}[ (figure \ref{fig5})]
\index{croissant}
  Un \res{croissant attractif} pour $h$ est  un disque topologique
 ferm\'e $C$ dans $\S^{2}$, v\'erifiant~: \nopagebreak
\begin{itemize}
\item $N, S \in \partial C$ ;
\item $C$ est un attracteur strict.
\end{itemize}
 Un \res{croissant r\'epulsif} est un croissant attractif pour $h^{-1}$.
\end{defi}

\begin{figure}[htp]
  \par
\centerline{\hbox{\input{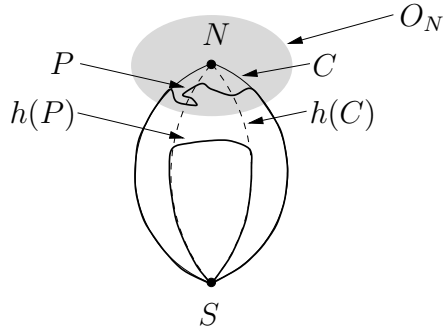}}}
\par
  \caption{\label{fig6}Un croissant attractif \`a dynamique Nord-Sud}
\end{figure}

Certains croissants ont la propri\'et\'e d'\^etre des limites de p\'etales~:
\begin{defi}[ (figure \ref{fig6})]\label{def.cans}
\index{croissant!\`a dynamique Nord-Sud}
 Le croissant attractif $C$ est dit \res{\`a dynamique Nord-Sud}
   si de plus pour tout voisinage $O_N$ de $N$, il existe un p\'etale
   attractif $P$ bas\'e en $S$ tel que 
\begin{itemize}
\item $P \subset C$ ;
\item $C \setminus P \subset O_N$.
\end{itemize} 

Un \res{croissant r\'epulsif \`a dynamique Sud-Nord} est un
 croissant attractif \`a dynamique Nord-Sud pour $h^{-1}$.
 On d\'efinit de m\^eme les croissants attractifs \`a dynamique
 Sud-Nord et les croissants r\'epulsifs \`a dynamique Nord-Sud
 en inversant les r\^oles des points  $N$ et $S$.
\end{defi}

\subsection{\'Enonc\'e du th\'eor\`eme principal}
\begin{figure}[hbtp]
\par
\centerline{\hbox{\input{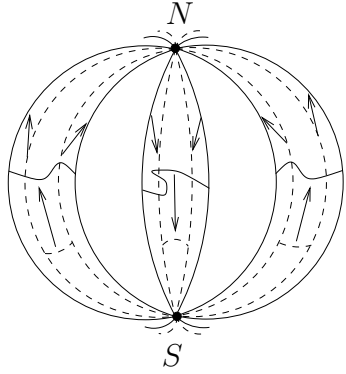}}}
\par
\caption{\label{fig14}Illustration du th\'eor\`eme principal}
\end{figure}

On peut maintenant \'enoncer le r\'esultat qui constitue le c\oe ur de
l'article (figure \ref{fig14}) :
\begin{theo}
\label{the.prin} 
 Soit $h$ un \homeo\ de la sph\`ere, pr\'eservant l'orientation,
fixant uniquement les deux points $N$ et $S$ (\textbf{hypoth\`ese
(H2)}), et tel que $\indi(N)=1-p <1$. 

Il existe alors $p$ croissants attractifs \`a dynamique Nord-Sud,
et $p$ croissants r\'epulsifs \`a dynamique Sud-Nord, deux \`a
deux d'intersection r\'eduite \`a $\{N,S\}$,
 les croissants attractifs et r\'epulsifs \'etant
cycliquement altern\'es autour de $N$ et $S$.
\end{theo}

 D'apr\`es le th\'eo\-r\`eme de Schoenflies-Homma (voir l'appendice),
\index{Schoenflies-Homma!th\'eo\-r\`eme de} on peut supposer que les
bords des croissants et leurs images sont des grands cercles
 (g\'eod\'esiques) de la sph\`ere euclidienne~; le
th\'eo\-r\`eme est alors illustr\'e par la figure~\ref{fig14}.
Remarquons que les hypoth\`eses du th\'eo\-r\`eme ne sont pas sym\'etriques en $N$ et $S$~: le point $N$
est d'indice $1-p<1$, tandis que l'indice de $S$ est $1+p>1$ (d'apr\`es
la formule de Lefschetz)~; la conclusion est tout aussi dissym\'etrique~: notamment, tous les p\'etales
associ\'es aux croissants sont bas\'es en $S$.

\subsection{Id\'ee de la preuve du th\'eo\-r\`eme \ref{the.prin}}\label{sub.idee}
Expliquons bri\`evement les grandes \'etapes de la preuve du th\'eo\-r\`eme.
On commence par construire une \emph{d\'ecom\-po\-sition en briques}~: il
s'agit d'une sorte de triangulation de $\S^2 \setminus\{N,S\}$,
suffisamment fine pour que chaque triangle (appel\'e ici
\emph{brique}) soit disjoint de son image par $h$ (section~\ref{sec.deco}).
 On cherche alors
des croissants et p\'etales \emph{simpliciaux}, c'est-\`a-dire qui sont
des r\'eunions de briques de la d\'ecom\-po\-sition. 
Une courte  \'etude
combinatoire permet tout d'abord de trouver un premier croissant simplicial
(en fait, on se contente d'en trouver un c\^ot\'e, proposition~\ref{pro.dbns} et section~\ref{sec.dbns}).
Par ailleurs,  on montre que tout croissant simplicial minimal pour l'inclusion
est \`a dynamique Nord-Sud ou Sud-Nord (proposition~\ref{pro.cam} et
section~\ref{sec.cam}).

On consid\`ere alors une famille maximale $\cal F$  de croissants attractifs et
r\'epulsifs simpliciaux minimaux deux \`a deux d'int\'erieurs disjoints
(section~\ref{sec.preu}). Il reste \`a  montrer que cette famille $\cal F$
contient suffisamment de croissants. Pour cela, on calcule l'indice du
point fixe $N$ \`a l'aide de la notion \emph{d'indice partiel}
(\'etudi\'ee \`a la section~\ref{sec.inpa}). Dans ce calcul
(proposition~\ref{pro.codb2} et section~\ref{sec.brsu}),
 on montre d'une part que les zones entre deux
croissants successifs de la famille $\cal F$ ont une contribution
nulle (c'est une cons\'equence de la maximalit\'e de $\cal F$).  D'autre
part, la contribution d'un croissant vaut $+1/2$ ou $-1/2$ selon qu'il
est attractif ou r\'epulsif, et \`a dynamique Nord-Sud ou Sud-Nord~: les
``bons'' croissants (ceux des deux types recherch\'es par le th\'eo\-r\`eme~\ref{the.prin})
ayant une contribution n\'egative et les ``mauvais'' une contribution
positive. Un petit argument combinatoire montre alors que la famille
de croissants $\cal F$ contient une sous-famille satisfaisant le
r\'esultat du th\'eo\-r\`eme.

Il reste \`a pr\'eciser que les r\'esultats interm\'ediaires sont racont\'es dans un
cadre diff\'erent, obtenu en relevant la dynamique au rev\^etement
universel de $\S^2 \setminus \{N,S\}$. L'existence d'une unique
mani\`ere non triviale\footnote{C'est-\`a-dire non conjugu\'ee \`a une
translation, voir la proposition~\ref{pro.comp}.} de relever la dynamique est prouv\'ee \`a la
section~\ref{sec.reca}. On peut r\'esumer l'int\'er\^et de ce passage au
rev\^etement par le fait suivant~: dans la sph\`ere, le compl\'ementaire
d'un croissant est connexe~; par contre, le relev\'e d'un croissant au
rev\^etement universel est une bande qui s\'epare le plan en deux
composantes connexes, ce qui permet de parler du ``c\^ot\'e droit'' et
du ``c\^ot\'e gauche'' de la bande.

\section{La th\'eorie de \hbox{Brouwer}}
\label{sec.thbr}
Dans cette section, on rappelle des objets et des r\'esultats autour de la th\'eorie
de \hbox{Brouwer} : l'indice, les cha\^ines de disques et une variante d'un
lemme de Franks qui joue un r\^ole important dans ce texte, les
droites de \hbox{Brouwer} et l'\'enonc\'e du th\'eo\-r\`eme de translation plane.

La th\'eorie de \hbox{Brouwer} traite des hom\'eo\-mor\-phismes de la sph\`ere $\S^2$  qui n'admettent pas
de courbe d'indice $1$. Le point cl\'e est que cette hypoth\`ese interdit
presque toute forme de r\'ecurrence non tri\-via\-le  dans la dynamique (section
\ref{sub.recu}). Cette propri\'et\'e remarquable
permet ensuite  de montrer que tous ces hom\'eo\-mor\-phismes peuvent \^etre
obtenus en recollant des translations : plus pr\'ecis\'ement, 
 le compl\'ementaire des points fixes est recouvert
par des ouverts totalement invariants, hom\'eomorphes au plan (autrement
dit : connexes et simplement connexes), sur lesquels la dynamique est
conjugu\'ee \`a une translation du plan ; c'est une partie de ce qu'on appelle le
th\'eor\`eme de translation plane (section~\ref{sub.entp}). Les r\'ef\'erences modernes sont
\cite{barg1,brow1,brow4,brow9,fath1,fran1,guil1,guil2,leca3}.

La th\'eorie construit notamment des courbes disjointes de leur image,
appel\'ees \emph{droites de \hbox{Brouwer}} (section~\ref{sub.drbr})~;
 de telles courbes formeront les bords des
croissants et des p\'etales dans la preuve du th\'eor\`eme principal. 

\textbf{On suppose d\'esormais que $h$ est un hom\'eo\-mor\-phisme de la
sph\`ere $\S^2$,
pr\'eservant l'orien\-ta\-tion, n'ayant qu'un nombre fini de points fixes}
($h$ poss\`ede au moins un point fixe d'apr\`es la formule de Lefschetz,
voir ci-dessous). \`{A} partir de 
la section \ref{sub.drbr}, on supposera de plus que $h$ n'a pas de courbe
d'indice $1$.

\subsection{Rappels sur l'indice}
\index{indice}\label{ss.rain}
 La formule de
Lefschetz est expliqu\'ee dans les livres \cite{gupo1} dans le cadre
diff\'e\-ren\-tiable, \cite{browrf1,dold2} dans le cadre topologique.

\subsubsection*{$\bullet$ Indice d'une courbe dans le plan}

\begin{defi}\index{courbe}\index{cercle topologique} \index{arc}
 Une \res{courbe} du plan est une application continue $\gamma$ de $[0,1]$ dans
 $\R^2$. Elle est \res{simple} si c'est une application injective (on
 dit aussi que c'est un \res{arc}) ;
 \res{ferm\'ee} si $\gamma(0)=\gamma(1)$ ; \res{ferm\'ee simple} si
 elle n'a pas d'autre point double que $\gamma(0)$ (dans ce cas, on
 parle aussi de \res{courbe de Jordan}, ou de \res{cercle
 topologique}). Une courbe \res{va de $X$ \`a $Y$} ou \res{joint $X$ et
 $Y$} si $\gamma(0) \in X$ et $\gamma(1) \in Y$.
Les \res{extr\'emit\'es} de $\gamma$ sont les points $\gamma(0)$ et
$\gamma(1)$~; son \res{int\'erieur} est $\gamma(]0,1[)$, not\'e $\inte(\gamma)$.

On utilisera le m\^eme vocabulaire pour les courbes de
 la sph\`ere. On confondra souvent une courbe et son image.
\end{defi}
Soit $\gamma$ une courbe du plan euclidien $\R^2$.
Si $\vec v$ est un champ de vecteurs continu d\'efini sur $\gamma$ et ne s'y annulant pas,
 on appelle indice 
 de $\vec v$ le long de $\gamma$, et on note $\indi(\vec v,\gamma)$,
la variation angulaire de  $\vec v$ quand on parcourt $\gamma$ :
plus pr\'ecis\'ement, soit $\phi$ le rev\^etement universel du cercle
unit\'e du plan :
 $$
\begin{array}{rcl}
\phi : \R & \longrightarrow & \S^1 \\
\theta & \longmapsto & \exp(2i\pi\theta). 
\end{array}
$$
L'application 
$$
f = \frac{\vec v \circ \gamma}{\|\vec v \circ \gamma\|} : [0,1]  \longrightarrow  \S^1 
$$
se rel\`eve par $\phi$  en une application $F : [0,1] \rightarrow \R$
(\ie telle que $\phi \circ F=f$), et on pose
$\indi(\vec v,\gamma)=F(1)-F(0)$.
Remarquons que si $\gamma$ est une courbe ferm\'ee, l'indice est un
nombre entier.

Soit $g$ une application continue d'un ouvert $U$ du plan dans le
plan, et $E$ l'ensemble des points fixes de
$g$ ; on  d\'efinit sur $U \setminus E$ le champ de vecteurs
$$
\vec g_x=\frac{g(x)-x}{\| g(x)-x\|}.
$$
Si maintenant $\gamma$ est une courbe dans $U \setminus E$,
on pose $\indi(g,\gamma)=\indi(\vec g, \gamma)$.

Si $\gamma$ est une courbe ferm\'ee, l'entier  $\indi(g,\gamma)$ ne
d\'epend que de la classe d'homotopie de $\gamma$ (en tant que courbe
ferm\'ee) dans le compl\'ementaire
des points fixes de $g$ : en particulier, une courbe qui y est
nulhomotope est d'indice nul. 
 L'indice est un invariant de conjugaison,
c'est-\`a-dire que pour tout hom\'eo\-mor\-phisme $\phi$ pr\'eservant
l'orientation, on a 
$$
\indi(\phi \circ g \circ \phi^{-1},\phi(\gamma))=\indi(g,\gamma).
$$
\index{espace d'hom\'eo\-mor\-phismes}
En effet, l'espace des hom\'eo\-mor\-phismes du plan pr\'eservant
l'orientation est connexe, et le nombre $\indi(\phi \circ g \circ
\phi^{-1},\phi(\gamma))$ est un entier qui d\'epend contin\^ument de
$\phi$, il est donc constant (voir par exemple \cite{lero1}
pour les propri\'et\'es topologiques des espaces d'hom\'eo\-mor\-phismes de surfaces).

\subsubsection*{$\bullet$ Indice d'un point fixe}
Soit $g$ comme ci-dessus, et $x_0$ un point fixe isol\'e de $g$,
 c'est-\`a-dire un point isol\'e de l'ensemble $E$.
 On d\'efinit alors l'indice de $x_0$
comme l'indice de n'importe quelle courbe de Jordan
 $\gamma$ dans $U$  qui entoure $x_0$ mais n'entoure  aucun autre point fixe.
Ce nombre ne d\'epend pas de la courbe $\gamma$ choisie, et ne d\'epend
 que de la classe de conjugaison du  germe de $g$ en $x_0$.

On a alors la formule suivante : l'indice d'une courbe de Jordan $\gamma$
qui n'entoure qu'un nombre fini de
points fixes de $g$ est \'egal \`a la somme des indices de ces
points fixes.

\subsubsection*{$\bullet$ Indices sur la sph\`ere}
Revenons \`a notre hom\'eo\-morphisme $h$ de la sph\`ere.
On peut maintenant d\'efinir l'indice d'un point fixe de $h$ en se
ramenant au plan au moyen d'une carte de la sph\`ere.

Comme $h$ est isotope \`a l'identit\'e, et comme la caract\'eristique d'Euler-Poincar\'e de la
sph\`ere vaut $2$, on a la formule de
Lefschetz (voir \cite{})~: 

\centerline{\emph{La somme des indices des points fixes
de $h$ vaut $2$.}}\index{Lefschetz!formule de}
\noindent  Notamment, $h$ a au moins un point fixe.

Soit $\gamma$ une courbe de Jordan de la sph\`ere \'evitant les points fixes de
$h$. Choisissons un point fixe
 de $h$, que l'on va noter $\infty$,
 et  identifions l'ouvert $\S^2 \setminus\{\infty\}$ au
plan. On peut montrer que l'indice de $\gamma$ dans cette carte, que
l'on note $\indi_\infty(h,\gamma)$, ne
d\'epend que de la position de $\infty$ par rapport \`a $\gamma$ : plus
pr\'ecis\'ement, le compl\'ementaire de $\gamma$ dans $\S^2$ a exactement deux
composantes connexes (th\'eor\`eme de Jordan) ; si $\infty$ et $\infty'$
sont dans la m\^eme composante connexe, on a
$\indi_{\infty'}(h,\gamma)=\indi_{\infty}(h,\gamma)$ ;
 dans le cas contraire (\ie si $\gamma$ s\'epare $\infty$ et
$\infty'$), $\indi_{\infty'}(h,\gamma)=2-\indi_{\infty}(h,\gamma)$.
Cette affirmation est en fait  \'evidente si $h$ n'a qu'un nombre fini
de points fixes (avec la formule de Lefschetz),
 puisque l'indice d'une courbe plane est la somme des
indices des points fixes qu'elle entoure.

Il ressort de ce qui pr\'ec\`ede que :
\begin{affi}
\label{aff.indi}
Soit $\gamma$ une courbe de Jordan sur la sph\`ere ;
 les propri\'et\'es suivantes sont \'equi\-va\-lentes :
\begin{enumerate}
\item il existe un point fixe $\infty$ de $h$ tel que
$\indi_\infty(h,\gamma)=1$ ;
\item pour tout point fixe $\infty$ de $h$, $\indi_\infty(h,\gamma)=1$
;
\item la somme des indices des points fixes de $h$ dans chacune des
deux composantes du compl\'ementaire de $\gamma$ vaut $1$. 
\end{enumerate}
\end{affi}
L'affirmation pr\'ec\'edente permet de dire ce qu'est une courbe de Jordan
d'indice $1$ dans la sph\`ere~:
\begin{defi}
\index{courbe!d'indice $1$}
Une courbe $\gamma$ dans $\S^2$ est dite \res{d'indice $1$} si les
propri\'et\'es \'equivalentes de l'affirmation \ref{aff.indi} sont v\'erifi\'ees.
\end{defi}

\subsubsection*{$\bullet$ Indice des bords des disques attractifs}

Finissons par un lemme qui nous servira \`a plusieurs reprises, et dont
la preuve est typique des calculs d'indices :
\begin{lemm}
\label{lem.diat}
\index{disque attractif}
Soit $D$  un disque topologique ferm\'e de $\S^2$ dont la fronti\`ere ne contient pas
de point fixe. Si $D$ est un attracteur ou un r\'epulseur, alors le bord de
$D$ est une courbe d'indice $1$. 
\end{lemm}
\begin{demo}
Quitte \`a changer $D$ en son compl\'ementaire, on suppose que $D$ est un
attracteur. Le disque compl\'ementaire de $D$ est un attracteur pour
$h^{-1}$, il contient donc un point fixe de $h$
(d'apr\`es le th\'eor\`eme de point fixe de \hbox{Brouwer}). 
Le compl\'ementaire de ce point fixe est identifi\'e au plan.
En utilisant le th\'eor\`eme de Schoenflies, 
on se ram\`ene alors  \`a la situation
o\`u $D$ est le disque unit\'e du plan. 

Pour tout $t \in [0,1]$, on note
$\phi_t$ l'homoth\'etie de centre $(0,0)$ et de rapport $t$, et on pose
$h_t(x)=\phi_t\circ h(x)$. L'hypoth\`ese $h(D) \subset D$ entra\^\i ne que
$h_t$ n'a pas de point fixe sur $\partial D$. Le nombre
$\indi(h_t,\partial D)$ est donc d\'efini pour tout $t$, c'est un entier
qui varie contin\^ument, il est donc constant. D'autre part, $h_1=h$,
et $h_0$ est l'application constante $x \mapsto (0,0)$,
 d'o\`u $\indi(h,\partial D)= \indi(h_0, \partial D)=1$.
\end{demo}

\subsection{R\'ecurrence entra\^\i ne indice 1}
\label{sub.recu}
\begin{figure}[hbtp]
\par
\centerline{\hbox{\input{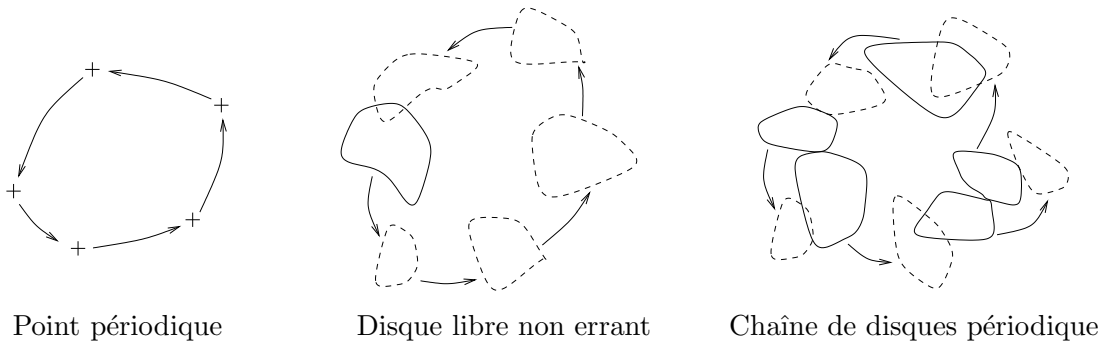}}}
\par
\caption{\label{fig9}Trois types de r\'ecurrence de plus en plus fins}
\end{figure}

Les trois lemmes de ce paragraphe annoncent l'existence d'une
courbe d'indice 1 sous des hypoth\`eses de ``quasi-p\'eriodicit\'e'' de 
plus en plus faibles, illustr\'ees par la figure~\ref{fig9}.

\subsubsection*{$\bullet$ Orbite p\'eriodique}
\index{orbite p\'eriodique}
\begin{lemm}\label{lem.peri}
Si $h$ a un point p\'eriodique non fixe, alors il existe une courbe
d'indice $1$.
\end{lemm}

Nous admettons ce lemme : pour la preuve en toute g\'en\'eralit\'e,
 voir \cite{guil1}, \cite{fath1}, \cite{brow1} ou  \cite{barg1}. 

Si la p\'eriode est $2$, voici un argument adapt\'e de la preuve de
 A. Fathi revisit\'ee par 
 M. Barge et J. Franks.
On consid\`ere le rev\^etement $\pi :
 \S^2 \rightarrow \S^2$ \`a deux feuillets
 ramifi\'e au dessus des deux points de l'orbite p\'eriodique. Soit $\tilde h$ un relev\'e
 de $h$ par $\pi$, et
 $F$ le sous-ensemble des points fixes de $h$ constitu\'e de ceux qui se rel\`event en deux
points fixes de $\tilde h$. Alors l'ensemble des points fixes de
$\tilde h$ est $\pi^{-1}(F)$ ; la somme des indices des points de
$\pi^{-1}(F)$ pour $\tilde h$ est donc \'egale \`a $2$ (formule de
 Lefschetz) ; mais c'est aussi
le double de la somme des indices des points de $F$ pour $h$. D'apr\`es
l'affirmation \ref{aff.indi}, toute courbe s\'eparant $F$ des autres points
fixes de $h$ est d'indice $1$.

 Pour une p\'eriode plus grande, on ne connaît pas de preuve de
``topologie alg\'ebrique''. L'id\'ee de \cite{barg1} consiste \`a  se ramener au cas d'une p\'eriode $2$
en effectuant une s\'erie de modifications libres de $h$ (d\'efinition \ref{def.moli}). 
Nous allons d\'eduire les deux lemmes suivants de ce lemme admis.

\subsubsection*{$\bullet$ Disque libre non errant}

\begin{defi}\index{libre}
Un ensemble connexe $C$ est \res{libre} s'il est disjoint de son image $h(C)$.
\end{defi}

\begin{lemm}[figure \ref{fig9}, milieu]
\label{lem.libr}
Soit $\ii{D}$ un ensemble connexe par arcs de la sph\`ere. On suppose que $\ii{D}$ est libre,
mais qu'il  rencontre l'un de ses it\'er\'es $h^n(\ii{D})$ pour $n
\neq 0$. Alors il existe une courbe d'indice $1$.
\end{lemm}

La preuve utilise la notion de modification libre (d\'efinition
\ref{def.moli}).\index{modification libre}

\begin{demo}[du lemme \ref{lem.libr}]
\emph{Dans un premier temps, on suppose que $D$ est un disque
topologique ouvert.}
On peut d'abord supposer que $n$ est positif (quitte \`a remplacer $D$ par
 $D'=h^n(D)$ et $n$ par $n'=-n$),
 puis que c'est le plus petit
entier positif v\'erifiant $D \cap h^n(D)
\neq \emptyset$.
Par hypoth\`ese, il existe un point $x$ de $D$ tel
que $h^n(x) \in D$. Soit $\phi$ un hom\'eo\-mor\-phisme qui envoie $h^n(x)$
sur $x$ et qui est l'identit\'e hors de $D$ ; on pose $h_1=\phi \circ
h$ (voir figure~\ref{fig22}).
\begin{figure}[htpb]
\par
\centerline{\hbox{\input{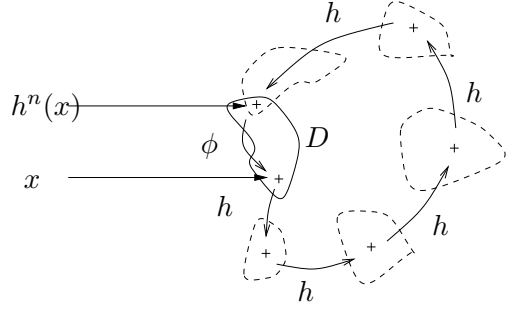}}}
\par
\caption{\label{fig22}Transformation d'une orbite ``quasi-p\'eriodique''
en orbite p\'eriodique}
\end{figure}

En utilisant la minimalit\'e de $n$, on voit qu'en restriction \`a $D$,
on a $h_1^n=\phi \circ h^n$. On en d\'eduit que
 $x$ est un point p\'eriodique de $h_1$. D'apr\`es le lemme \ref{lem.peri},
il existe une courbe $\gamma$ d'indice 1 pour $h_1$. D'autre part, les
modifications libres ne changent pas l'indice des courbes, donc
$\gamma$ est aussi une courbe d'indice $1$ pour $h$.

\emph{On traite maintenant le cas g\'en\'eral, o\`u $D$ est n'importe quel
ensemble connexe par arcs.}
Soit $x$ un point de $D$ tel que $h^n(x)$ soit aussi dans $D$. Soit
$\gamma$ un arc inclus dans $D$ allant de $x$
\`a $h^n(x)$. Cet arc est clairement libre~; en utilisant le th\'eo\-r\`eme de
Schoenflies, on trouve un disque topologique ouvert $D'$, contenant
$\gamma$, qui est encore libre~; or $D'$ rencontre son it\'er\'e $h^n(D')$.
 On conclut alors en appliquant le premier cas.
\end{demo}

Rappelons
 que \textit{l'ensemble $\omega$-limite} d'un point $x$ est
$$\omega(x)=\bigcap_{n_0 \geq 0}\adhe(\{h^n(x) \mid n \geq n_0\}).$$

\begin{coro}\label{cor.erre}
Soit $x \in \S^2$ ; si il n'y a pas de courbe d'indice $1$,
alors $\omega(x)$ est r\'eduit \`a un point fixe de $h$ : autrement dit,
la suite $(h^n(x))_{n \geq 0}$ converge vers un point fixe de $h$.
\end{coro}

\begin{demo}
Soit $x_0$ un point de $\S^2$ qui n'est pas fixe. Alors tout disque
$D$ centr\'e en $x_0$ et assez petit est libre. D'apr\`es le lemme
\ref{lem.libr}, $D$ est disjoint de tous ses it\'er\'es (on dit que $x_0$
est un point errant), autrement dit aucune orbite ne peut rencontrer
$D$ plus d'une fois ; donc $x_0$ n'est dans l'ensemble $\omega$-limite
d'aucun point $x$.

Ceci montre que pour tout $x$, $\omega(x) \subset \fixe(h)$ ; il reste
\`a voir que $\omega(x)$ ne peut pas contenir deux points fixes distincts.
Soit $k$ le nombre de points fixes de $h$, soit
$V_1, \cdots, V_k$ des petits voisinages ouverts  (disjoints) de chacun des points fixes,
et $K$ le compact compl\'ementaire de la r\'eunion des $V_i$.
On peut choisir les $V_i$ assez petits pour que la propri\'et\'e de
transition suivante soit respect\'ee : tout point $x$ dans un
$V_i$ qui a un it\'er\'e futur $h^n(x)$ dans un autre $V_j$ doit avoir un
it\'er\'e interm\'ediaire $h^{m}(x)$  ($0<m<n$) dans $K$. D'autre part, $K$
peut \^etre recouvert par un nombre fini de disques libres, et
 chaque disque libre ne peut contenir qu'un seul it\'er\'e de $x$ (lemme
\ref{lem.libr}). Le point $x$ n'a donc qu'un nombre fini d'it\'er\'es
positifs dans $K$ ; d'apr\`es la propri\'et\'e de transition, 
tous ses it\'er\'es assez grands appartiennent alors au
m\^eme $V_i$, et $\omega(x)$ contient un point fixe.
\end{demo}

\subsubsection*{$\bullet$ Cha\^\i ne de disques p\'eriodique}
Le concept de cha\^{\i}ne de disques a \'et\'e introduit par J. Franks
pour g\'en\'eraliser le th\'eor\`eme de Poincar\'e-Birkhoff (\cite{fran4}).
Il permet d'affaiblir encore l'hypoth\`ese de quasi-p\'eriodicit\'e du
lemme~\ref{lem.libr}. 
Le lemme de J. Franks concerne les cha\^{\i}ne de
\emph{disques ouverts} (d\'efinition~\ref{def.cddo} ci-dessous). Nous
donnons ici une g\'en\'eralisation aux cha\^{\i}nes de \emph{pseudo-disques}
(d\'efinitions~\ref{def.psdi} et \ref{def.cdpd}).
Cette petite am\'elioration
technique nous permettra de simplifier les preuves utilisant les
d\'ecom\-po\-si\-tions en briques, en rendant
inutile la propri\'et\'e de transversalit\'e (voir  \cite{leca3},
\cite{sauz1} et la section~\ref{sec.deco} du pr\'esent texte).

\begin{defi}\index{pseudo-disque}\label{def.psdi}
Un \res{pseudo-disque} est une partie $D$ de la sph\`ere telle que pour tout $x,y \in D$, il existe
un arc $\gamma$ qui va de $x$ \`a $y$, tel que l'int\'erieur de $\gamma$
est inclus dans l'int\'erieur de $D$.
\end{defi}

\begin{defi}[ (figure \ref{fig9}, droite)]
\label{def.cddo}
\index{cha\^{\i}ne!de disques ouverts}
\index{cha\^{\i}ne!p\'eriodique}
  Une \res{cha\^\i ne de disques ouverts} pour $h$ est une suite
  $(\ii{B_i})_{i=1..k}$, $(k \geq 1)$, de disques topologiques ouverts
de la sph\`ere tels que~:
\begin{enumerate}
  \item les disques $B_i$ sont disjoints deux \`a deux~;
  \item chaque disque $B_i$ est libre~; 
  \item il existe des entiers $n_i >0$ tels
  que $h^{n_i}(\ii{B_i}) \cap \ii{B_{i+1}} \neq \emptyset$ pour $i=1, \cdots,
  k-1$.
\end{enumerate}
La cha\^\i ne de disques ouverts est dite \res{p\'eriodique}
 si de plus il existe $n_k>0$ tel
que $h^{n_k}(\ii{B_k}) \cap \ii{B_{1}} \neq \emptyset$.
\end{defi}

\begin{defi}
\label{def.cdpd}
\index{cha\^{\i}ne!de pseudo-disques}
  Une \res{cha\^\i ne de pseudo-disques} pour $h$ est une suite
  $(\ii{B_i})_{i=1..k}$, $(k \geq 1)$, de pseudo-disques
de la sph\`ere tels que~:
\begin{itemize}
  \item[1'.] les \textbf{int\'erieurs} des pseudo-disques $B_i$ sont disjoints deux \`a deux~;
  \item[2'.] chaque pseudo-disque $B_i$ est libre~; 
  \item[3'.] il existe des entiers $n_i >0$ tels
  que $h^{n_i}(\ii{B_i}) \cap \ii{B_{i+1}} \neq \emptyset$ pour $i=1, \cdots,
  k-1$.
\end{itemize}
La cha\^\i ne de pseudo-disques est dite \res{p\'eriodique}
 si de plus il existe $n_k>0$ tel
que $h^{n_k}(\ii{B_k}) \cap \ii{B_{1}} \neq \emptyset$.
\end{defi}

\begin{defi}\index{transition!temps de}\index{transition!points de}
Soit $(\ii{B_i})_{i=1..k}$ une cha\^ine de pseudo-disques. Pour chaque
entier $i$  entre $1$ et $k-1$, soit $n_i$ \textbf{le plus petit} entier
positif v\'erifiant  $h^{n_i}(\ii{B_i}) \cap \ii{B_{i+1}} \neq \emptyset$.
Les entiers $(n_i)_{i=1..k-1}$ seront dits \emph{temps de transition} de
la cha\^{\i}ne de pseudo-disques, et on appellera \emph{points de transition}
des points $(x_i)_{i=1..k-1}$ avec $x_i \in B_i \cap
h^{-n_i}(B_{i+1})$.
Pour une cha\^{\i}ne p\'eriodique, on d\'efinira de m\^eme un temps de transition 
 $n_k$ et un point de transition $x_k \in B_k \cap h^{-n_k}(B_1)$.
\end{defi}

\begin{lemm}[``lemme de Franks'']
\label{lem.fran}
\index{lemme de Franks}
       Si il existe une  cha\^\i ne  p\'eriodique de pseudo-disques pour $h$,
       alors il existe une courbe d'indice $1$.     
\end{lemm}

\begin{demo}[du lemme \ref{lem.fran}]
L'id\'ee consiste \`a montrer que l'existence d'une
cha\^\i ne p\'eriodique de pseudo-disques entra\^\i ne l'existence
d'une cha\^\i ne  p\'eriodique de disques ouverts, puis
l'existence d'une orbite p\'eriodique (apr\`es une s\'erie de
modifications libres, en g\'en\'eralisant la preuve du
lemme~\ref{lem.libr}).\index{modification libre}

 Soit
$(\ii{B_i})_{i=1..k}$ une cha\^\i ne p\'eriodique de pseudo-disques,  telle
que $k$ soit minimal
parmi toutes les cha\^\i nes p\'eriodiques  pour $h$.
 Soient $(n_i)$ les temps de transition de cette cha\^{\i}ne, et $(x_i)$ des
points de transition.

\emph{Supposons qu'il n'existe pas de courbe d'indice $1$.}
Montrons d'abord que $k \neq 1$~: en effet, sinon, $B_1$ serait un
ensemble libre, connexe par arcs,
qui rencontre l'un de ses it\'er\'es par $h$~; ceci contredirait le
lemme~\ref{lem.libr}.

La minimalit\'e de $k$ implique alors\footnote{Dans
tout ce paragraphe, les indices $i$ et $j$ sont consid\'er\'es
modulo $k$.}~:
$$
(*) \ \ \ \forall i,j, \forall n >0, \ \ \ ( \ h^{n}(\ii{B_i}) \cap \ii{B_j}
\neq \emptyset \Longrightarrow j=i+1 \ ).
$$
 On en d\'eduit que les points $x_1,h^{n_1}(x_1), x_2, h^{n_2}(x_2),
 \dots , x_k,h^{n_k}(x_k)$ sont deux \`a deux distincts.
En effet~:
 \begin{enumerate}
\item $x_i \neq h^{n_i}(x_i)$ d'apr\`es le lemme \ref{lem.peri}~;
\item si $i \neq j$, alors $x_i \neq x_j$ (sinon
$h^{n_i}(B_j)$ rencontre $B_{i+1}$, ce qui contredit $(*)$)~;
\item  si $i \neq j$, alors $x_i \neq h^{n_j}(x_j)$  (sinon
$h^{n_i+n_j}(B_j)$ rencontre $B_{i+1}$, ce qui contredit encore $(*)$).
 \end{enumerate}

On utilise maintenant la d\'efinition des pseudo-disques~: pour chaque
$i$, on choisit un arc $\gamma_i$ d'int\'erieur
inclus dans l'int\'erieur de $\ii{B_i}$ et reliant $x_i$ et
$y_i=h^{n_{i-1}}(x_{i-1})$. 
D'apr\`es le point~2' de la d\'efinition d'une cha\^{\i}ne de pseudo-disques, 
 chaque arc $\gamma_i$ est libre. D'apr\`es le point 1', les int\'erieurs
de ces arcs sont disjoints deux \`a deux. D'autre part, les consid\'erations
pr\'ec\'edentes sur leurs extr\'emit\'es montrent que les arcs $\gamma_i$
sont \'egalement disjoints deux \`a deux. 

En \'epaississant
l\'eg\`erement ces arcs (\`a l'aide du th\'eo\-r\`eme de Schoenflies), 
on obtient des disques topologiques ouverts $B'_i$ (chaque $B'_i$ contient
$\gamma_i$),
qui sont encore libres et disjoints deux \`a deux. La suite
$(B'_i)_{i=1, \dots, k}$ est
alors une cha\^\i ne p\'eriodique de disques ouverts. 

Quitte \`a changer de cha\^{\i}ne, on peut supposer que $k$ est encore
minimal et que les $(n_i)$ sont les temps de transition de $(B'_i)_{i=1, \dots, k}$.
Dans ce cas, pour chaque $i$, les it\'er\'es $h(x_i), \dots,
h^{n_i-1}(x_i)$ sont en dehors de la r\'eunion de tous les disques
de la cha\^{\i}ne.  
Pour chaque $i$, on consid\`ere alors un hom\'eo\-morphisme $\phi_i$, \`a support
dans $\adhe(B'_i)$, qui envoie $h^{n_{i-1}}(x_{i-1})$ sur
$x_i$. Remarquons que ces $k$
hom\'eo\-morphismes commutent (voir figure~\ref{fig23}). On pose 
 $\Phi=\phi_1 \circ \cdots \circ \phi_k$, et $h_1=\Phi \circ h$. Le
point $x_1$ est un point p\'eriodique de $h_1$ (sa p\'eriode est la somme
des temps de transition).
\begin{figure}[htpb]
\par
\centerline{\hbox{\input{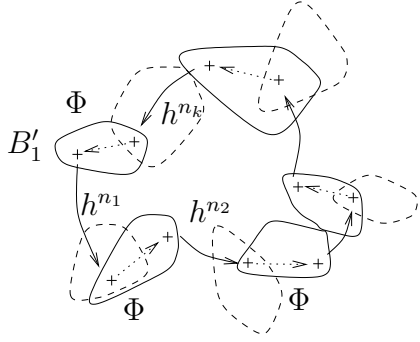}}}
\par
\caption{\label{fig23}Transformation d'une ``quasi-orbite p\'eriodique''
en orbite p\'eriodique}
\end{figure}
\end{demo}

\subsection{Droites et domaines de \hbox{Brouwer}}
\label{sub.drbr}
\textbf{On suppose dor\'enavant qu'il n'existe pas de courbe d'indice $1$ dans
$\S^2\setminus \fixe(h)$}, ce qui revient \`a dire qu'il n'existe pas de
 sous-ensemble de points fixes dont la somme des indices fait 1
(affirmation \ref{aff.indi}).

\subsubsection*{$\bullet$ D\'efinitions}
\begin{defi}
\index{plongement}
\index{droite topologique}
Un \res{plongement} est une application continue injective. Un
plongement est \res{propre} si l'image r\'eciproque de tout compact est
compacte. Une \res{droite topologique} est l'image d'un plongement propre $\phi$ de
la droite r\'eelle $\R$ dans $\S^2 \setminus \fixe(h)$.
\end{defi}

Pour un plongement de $\R$ dans $\S^2 \setminus \fixe(h)$,
 le fait d'\^etre propre \'equivaut \`a se
prolonger contin\^ument \`a la droite achev\'ee $\{-\infty \}\cup \R \cup \{+\infty \}$ en
envoyant les extr\'emit\'es
$+\infty$ et $-\infty$ dans $\fixe(h)$ (ou encore, $\phi$ est un
hom\'eo\-mor\-phisme sur son image, qui est ferm\'ee dans $\S^2 \setminus
\fixe(h)$).
 On confondra souvent le plongement
$\phi$ et son image. 
D'autre part, remarquons que  l'adh\'erence (dans $\S^2$) d'une droite topologique
 est un cercle topologique ou un arc, selon qu'elle contient $1$ ou $2$
points fixes.

\begin{defi}\index{s\'eparer}
Si $M$ est un espace topologique et $A$, $B$ et $C$ sont trois
sous-espaces disjoints, on dit que $C$ \res{s\'epare} $A$ et $B$ si $A$ et
$B$ sont inclus  dans deux composantes connexes distinctes du
compl\'ementaire de $C$.\footnote{La relation ``ne pas \^etre s\'epar\'es
par $C$'' est \'evidemment
une relation d'\'equivalence sur l'ensemble des parties connexes de $M$
disjointes de $C$.}
\end{defi}

\begin{defi}
\label{def.drbr}
\index{droite de \hbox{Brouwer}}
Une  droite topologique $\Delta$ est une \res{droite de \hbox{Brouwer}} si
\begin{enumerate}
\item  elle est libre par $h$ ;

\item dans le cas o\`u son adh\'erence  est un cercle topologique, elle
s\'epare son image $h(\Delta)$ et
sa pr\'eimage $h^{-1}(\Delta)$.
\end{enumerate}
\end{defi}

\begin{figure}[htp]
\par
\centerline{\hbox{\input{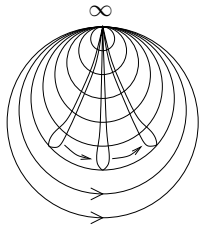}}}
\par
\caption{\label{fig15p}Une fausse droite de \hbox{Brouwer}}
\end{figure}
On peut voir des exemples de droites de \hbox{Brouwer} sur les figures
\ref{fig35} et \ref{fig18}. La figure \ref{fig15p} repr\'esente une droite
topologique libre pour la translation $\tau: (x,y) \mapsto (x+1,y)$
sur la sph\`ere $\R^2 \cup \{ \infty \}$, mais qui n'est pas une droite de \hbox{Brouwer}, montrant la
n\'ecessit\'e de la deuxi\`eme hypoth\`ese.

\subsubsection*{$\bullet$ Propri\'et\'es}
\begin{rema}\label{rem.orie}
Soit $\Delta_0$ une droite topologique libre dont l'adh\'erence est un
cercle topologique, et notons $\vec
\Delta_0$ la droite $\Delta_0$ orient\'ee de mani\`ere \`a ce que
$h(\Delta_0)$ soit situ\'e \`a droite de $\vec \Delta_0$. Alors la
deuxi\`eme condition de la d\'efinition~\ref{def.drbr} est v\'erifi\'ee (et
$\Delta_0$ est une droite de \hbox{Brouwer}) si et seulement si $\Delta_0$
est \`a gauche de $h(\vec \Delta_0)$. Ceci est d\^u au fait que $h$
pr\'eserve l'orientation.
\end{rema}
\begin{rema}\label{rem.attr}
Soit $\Delta_0$ une droite topologique libre dont l'adh\'erence est un
cercle topologique, et notons $P^+(\Delta_0)$ l'unique disque
topologique ferm\'e de la sph\`ere, de fronti\`ere $\adhe(\Delta_0)$, qui
contient $h(\Delta_0)$. Alors $\Delta_0$  est une droite de \hbox{Brouwer} si
et seulement si $P^+(\Delta_0)$ est un attracteur strict. Cet
attracteur sera baptis\'e  \emph{p\'etale attractif de $\Delta_0$}.
\end{rema}

Le lemme~\ref{lem.libr} entra\^{\i}ne imm\'ediatement~:
\begin{rema}
\label{rem.dbli}
Une droite de \hbox{Brouwer} est disjointe de tous ses it\'er\'es. En
particulier, avec la remarque~\ref{rem.attr}, on voit qu'une
 droite de \hbox{Brouwer} pour $h$ est une droite de \hbox{Brouwer} pour
$h^n$ pour tout entier $n \neq 0$.
\end{rema}

L'int\'er\^et des droites de \hbox{Brouwer} r\'eside dans l'existence d'une zone
o\`u la dynamique est conjugu\'ee \`a une translation~:
\begin{figure}[hbtp]
\par
\centerline{\hbox{\input{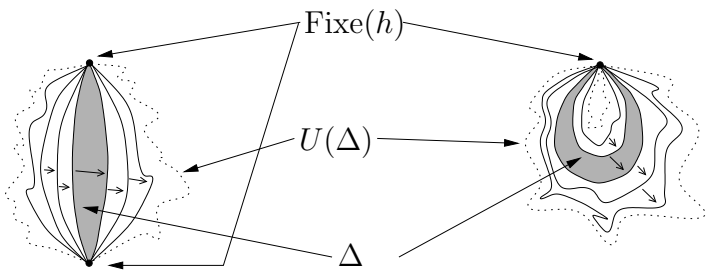}}}
\par
\caption{\label{fig15}Droites et domaines de \hbox{Brouwer}}
\end{figure}

\begin{affi}[figure~\ref{fig15}] 
\label{aff.dofo}
Soit $\Delta$ une droite de \hbox{Brouwer}.
\begin{enumerate}
\item Il existe un unique disque
topologique ouvert $D$ tel que 
\begin{itemize}
\item $\partial D \setminus\fixe(h)=\Delta \cup h(\Delta) $ ;
\item $D \cap h^{-1}(\Delta)=\emptyset$.
\end{itemize}
De plus, ce disque  est libre.

\item 
\label{aff.cotr}
\index{translation}Soit $$U(\Delta)=\bigcup_{n \in \Z}h^n(D\cup \Delta).$$
Il existe un hom\'eo\-mor\-phisme $\phi : \R^2 \rightarrow U(\Delta)$ qui conjugue la
translation $\tau : (x,y) \mapsto (x+1,y)$ et la restriction $h :
U(\Delta) \rightarrow U(\Delta)$. De plus, 
\begin{itemize}
\item les ensembles $\phi(\{x\} \times \R)$ sont des droites de \hbox{Brouwer} ;
\item  $\phi(\{0\} \times \R)=\Delta$. 
\end{itemize}
\end{enumerate}
\end{affi}
\begin{defi}
\index{domaine fondamental}
\index{domaine de \hbox{Brouwer}}
Le disque $D$ de l'affirmation pr\'ec\'edente est appel\'e \res{domaine fondamental de
$\Delta$}, et not\'e $D(\Delta,h(\Delta))$. L'ouvert $U(\Delta)$
est appel\'e \res{domaine de \hbox{Brouwer} engendr\'e par $\Delta$}.
\end{defi}
L'ouvert $U(\Delta)$ est hom\'eomorphe au plan, mais 
sa fronti\`ere n'est en g\'en\'eral pas localement connexe.

\pagebreak
\begin{demo}[de l'affirmation \ref{aff.dofo}]
\paragraph{Premier point}
Soit $\Delta$ une droite de \hbox{Brouwer} ; on fait la preuve dans le cas o\`u
l'adh\'erence de $\Delta$ contient deux points fixes
distincts $N$ et $S$ (l'autre cas \'etant similaire). L'ensemble $\gamma=\Delta \cup h(\Delta) \cup
\{N,S\}$ est une courbe de Jordan, et $h^2(\Delta)$ est disjointe de
$\gamma$ (remarque \ref{rem.dbli}). Appelons $D$ celui des 
deux disques topologiques ouverts de
fronti\`ere $\gamma$ qui ne contient pas $h^2(\Delta)$.

Montrons que $D$ est libre. La fronti\`ere de $h(D)$ est
$h(\gamma)$, elle ne rencontre pas $D$. Il y a donc deux
possibilit\'es : ou bien $D \subset h(D)$, ou bien
$D \cap h(D)=\emptyset$ ; il  reste \`a montrer que le
premier cas est impossible.

Ceci est d\^u au fait que $h$ pr\'eserve l'orientation. En effet,
notons $\vec \Delta$ la droite $\Delta$ orient\'ee de $S$ vers $N$ ; ceci d\'efinit localement
un c\^ot\'e droit et un c\^ot\'e gauche de $\vec \Delta$. L'image
$h(\vec \Delta)$ \'egalement orient\'ee de $S$ vers $N$
(parce que $h$ fixe $N$ et
$S$). Supposons par exemple que $D$ soit du c\^ot\'e droit de 
 $\vec \Delta$. Alors $D$ est du c\^ot\'e gauche de $h(\vec \Delta)$
(regarder l'allure au voisinage de $N$), mais $h(D)$ est du c\^ot\'e
droit de $h(\vec \Delta)$ puisque $h$ pr\'eserve l'orientation. Ceci montre
que le premier cas est impossible (m\^eme raisonnement si $D$ est du
c\^ot\'e gauche de $\Delta$).

\paragraph{Deuxi\`eme point}
On d\'efinit le plongement $\phi$ en envoyant le domaine fondamental de
la translation $[0,1] \times \R$ sur 
$\Delta \cup D \cup h(\Delta)$, en respectant la condition de
conjugaison sur le  bord du domaine ; il y a ensuite une unique
mani\`ere de prolonger $\phi$ au plan en remplissant la condition de
conjugaison.
 L'affirmation pr\'ec\'edente et le lemme \ref{lem.libr} sur les it\'er\'es
des disques libres 
montrent alors que $\phi$ est injective. Les
d\'etails sont laiss\'es au lecteur.
\end{demo}

\subsection{Th\'eor\`eme de translation plane}
\label{sub.entp}
On peut maintenant \'enoncer le th\'eor\`eme de translation plane de
\hbox{Brouwer}. L'\'enonc\'e d'origine de \hbox{Brouwer} concerne le cas avec un seul
point fixe ; l'\'enonc\'e g\'en\'eralis\'e ci-dessous est d\^u \`a Slaminka 
(voir \cite{slam3}). Par ailleurs, rappelons que L. Guillou a  encore \'etendu l'\'enonc\'e \`a tout
hom\'eo\-mor\-phisme libre d'une surface compacte (\cite{guil2}).
\begin{theo}[\hbox{Brouwer}]
\index{translation plane!\'enonc\'e}
Soit $h$ un \homeo\ de la sph\`ere, pr\'eservant l'orien\-ta\-tion, ayant un
nombre fini de point fixe, sans courbe d'indice $1$. Alors pour tout
$x \in \S^2 \setminus \fixe(h)$ il existe une droite de \hbox{Brouwer}
$\Delta$ contenant $x$.
\end{theo}
De mani\`ere \'equivalente, tout point qui n'est pas fixe est dans un
domaine de \hbox{Brouwer} (\emph{via} l'affirmation \ref{aff.dofo}, deuxi\`eme point).
Une des preuves modernes, due \`a P. Le Calvez et A. Sauzet, sera donn\'ee
\`a la section \ref{sub.trpl} comme premi\`ere application des d\'ecom\-po\-si\-tions
en briques.


\section{Indices partiels}
\label{sec.inpa}
\index{indice partiel}
\emph{\small Les sections \ref{sec.deco},
\ref{sec.cam}, \ref{sec.dbns}, \ref{sec.reca} 
 sont ind\'ependantes de cette section du texte.}

Dans cette section, \textbf{on suppose que $h$ est un hom\'eo\-morphisme
de la sph\`ere, pr\'eservant l'orientation, et ayant un unique point
fixe, qu'on note $\infty$}. L'ensemble de ces donn\'ees sera appel\'e
``\textbf{hypoth\`ese (H1)}''.\index{hypoth\`ese (H1)}

D'apr\`es la formule de
Lefschetz,
 le point fixe est d'indice $2$, et les r\'esultats des paragraphes 
pr\'ec\'edents (\ref{sub.recu}, \ref{sub.drbr}, \ref{sub.entp}) s'appliquent \`a
$h$. On identifiera $\S^2 \setminus \{\infty\}$ au plan orient\'e $\R^2$ \emph{via} un
hom\'eo\-mor\-phisme pr\'eservant l'orientation, ce qui fait de  $h$
un \emph{hom\'eo\-mor\-phisme de \hbox{Brouwer}} :
\begin{defi}
\index{hom\'eo\-mor\-phisme de \hbox{Brouwer}}
Un \res{hom\'eo\-mor\-phisme de \hbox{Brouwer}} est un hom\'eo\-mor\-phisme du plan,
pr\'eservant l'orientation, sans point fixe.
\end{defi}

Le th\'eor\`eme principal relie l'indice du point fixe $N$ au nombre de
croissants. Pour sa preuve, lorsqu'il  faudra montrer qu'on a obtenu
le  nombre de croissants voulu, nous devrons calculer la contribution
de chaque croissant et de chaque zone entre deux croissants adjacents
\`a l'indice du point $N$. A cet effet, nous introduisons et \'etudions
dans cette section une notion d'indice partiel.

En r\'ealit\'e, nous effectuerons le calcul apr\`es
un passage au rev\^etement universel de $\S^2 \setminus\{N,S\}$~; ceci
explique le cadre de cette section, celui des hom\'eo\-morphismes de \hbox{Brouwer}. 
 Dans ce cadre, l'indice partiel va \^etre d\'efini comme
un invariant de conjugaison associ\'e \`a n'importe quel couple de droites de \hbox{Brouwer}
disjointes.

\subsection{Topologie des couples de droites de \hbox{Brouwer} disjointes}
\index{droite de \hbox{Brouwer}!couple}
Soient $\Delta_0$ et $\Delta_1$ deux droites de \hbox{Brouwer}
disjointes.
D'apr\`es le th\'eo\-r\`eme de Schoenflies-Homma (voir l'appendice),
 \index{Schoenflies-Homma!th\'eo\-r\`eme de}
 il n'y a qu'une seule
configuration topologique possible pour le couple $(\Delta_0,\Delta_1)$.
Nous faisons quelques consid\'erations pr\'eliminaires sur la topologie de
ces deux droites et de leurs it\'er\'es.

\begin{defi}\index{$D(\Delta_0,\Delta_1)$}
On appellera \res{disque d\'elimit\'e par $\Delta_0$ et $\Delta_1$}, et on
notera $D(\Delta_0,\Delta_1)$, l'unique disque topologique ouvert de
fronti\`ere $\Delta_0 \cup \Delta_1$.
\end{defi}

\begin{figure}[hp]
\par
\centerline{\hbox{\input{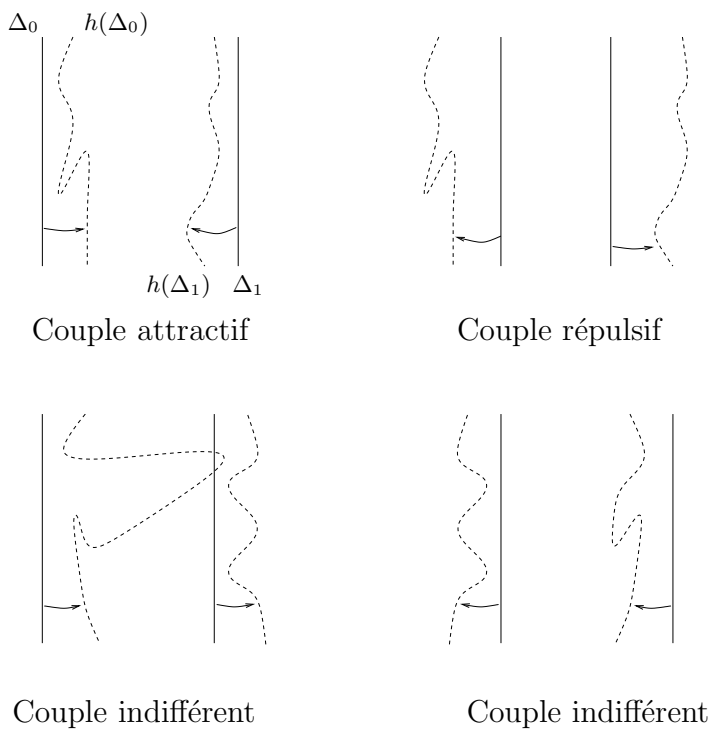}}}
\par
\caption{\label{fig25} Couples de droites de \hbox{Brouwer} disjointes}
\end{figure}

\begin{defi}[ (figure~\ref{fig25})]
\label{def.coat}
\index{droite de \hbox{Brouwer}!attractive}
On dit que $\Delta_0$ est \res{de type attractif} ou \res{attractive}
(sous-entendu, \res{relativement \`a $\Delta_1$})
si elle ne s\'epare pas $\Delta_1$ et $h(\Delta_0)$ (de mani\`ere
\'equivalente, si elle s\'epare $\Delta_1$ et $h^{-1}(\Delta_0)$). Dans le
cas contraire, on dit que $\Delta_0$ est \res{r\'epulsive}. Les m\^emes
d\'efinitions s'appliquent \`a $\Delta_1$.

Si les types 
des deux droites sont oppos\'es,  le couple
$(\Delta_0,\Delta_1)$ est dit \res{indiff\'erent} ; sinon, il est
qualifi\'e d'\res{attractif}
ou de \res{r\'epulsif} selon le type commun.
\end{defi}

\begin{figure}[htp]
\par
\centerline{\hbox{\input{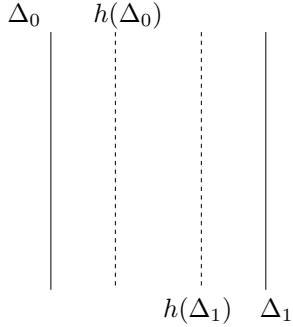}}}
\par
\caption{\label{fig.couple-attractif} Topologie d'un couple attractif}
\end{figure}

\begin{affi}[Topologie des couples attractifs] \label{aff.topo}
Soit $(\Delta_0,\Delta_1)$ un couple attractif de droites de \hbox{Brouwer}
disjointes. Alors l'ensemble $\adhe(D(\Delta_0,\Delta_1))$ est un attracteur
strict, et la configuration des droites $\Delta_0, h(\Delta_0),
\Delta_1$ et $h(\Delta_1)$ est hom\'eomorphe \`a celle de la
figure~\ref{fig.couple-attractif}. 
\end{affi}

\begin{demo}
Avec les notations de la remarque~\ref{rem.attr},  on a 
$\adhe(D(\Delta_0,\Delta_1))=P^+(\Delta_0) \cap P^+(\Delta_1)$. Cet
ensemble  est donc un attracteur strict. En particulier, on a
$h(\Delta_0) \subset D(\Delta_0,\Delta_1)$, et les quatre droites
  $\Delta_0, h(\Delta_0), \Delta_1$ et $h(\Delta_1)$ sont deux \`a deux
disjointes.

Il reste \`a \'etudier la configuration topologique.
D'apr\`es le th\'eo\-r\`eme de Schoenflies-Homma (voir l'appendice),
\index{Schoenflies-Homma!th\'eo\-r\`eme de} il suffit de montrer  que $h(\Delta_0)$ s\'epare
$\Delta_0$ de $\Delta_1$, et que $h(\Delta_1)$ s\'epare $h(\Delta_0)$ de
$\Delta_1$.

Montrons que $h(\Delta_0)$ s\'epare $\Delta_0$ et
$\Delta_1$. Supposons le  contraire~:
\begin{enumerate}
\item la situation est hom\'eomorphe \`a la
figure~\ref{fig.bd} (a) (d'apr\`es le th\'eor\`eme de Schoenflies-Homma,
\index{Schoenflies-Homma!th\'eo\-r\`eme de} quitte \`a renverser
l'orientation)~; 
\item on oriente
$\Delta_0$ et $\Delta_1$ comme sur la figure~\ref{fig.bd} (b), et on
en d\'eduit l'orientation de $h(\Delta_0)$ (d'apr\`es la
remarque~\ref{rem.orie})~; 
\item comme $h$ pr\'eserve l'orientation, et que $\Delta_1$ est \`a droite
de $\vec \Delta_0$, la droite $h(\Delta_1)$ est aussi \`a droite de
$h(\vec \Delta_0)$, et la situation est hom\'eomorphe \`a
la figure~\ref{fig.bd} (c)~;
\item l'orientation de $h(\vec \Delta_1)$ est alors celle de la figure~\ref{fig.bd}~(d)
 (d'apr\`es la remarque~\ref{rem.orie}).
\end{enumerate}
\begin{figure}[htp]
\par
\centerline{\hbox{\input{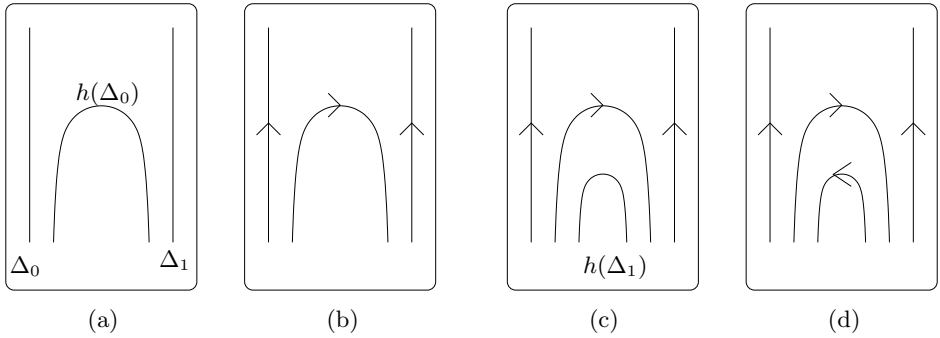}}}
\par
\caption{\label{fig.bd}Si $h(\Delta_0)$ ne s\'epare pas $\Delta_0$ et
$\Delta_1$...}
\end{figure}
Mais cette derni\`ere figure est contradictoire, puisque $\Delta_0$ est
\`a gauche de $\vec \Delta_1$, alors que $h(\Delta_0)$ est \`a droite de
$h(\vec \Delta_1)$, et que $h$ pr\'eserve l'orientation.

On montre de la m\^eme mani\`ere que $h(\Delta_1)$ s\'epare $h(\Delta_0)$
et $\Delta_1$, ce qui termine la preuve.
\end{demo}

Remarquons qu'il n'y a \'egalement qu'une seule configuration topologique possible
pour les $n$ premiers it\'er\'es de $\Delta_0$ et $\Delta_1$, pour tout
entier $n$.

\subsection{D\'efinition de l'indice partiel}
Soient $\Delta_0$ et $\Delta_1$ deux droites de \hbox{Brouwer}
disjointes. Nous allons d\'efinir 
l'indice partiel de $h$ entre $\Delta_0$ et
 $\Delta_1$ en examinant la situation au moyen de
``cartes'' qui redressent ces deux droites ; le jeu consistera
alors \`a montrer que les objets construits ne d\'ependent pas du choix
des cartes. Nous ferons ici un usage intensif des coordonn\'ees
euclidiennes du plan $\R^2$.

 D'apr\`es le th\'eor\`eme de Schoenflies-Homma (voir l'appendice),
\index{Schoenflies-Homma!th\'eo\-r\`eme de} il existe un hom\'eo\-mor\-phisme
$g'$ du plan, pr\'eservant l'orientation, tel que les droites
$\Delta_0'=g'(\Delta_0)$ et $\Delta_1'=g'(\Delta_1)$ sont
verticales, et $\Delta_0'$ est \`a gauche de $\Delta_1'$ (\emph{i. e.}
$\Delta_0'=\{x_0\} \times \R$, $\Delta_1'=\{x_1\} \times \R$, avec
$x_0 < x_1$, ce que dor\'enavant nous \'ecrirons $\Delta_0' <
 \Delta_1'$). On pose  $h'=g' \circ h \circ g'^{-1}$.

Soit $\gamma'$
 une courbe allant de $\Delta'_0$ \`a $\Delta'_1$ (on ne suppose
 pas que $\gamma$ est simple).
On distingue deux cas (figure~\ref{fig24})~:
\begin{figure}[hbtp]
  \par
\centerline{\hbox{\input{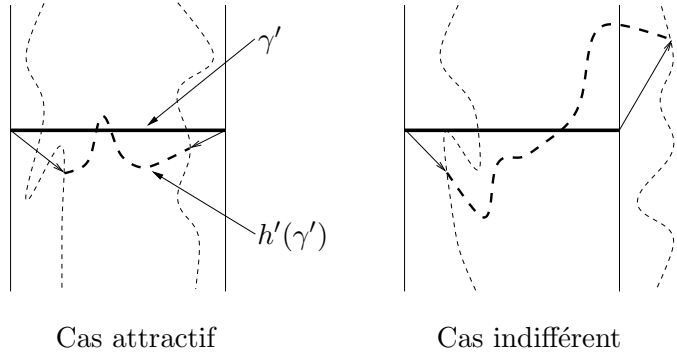}}}
\par
  \caption{\label{fig24}D\'efinition du nombre $I$}
\end{figure}
\begin{itemize}
\item[Premier cas.] Le couple $(\Delta_0, \Delta_1)$ est indiff\'erent. Dans ce cas,
 les vecteurs $h'(\gamma'(0))-\gamma'(0)$ et $h'(\gamma'(1))-\gamma'(1)$ 
 sont tous les deux d'abscisses strictement positive ou bien tous les
deux d'abscisses strictement n\'egative ; en tous cas, ils ne peuvent
pas avoir des directions oppos\'ees, et
on a $\indi(h',\gamma') \not \in (1/2 + \Z)$. On d\'efinit
le nombre $I$ comme l'entier le plus proche de $\indi(h',
\gamma')$.
\item[Second cas. ] Le couple $(\Delta'_0, \Delta'_1)$ est attractif ou r\'epulsif.
Cette fois-ci,
on a $\indi(h',\gamma') \not \in \Z$,
et on d\'efinit $I$ comme le demi-entier le plus proche de
$\indi(h',\gamma')$.
\end{itemize}

\begin{defi}\label{def.inpa}
\index{$\inpa(h,\Delta_0,\Delta_1)$}
 On appelle \res{indice partiel de $h$ entre $\Delta_0$ et $\Delta_1$},
 et on note $\inpa(h,\Delta_0,\Delta_1)$, le nombre $I$.
\end{defi}
Nous allons montrer que cette d\'efinition est correcte, c'est-\`a-dire
que le nombre $I$ ne d\'epend pas des choix de la courbe  $\gamma'$ et
de l'hom\'eo\-morphisme $g'$.

\begin{affi}
Le nombre $I$ ne d\'epend pas du choix de la courbe $\gamma'$.
\end{affi}
\begin{demo}
L'ensemble des choix possibles pour $\gamma'$ est un espace de
courbes qui est connexe (pour la topologie de la convergence
uniforme) : c'est m\^eme un espace convexe pour la structure
affine naturelle.
 Or $\indi(h', \gamma')$ d\'epend contin\^ument de 
$\gamma'$ (pour la m\^eme topologie). Dans le cas d'un couple 
indiff\'erent, la composante connexe de $\R \setminus(\Z+1/2)$ qui
contient $\indi(h',\gamma')$ ne d\'epend donc pas de
$\gamma'$, et le nombre $I$ non plus. Le raisonnement est analogue 
dans le cas attractif ou r\'epulsif.
\end{demo}

\begin{affi}\label{aff.inco}
Le nombre $I$ ne d\'epend pas du choix de l'hom\'eo\-mor\-phisme  $g'$.
\end{affi}
\index{espace d'hom\'eo\-mor\-phismes}
Il s'agit encore d'un argument de con\-nexi\-t\'e. Appelons $\cal G$
l'ensemble des hom\'eo\-mor\-phismes  $G$ du plan, pr\'eservant l'orientation,
tels que $G(\Delta'_0)$ et $G(\Delta'_1)$ soient encore deux droites
verticales avec $G(\Delta'_0) < G(\Delta'_1)$. On a :
\begin{lemm}\label{lem.tcvu}
L'espace $\cal G$, muni de la topologie de la convergence uniforme sur
les compacts, est connexe par arcs.
\end{lemm}
\begin{demo}[du lemme~\ref{lem.tcvu}]
Soit $G \in {\cal G}$, on cherche un chemin dans $\cal G$ de $G$ \`a l'identit\'e.
On peut d'abord se ramener  au cas ou $G$ fixe
globalement $\Delta'_0$ et $\Delta'_1$ (\`a l'aide d'une isotopie
affine), puis au cas
o\`u $G$ fixe tous les points de ces deux droites.
 On termine par une variante de l'isotopie d'Alexander (\cite{alex1},
ou bien \cite{lero1}).
\end{demo}

\begin{demo}[de l'affirmation \ref{aff.inco}]
On fait la preuve dans le cas indiff\'erent, l'argument \'etant le m\^eme
dans l'autre situation.
Soit $g''$ un autre hom\'eo\-mor\-phisme ayant les propri\'et\'es requises
(pr\'eserver l'orientation et envoyer les deux droites de \hbox{Brouwer}
$\Delta_0$ et $\Delta_1$ sur
deux droites verticales avec $g''(\Delta_0) < g''(\Delta_1)$).
 Posons $G=g'' \circ g'^{-1}$ ; alors $G \in
{\cal G}$. D'apr\`es le lemme, il existe une isotopie $(G_t)_{t \in
[0,1]}$ telle que $G_0=G$ et $G_1$ est l'identit\'e. 
Soit $\gamma$ une courbe allant de $\Delta_0$ \`a $\Delta_1$,
$\gamma'=g'(\gamma)$ et  $\gamma''=g''(\gamma)$. Le nombre 
$I(t)=\indi(G_t h' G_t^{-1}, G_t(\gamma'))$ varie contin\^ument, il ne prend
pas de valeurs dans $\Z + 1/2$, donc les deux nombres 
$I(0)=\indi(g'' h g''^{-1},\gamma''))$ et  $I(1)=\indi(g' h
g'^{-1},\gamma')$ sont dans la m\^eme composante connexe de $\R
\setminus (\Z + 1/2)$, ce que l'on voulait montrer. 
\end{demo}

L'affirmation \ref{aff.inco} entra\^\i ne imm\'ediatement :
\begin{coro}\label{cor.inco}
L'indice partiel est un invariant de conjugaison (orient\'ee)~:
autrement dit, si $g$ est un
hom\'eo\-mor\-phisme du plan pr\'eservant l'orientation, on a
$$
\inpa(h,\Delta_0,\Delta_1)=\inpa(ghg^{-1},g(\Delta_0),g(\Delta_1)).
$$ 
\end{coro}

On a aussi, de mani\`ere imm\'ediate~:
\begin{affi}\hspace{1cm}
\label{aff.ipfa}
\begin{enumerate}
\item $\inpa(h,\Delta_1,\Delta_0)=-\inpa(h,\Delta_0,\Delta_1)$ ;
\item $\inpa(h^{-1},\Delta_0,\Delta_1)=\inpa(h,\Delta_0,\Delta_1)$.
\end{enumerate}
\end{affi}

\subsection{Exemples}\label{sub.exip}
On peut voir facilement que pour la translation $\tau : (x,y) \mapsto
(x+1,y)$, tout couple de droites de \hbox{Brouwer} disjointes est d'indice
partiel $0$, $1/2$ ou $-1/2$. Ceci sera d'ailleurs une cons\'equence du
lemme \ref{lem.arli}.

\begin{figure}[htpb]
\par
\centerline{\hbox{\input{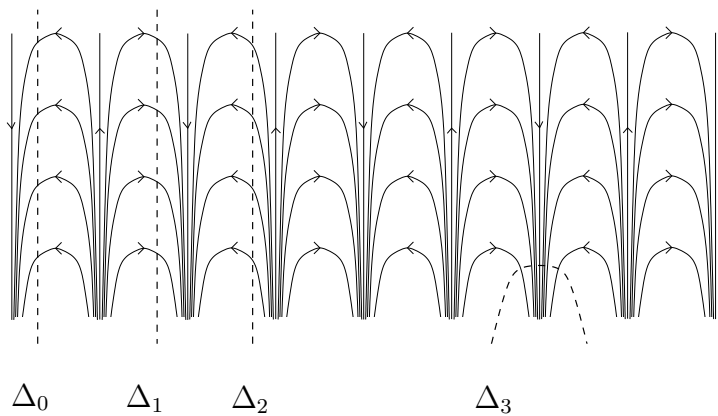}}}
\par
\caption{\label{fig35}Exemples d'indices partiels}
\end{figure}

Pour l'hom\'eo\-mor\-phisme ``multi-Reeb'' repr\'esent\'e sur la figure~\ref{fig35}~:
\begin{enumerate}
\item le couple $(\Delta_0, \Delta_1)$ est r\'epulsif d'indice
partiel $-1/2$ ;
\item  le couple $(\Delta_1, \Delta_2)$ est attractif d'indice
partiel $-1/2$ ;
\item le couple $(\Delta_0, \Delta_2)$ est indiff\'erent d'indice
partiel $-1$ ;
\item le couple $(\Delta_3, \Delta_2)$ est r\'epulsif d'indice
partiel $3/2$.
\end{enumerate}

Ce type d'exemple montre que l'indice partiel peut prendre 
toute  valeur $n/2$ avec $n$ entier. Un des principaux r\'esultats de
cet article montrera que si $\mid \inpa(\Delta_0,\Delta_1)\mid >1/2$,
  alors on peut d\'ecouper la bande $D(\Delta_0,\Delta_1)$ en tranches, \`a l'aide
d'une famille de droites de \hbox{Brouwer}, de mani\`ere \`a ce que l'indice
partiel dans chaque tranche soit \'egal \`a $0$, $1/2$ ou $-1/2$.

\subsection{Relation d'additivit\'e}
\index{relation de Chasles}
L'efficacit\'e de l'indice partiel dans les calculs d'indice repose de
mani\`ere cruciale sur la propri\'et\'e d'additivit\'e suivante~:
\begin{lemm}[relation de Chasles]
\label{lem.chas}
 Soient $\Delta_0$, $\Delta_1$ et $\Delta_2$ trois droites de
 \hbox{Brouwer} disjointes, et supposons que $\Delta_1$ s\'epare $\Delta_0$ et
 $\Delta_2$. On a alors :
$$
\inpa(h,\Delta_0,\Delta_2)=\inpa(h,\Delta_0,\Delta_1)+\inpa(h,\Delta_1,\Delta_2).
$$
\end{lemm}
\begin{figure}[h!tpb]
  \par
\centerline{\hbox{\input{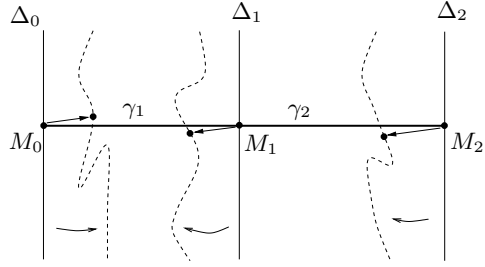}}}
\par
  \caption{\label{fig26}La relation de Chasles}
\end{figure}

\begin{demo}
Commençons par une remarque banale~: dans la construction conduisant \`a
la d\'efinition de l'indice partiel (d\'efinition~\ref{def.inpa},
figure~\ref{fig24}), si chacune des extr\'emit\'es de la
courbe $\gamma'$ est sur la m\^eme droite horizontale que son image
par $h'$, alors l'indice partiel est exactement \'egal \`a $\indi(h',\gamma')$.

Puisque l'indice partiel est un invariant de conjugaison,  on peut
supposer que $\Delta_0$, $\Delta_1$ et $\Delta_2$ sont trois droites
verticales v\'erifiant $\Delta_0 < \Delta_1 < \Delta_2$ (quitte \`a conjuguer en utilisant
le th\'eor\`eme de Schoenflies-Homma \index{Schoenflies-Homma!th\'eo\-r\`eme de}).
Soient $M_0$, $M_1$ et $M_2$ les points d'intersection respectifs des
droites $\Delta_0$, $\Delta_1$, $\Delta_2$ avec l'axe des
abscisses. En conjuguant \`a nouveau  $h$ par une application $(x,y) \mapsto
(x, y/k)$ avec $k$ tr\`es grand, on peut supposer que les trois vecteurs
$h(M_i)-M_i$ forment des angles arbitrairement petits avec la
direction horizontale (figure~\ref{fig26}).

Soient alors $\gamma_1$ une courbe allant de $M_0$ \`a $M_1$,
$\gamma_2$ une courbe allant de $M_1$ \`a $M_2$, et $\gamma$ la courbe
obtenue en mettant les deux pr\'ec\'edentes bout \`a bout. On a bien s\^ur
$$
\indi(h,\gamma)=\indi(h,\gamma_1)+\indi(h,\gamma_2).
$$
On utilise maintenant la remarque initiale~: dans la carte choisie,
\begin{itemize}
\item le nombre $\inpa(h,\Delta_0,\Delta_2)$ est
arbitrairement proche de $\indi(h,\gamma)$~; 
\item  le nombre 
$\inpa(h,\Delta_0,\Delta_1)+\inpa(h,\Delta_1,\Delta_2)$  est
arbitrairement proche de $\indi(h,\gamma_1)+\indi(h,\gamma_2)$.
\end{itemize}
Or ces deux nombres sont des demi-entiers~:
si l'\'egalit\'e du lemme \'etait fausse, ils diff\'ereraient
d'au moins $1/2$, ce qui est absurde.
\end{demo}

\subsection{Indice partiel et arcs libres}
\label{ss.ipal}
\index{arc libre}

Nous \'enonçons ici un lemme qui  permettra notamment de calculer
l'indice partiel \`a travers un croissant attractif ou r\'epulsif \`a
dynamique Nord-Sud ou Sud-Nord.

Soient $\Delta_0$ et $\Delta_1$ deux droites de \hbox{Brouwer} disjointes. 
On suppose qu'il existe un arc libre $\gamma$ allant de $\Delta_0$
\`a $\Delta_1$. Puisqu'un sous-arc d'un arc libre est
encore libre, on peut en fait demander de plus que l'int\'erieur de  $\gamma$
soit disjoint des deux droites.

A l'aide du th\'eor\`eme de Schoenflies-Homma (quitte \`a conjuguer la
situation par un hom\'eo\-morphisme pr\'eservant l'orientation), on peut alors supposer
que les deux droites de \hbox{Brouwer} sont verticales, avec $\Delta_0 <
\Delta_1$, et aussi que $\gamma$ est un segment horizontal.

Si le couple $(\Delta_0,\Delta_1)$ est attractif, la disque
$D(\Delta_0,\Delta_1)$ qu'il
d\'elimite est un attracteur (d'apr\`es l'affirmation~\ref{aff.topo}), et
$h(\gamma)$ est un arc disjoint de $\gamma$ et inclus dans  ce
disque. Il y a donc deux possibilit\'es :  ou bien $h(\gamma)$ est situ\'e
au-dessus de $\gamma$, ou bien  $h(\gamma)$ est situ\'e
au-dessous de $\gamma$. On montre que la situation ne d\'epend pas du
choix de la ``carte'' (pr\'eservant l'orientation)
 donn\'ee par le th\'eor\`eme de Schoenflies-Homma.
Ceci vaut \'egalement dans le cas r\'epulsif en
rempla\c{c}ant $h$ par $h^{-1}$.

\begin{lemm}
\label{lem.arli}Soient $\Delta_0$ et $\Delta_1$ deux droites de
\hbox{Brouwer} disjointes, et $\gamma$ un arc libre allant de $\Delta_0$
\`a $\Delta_1$, dont l'int\'erieur est disjoint des deux droites.

Si $(\Delta_0,\Delta_1)$ est indiff\'erent, l'indice partiel entre   $\Delta_0$
et $\Delta_1$ est nul. 

Si $(\Delta_0,\Delta_1)$ n'est pas indiff\'erent, il y a quatre cas :
\begin{itemize}
\item si le couple est attractif et que $h(\gamma)$ est au-dessus de
$\gamma$, l'indice partiel vaut $+1/2$ ;
\item si le couple est attractif et que $h(\gamma)$ est au-dessous de
$\gamma$, l'indice partiel vaut $-1/2$ ;
\item si le couple est r\'epulsif et que $h^{-1}(\gamma)$ est au-dessus de
$\gamma$, l'indice partiel vaut $+1/2$ ;
\item si le couple est r\'epulsif et que $h^{-1}(\gamma)$ est au-dessous de
$\gamma$, l'indice partiel vaut $-1/2$.
\end{itemize}
\end{lemm}
\begin{demo}
On fait la preuve dans le cas indiff\'erent, les autres cas \'etant
similaires,
et on suppose par exemple
que $\Delta_0$ est attractive et $\Delta_1$ 
r\'epulsive~; autrement dit, pour tout point $x$ de $\Delta_0$ ou
$\Delta_1$,  le vecteur $h(x)-x$ est d'abscisse strictement positive
(figure \ref{fig27}).
Pour tout $t \geq 0$, on pose 
$$
i_t=\indi \left(h(x)+(t,0)-x, \gamma  \right).
$$
\begin{figure}[hbtp]
  \par
\centerline{\hbox{\input{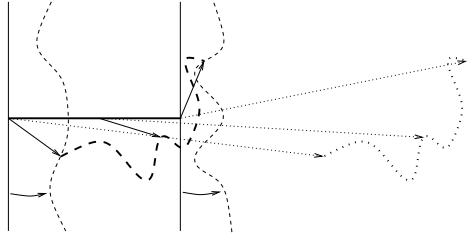}}}
\par
  \caption{\label{fig27}L'indice partiel est nul}
\end{figure}

 Comme $\gamma$ est libre pour $h$ et que $\Delta_0$ est attractive,
 quand $x$  est un point de $\gamma$, le point $h(x)$ n'est jamais situ\'e
 sur la demi-droite horizontale \`a gauche
 de $x$, donc le vecteur $h(x)+(t,0)-x$
 n'est jamais nul, et le nombre $i_t$ est bien d\'efini. D'autre part
 $i_t \not \in (1/2+\Z)$, $i_t$ varie contin\^ument avec $t$, et
 $\lim_{t \rightarrow +\infty}i_t=0$, donc $i_t \in ]-1/2,+1/2[$ pour
 tout $t \geq 0$ ; comme $\inpa(h,\Delta_0,\Delta_1)$ est l'entier le
 plus proche de $i_0$, il est bien nul.
\end{demo}

\begin{coro}
\label{cor.dotr}
Soient $\Delta_0$ et $\Delta_1$  deux droites
de \hbox{Brouwer} disjointes.
\begin{enumerate}
\item Si le domaine de
\hbox{Brouwer} engendr\'e par $\Delta_0$ rencontre  $\Delta_1$, alors
$\inpa(h,\Delta_0, \Delta_1)=0$ (en particulier, pour  tout entier $n$ non
nul, $\inpa(h,\Delta_0, h^n(\Delta_0))=0$). 
\item Dans le cas contraire, l'indice partiel $\inpa(h,h^p(\Delta_0),
h^q(\Delta_1))$ est bien d\'efini pour tous entiers $p$ et $q$, et est
\'egal \`a $\inpa(h,\Delta_0, \Delta_1)$.
\end{enumerate}
\end{coro}
Remarquons que si le couple $(\Delta_0,\Delta_1)$ est attractif ou
r\'epulsif, on est toujours dans le deuxi\`eme cas (d'apr\`es l'affirmation 
\ref{aff.topo}).

\begin{demo}
\paragraph{Premier point}
D'apr\`es la remarque pr\'ec\'edente, dans ce cas, le couple
$(\Delta_0,\Delta_1)$ est indiff\'erent. Notons ensuite que  
si deux points du plan  appartiennent \`a un m\^eme domaine de \hbox{Brouwer}, alors ou
bien ils sont dans la m\^eme orbite, ou bien il existe un arc libre
les joignant : en effet, puisque sur un domaine de \hbox{Brouwer} la
dynamique est conjugu\'ee \`a une translation, il suffit de prouver
la m\^eme propri\'et\'e pour la translation, ce qui n'est pas
difficile (par exemple en raisonnant dans l'anneau quotient de
l'action du plan par la translation).
Ceci permet de montrer que si le domaine de
\hbox{Brouwer} engendr\'e par $\Delta_0$ rencontre  $\Delta_1$, alors
il existe un arc libre allant de $\Delta_0$ \`a $\Delta_1$, puis on
conclut \`a l'aide du lemme \ref{lem.arli}.

\paragraph{Deuxi\`eme point}
Consid\'erons d'abord le cas o\`u le couple $(\Delta_0,\Delta_1)$ est
attractif. Comme l'indice partiel est un invariant de conjugaison,
quitte \`a composer par $h^{-p}$,  on peut supposer que $p=0$. On
traite le cas o\`u $q$ est positif (le cas n\'egatif \'etant similaire). 
 En utilisant de mani\`ere it\'er\'ee
l'affirmation~\ref{aff.topo} sur la topologie des couples attractifs,
 on montre que $h^q(\Delta_1)$ s\'epare $\Delta_0$ et
$\Delta_1$. Ceci permet d'appliquer la relation de Chasles
(lemme~\ref{lem.chas}) \`a ces trois droites~; comme l'indice partiel
entre $\Delta_1$ et $h^q(\Delta_1)$ est nul d'apr\`es le premier point,
on obtient le r\'esultat.

On obtient le cas r\'epulsif de la m\^eme mani\`ere~; il reste le cas indiff\'erent.
On peut encore supposer que $p=0$, et on se place dans le cas o\`u $q$
est positif (le cas n\'egatif est similaire).
 Si la droite $\Delta_1$ est r\'epulsive, alors $\Delta_1$
s\'epare $\Delta_0$ et $h^q(\Delta_1)$, et on peut appliquer la relation
de Chasles et conclure comme dans le cas pr\'ec\'edent. 

Supposons enfin que $\Delta_1$ est
attractive.  Si  $h^q(\Delta_1)$  s\'epare  $\Delta_0$
et $\Delta_1$, on conclut comme avant \`a l'aide de la relation de Chasles.
Si $h^q(\Delta_1)$ ne s\'epare pas $\Delta_0$ et $\Delta_1$,
on voit facilement que $\Delta_0$ est inclus dans le domaine de
\hbox{Brouwer} engendr\'e par $\Delta_1$, et les deux indices partiels sont
nuls d'apr\`es le premier point.
\end{demo}

\subsection{Indice partiel entre deux droites qui se touchent}
\label{ss.inpat}
Dans certains calculs, il sera pratique de pouvoir parler de l'indice
partiel entre deux droites de \hbox{Brouwer} non disjointes, mais ``qui ne se
traversent pas''.

Soit $\Delta_0$ et $\Delta_1$ deux droites topologiques, et $P_0$ et $P_2$ les
deux composantes connexes du compl\'ementaire de $\Delta_1$. 
\begin{defi}\index{droites topologiques!qui se
traversent}\index{droites topologiques!qui se touchent}\label{def.trav}
On dira que les droites $\Delta_0$ et $\Delta_1$ \res{ne se traversent
pas} si $\Delta_0$ est contenue dans l'adh\'erence de $P_0$ ou de $P_2$. On
dira que deux droites \res{se touchent} si elles ne se traversent pas,
mais ne sont pas disjointes.
\end{defi}
On peut montrer que ces deux relations sont sym\'etriques. 

On \'etend alors la d\'efinition de l'indice partiel au cas de deux
droites qui se touchent, en posant simplement
$\inpa(\Delta_0,\Delta_1)=0$ dans ce cas.

\begin{defi}\label{def.sepa2}\index{s\'eparer}
Soient $\Delta_0$, $\Delta_1$, $\Delta_2$ trois droites topologiques
qui, deux \`a deux, ne se traversent pas. Soient $P_1$ et
$P_2$ les deux composantes connexes du compl\'ementaire de
$\Delta_1$.
On dira que $\Delta_1$
\res{s\'epare $\Delta_0$ et $\Delta_2$} si $\Delta_0 \subset \adhe(P_0)$
et $\Delta_2 \subset \adhe(P_2)$ (ou \textit{vice-versa})
\end{defi}

Ces d\'efinitions permettent de g\'en\'eraliser la relation de Chasles
(lemme~\ref{lem.chas})~:
\begin{lemm}[relation de Chasles g\'en\'eralis\'ee]
\label{lem.chas2}\index{relation de Chasles!g\'en\'eralis\'ee}
 On se donne trois droites de
 \hbox{Brouwer}  $\Delta_0$, $\Delta_1$ et $\Delta_2$ qui, deux \`a deux, ne se traversent pas,
 et on suppose que $\Delta_1$ s\'epare $\Delta_0$ et
 $\Delta_2$. On a alors :
$$
\inpa(h,\Delta_0,\Delta_2)=\inpa(h,\Delta_0,\Delta_1)+\inpa(h,\Delta_1,\Delta_2).
$$
\end{lemm}
\begin{demo}
On reprend les notations de la d\'efinition~\ref{def.sepa2}.
Il est clair que l'une des deux droites $h(\Delta_0)$ et $h^{-1}(\Delta_0)$ est incluse dans
$P_0$, on la note $\Delta'_0$. De m\^eme, on note $\Delta'_2$ celle des
deux droites $h(\Delta_2)$ et $h^{-1}(\Delta_2)$ qui est incluse dans
$P_2$. Les trois droites de \hbox{Brouwer} $\Delta'_0$, $\Delta_1$ et
$\Delta'_2$ sont maintenant disjointes deux \`a deux, et $\Delta_1$
s\'epare les deux autres~: on peut alors leur appliquer la premi\`ere relation de
Chasles (lemme~\ref{lem.chas}). La formule en d\'ecoule, \`a l'aide du
corollaire~\ref{cor.dotr}. Les d\'etails sont laiss\'es au lecteur.

\end{demo}

\subsection{Indice partiel et cha\^\i nes de disques}
Nous allons g\'en\'eraliser le lemme~\ref{lem.arli} liant indice
partiel et arcs libres \`a l'aide de la notion de
cha\^\i ne de pseudo-disques (d\'efinition \ref{def.cdpd}).

\begin{figure}[htpb]
\par
\centerline{\hbox{\input{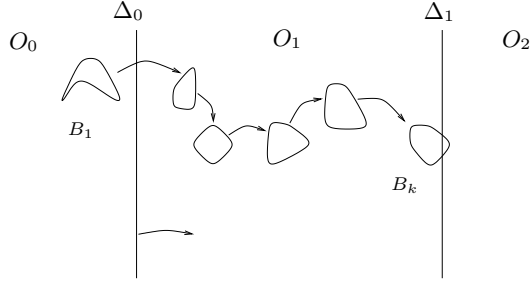}}}
\par
\caption{\label{fig28}Une cha\^\i ne de pseudo-disques de $\Delta_0$ \`a
$\Delta_1$}
\end{figure}

Soit $(\Delta_0,\Delta_1)$ un couple de  droites de \hbox{Brouwer}
disjointes ; on suppose que $\Delta_0$ est attractive. Le
com\-pl\'e\-men\-taire de $\Delta_0 \cup \Delta_1$ a trois composantes
connexes, dont les fronti\`eres sont respectivement $\Delta_0$,
$\Delta_0 \cup \Delta_1$ et $\Delta_1$ ; on les note $O_0$, $O_1$ et
$O_2$.
\begin{defi}[ (figure \ref{fig28})]\label{def.chtr}
\index{cha\^{\i}ne!entre deux droites de \hbox{Brouwer}}
Une \res{cha\^\i ne de pseudo-disques de $\Delta_0$ \`a $\Delta_1$} est une
cha\^\i ne de pseudo-disques $(\ii{B_i})_{i=1..k}$, ($k \geq 1$), 
telle que $\adhe(\ii{B_1})$ rencontre $\adhe(O_0)$ et $\adhe(\ii{B_k})$
rencontre $\adhe(O_2)$.
\end{defi}

Remarquons que si $k=1$, $B_1$ contient un arc libre allant de
$\Delta_0$ \`a $\Delta_1$~; dans le cas indiff\'erent, la proposition suivante est donc bien une
g\'en\'eralisation du lemme~\ref{lem.arli}~:
\begin{prop}[``lemme de Franks'' pour l'indice partiel]
\label{pro.frip}
\index{lemme de Franks!pour l'indice partiel}
On se donne deux droites de \hbox{Brouwer} $\Delta_0$ et $\Delta_1$  disjointes,
avec $\Delta_0$ attractive.
S'il existe une cha\^\i ne de pseudo-disques de $\Delta_0$ \`a $\Delta_1$,
alors~:
\begin{itemize}
\item si le couple $(\Delta_0,\Delta_1)$ est indiff\'erent,
$\inpa(h,\Delta_0,\Delta_1)$ vaut $0$~;
\item  si le couple $(\Delta_0,\Delta_1)$ est attractif,
$\inpa(h,\Delta_0,\Delta_1)$ vaut $+1/2$ ou $-1/2$.
\end{itemize}
\end{prop}
La suite du texte n'utilisera que le cas o\`u le couple
$(\Delta_0,\Delta_1)$ est indiff\'erent.

\begin{demo}

\paragraph{Cas indiff\'erent}
Supposons qu'il existe une cha\^\i ne de pseudo-disques $(\ii{B_i})_{i=1..k}$ de $\Delta_0$ \`a
$\Delta_1$, et choisissons-la de façon \`a ce que $k$ soit minimal parmi toutes les cha\^\i nes de
$\Delta_0$ \`a $\Delta_1$.
 Comme dans la preuve du lemme de Franks (lemme
\ref{lem.fran}), l'id\'ee est de transformer la ``quasi-orbite'' fournie
par l'hypoth\`ese en une vraie orbite, \`a l'aide de modifications
libres~; le domaine de \hbox{Brouwer} engendr\'e par $\Delta_0$ rencontrera
alors $\Delta_1$, ce qui permettra de conclure \`a l'aide du
corollaire~\ref{cor.dotr}.
\index{modification libre}

On commence par montrer l'existence  d'une cha\^\i ne 
de disques ouverts pour $h$, encore not\'ee $(\ii{B_i})_{i=1..k}$, qui est une
cha\^{\i}ne de pseudo-disques de $\Delta_0$ \`a
$\Delta_1$ : pour cela on recopie les 
arguments de minimalit\'e de la preuve du lemme \ref{lem.fran} (les
d\'etails correspondant sont laiss\'es au lecteur). 

Dans un deuxi\`eme  temps, expliquons comment modifier cette cha\^{\i}ne
de disques ouverts en une cha\^{\i}ne dont tous les disques ouverts sont
inclus dans $O_1$.
Notons $(n_i)_{i=1..k-1}$ les temps de  transition de la cha\^{\i}ne et $(x_i)$ des
points de transition.
On peut supposer que l'entier
$k$ est minimal (parmi toutes les cha\^{\i}nes de disques ouverts de
$\Delta_0$ \`a $\Delta_1$). Tout d'abord, cette minimalit\'e entra\^ine clairement que
les disques $B_2, \dots B_{k-1}$ sont tous inclus dans $O_1$. 
La suite $(h^{-n_1}(B_2),
B_3, \dots, B_k)$ est encore une cha\^ine de disques ouverts (en particulier,
$h^{-n_1}(B_2)$ est disjoint des autres disques de la cha\^{\i}ne, sans quoi
 il existerait une cha\^{\i}ne de disques p\'eriodique, ce qui
contredirait le lemme de Franks~\ref{lem.fran}). D'autre part, $x_1$
appartient \`a $h^{-n_1}(B_2)$.
Si $x_1$ \'etait dans $\adhe(O_0)$, alors cette nouvelle cha\^{\i}ne de disques  contredirait la
minimalit\'e de $k$.
Ceci montre que $x_1$ appartient \`a $O_1$. Soit $B'_1$ la composante
connexe de $B_1 \cap O_1$ qui contient $x_1$.
Puisque $\adhe(B_1)$ rencontre
$\adhe(O_0)$,  l'adh\'erence de $B'_1$ rencontre
$\Delta_0$.
En remplaçant de m\^eme $B_k$ par un disque $B'_k$ plus petit, on
obtient une cha\^{\i}ne de disques ouverts $(B'_1, B'_2=B_2, \dots, B'_{k-1}=B_{k-1},
B'_k)$ qui sont tous inclus dans $O_1$.

Puisque $\adhe(B'_k)$ rencontre
$\adhe(O_2)$,  et que la droite  $\Delta_2$ est r\'epulsive,
 il existe  un point $x_k$ dans
 $B'_k$ tel que $h(x_k)$
soit dans $O_2$.
De m\^eme, il existe un point $x_0$ dans $O_0$ tel que $h(x_0)$ soit
dans  $B'_1$. On pose $n_0=1$.

On effectue maintenant $k$ modifications libres de $h$~: plus pr\'ecis\'ement,
on choisit  pour chaque entier $i$ entre $1$ et $k$ un hom\'eo\-morphisme
$\phi_i$, \`a support dans $\adhe(B'_i)$, et qui envoie $h^{n_{i-1}}(x_{i-1})$ sur $x_i$.
On consid\`ere alors l'hom\'eo\-morphisme $h_1=\phi \circ h$, avec
 $\phi=\phi_1 \circ \cdots \circ \phi_k$ (figure~\ref{fig29}).
Pour cet hom\'eo\-morphisme, le point $h(x_k)$ est un it\'er\'e du point $x_0$.

\begin{figure}[htpb]
\par
\centerline{\hbox{\input{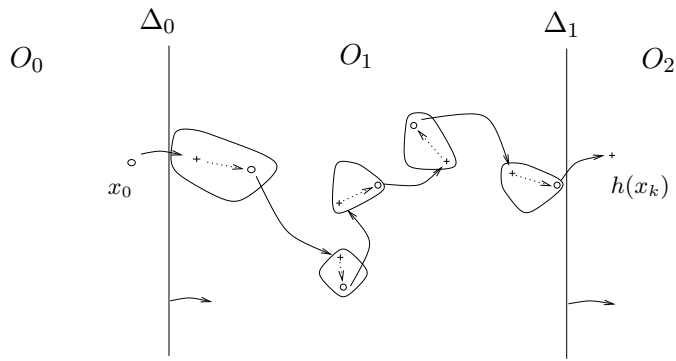}}}
\par
\caption{\label{fig29}Modifications libres pour obtenir une orbite qui
traverse la bande}
\end{figure}

Puisque les disques $B'_i$ sont tous inclus dans $O_1$, le support de
$\phi$ est inclus dans $\adhe(O_1)$. Par cons\'equent, on voit
facilement que  les droites $\Delta_0$ et
$\Delta_1$ sont aussi des droites de \hbox{Brouwer} pour l'hom\'eo\-morphisme
$h_1$. Ceci montre que l'indice partiel $\inpa(h_1,\Delta_0,\Delta_1)$ est bien d\'efini.
Comme il existe une orbite de $h_1$ qui rencontre \`a la fois $O_0$ et $O_2$, le
domaine de \hbox{Brouwer} engendr\'e par $\Delta_0$  (pour $h_1$) doit
rencontrer $\Delta_1$.
D'apr\`es le corollaire~\ref{cor.dotr}, on a alors $\inpa(h_1,\Delta_0,\Delta_1)=0$.
D'autre part, les indices partiels pour $\inpa(h_1,\Delta_0,\Delta_1)$
et $\inpa(h,\Delta_0,\Delta_1)$ sont \'egaux~: en effet, on peut
calculer ces deux indices dans une carte o\`u les deux droites sont
verticales (comme dans la d\'efinition de l'indice partiel),
 au moyen d'une courbe $\gamma$ disjointe du compact $K=\adhe(B'_1 \cup \cdots
\cup B'_k)$~; et puisque les deux hom\'eo\-mor\-phismes co\"\i ncident hors de
$K$, leurs indices le long de $\gamma$ sont \'egaux.

\paragraph{Cas attractif (sch\'ema de preuve)}
Les id\'ees sont les m\^emes que dans le cas indiff\'erent. La seule
diff\'erence significative appara\^{\i}t apr\`es avoir obtenu la
cha\^{\i}ne de disques ouverts inclus dans $O_1$~: cette fois-ci,
  il n'existe pas de
point de $B'_k$ dont l'image est dans $O_2$. On se contente donc
d'effectuer $k-1$ modifications libres, \`a support dans $B'_1, \dots, B'_{k-1}$.
On obtient ainsi un hom\'eo\-morphisme $h_1$
tel que, pour un certain $n$ positif, $h_1^n(\Delta_0)$ rencontre $B'_k$~;
d'autre part, $B'_k$ est libre pour $h_1$, il y a donc un arc libre
pour $h_1$  de
$h_1^n(\Delta_0)$ \`a $\Delta_1$. D'apr\`es le lemme~\ref{lem.arli},
l'indice  $\inpa(h_1,h_1^n(\Delta_0),\Delta_1)$ vaut $+1/2$ ou
$-1/2$. Le deuxi\`eme point du corollaire~\ref{cor.dotr} donne alors la
m\^eme valeur pour  $\inpa(h_1,\Delta_0,\Delta_1)$.
Comme dans le cas indiff\'erent, on a
$\inpa(h,\Delta_0,\Delta_1)=\inpa(h_1,\Delta_0,\Delta_1)$, ce qui
permet de conclure.
\end{demo}
\incn[Autre d\'emo : chirurgie...]

\section{D\'ecom\-po\-si\-tion en briques}
\label{sec.deco}
\index{d\'ecom\-po\-si\-tion}
\emph{Les r\'esultats de cette section sont tous adapt\'es de la th\`ese de
A. Sauzet (\cite{sauz1}) et de l'article de P. Le Calvez et
A. Sauzet \cite{leca3}}. 
 La possibilit\'e de se passer de la
condition de transversalit\'e pour la preuve du th\'eor\`eme de
translation plane est apparue lors d'une discussion avec
Lucien Guillou.

Dans toute cette section, \textbf{on suppose que $h$ est un hom\'eo\-mor\-phisme de
la sph\`ere\footnote{En r\'ealit\'e, la construction expliqu\'ee \`a la section~\ref{sec.code} est valable sur
n'importe quelle surface compacte sans bord.},
 pr\'eservant l'orientation, dont l'ensemble des points fixes
est fini}. Rappelons que les hom\'eo\-morphismes de \hbox{Brouwer} s'identifient
\`a des hom\'eo\-morphismes de la sph\`ere $\R^2 \cup \{\infty\}$ v\'erifiant
cette hypoth\`ese. 
 Nous rappelons d'abord la d\'efinition des
d\'ecom\-po\-si\-tions en briques et leur construction~; puis, \textbf{en supposant
l'absence de courbe d'indice $1$}, nous montrons l'existence d'une
droite de \hbox{Brouwer} associ\'ee \`a chaque brique, qui fournit une preuve
du th\'eo\-r\`eme de translation plane.

\subsection{Introduction aux d\'ecom\-po\-sitions en briques}
En premi\`ere approximation, une d\'ecom\-po\-si\-tion en briques est une sorte
de triangulation (localement finie)
du compl\'ementaire des points fixes de $h$ dans la sph\`ere. L'id\'ee principale
consiste \`a  choisir les
triangles assez petits pour qu'ils soient libres. Ceci permet
d'appliquer le lemme de Franks aux suites de triangles~: on pourrait
dire qu'on a ainsi obtenu une discr\'etisation de la dynamique qui garde
la propri\'et\'e de ne pas avoir d'orbite p\'eriodique.
Une autre propri\'et\'e importante des d\'ecom\-po\-si\-tions, appel\'ee minimalit\'e, est que
 les triangles ne sont pas trop petits, plus
pr\'ecis\'ement que la r\'eunion de deux triangles adjacents n'est jamais
libre.
Sous ces hypoth\`eses, un proc\'ed\'e dynamique simple (et automatique)
associe \`a chaque triangle $B$ un attracteur $A^+(B)$~: celui-ci
est simplement obtenu  en prenant le plus petit
attracteur, r\'eunion de triangles, qui contient l'image de $B$.
La minimalit\'e entra\^{\i}nera la connexit\'e de cet ensemble. On
montrera qu'il existe alors une unique 
composante connexe $\Delta(B)$ du bord de $A^+(B)$ qui rencontre $B$, et
que $\Delta(B)$ est une droite de \hbox{Brouwer}.

 Dans sa preuve du th\'eor\`eme de \hbox{Brouwer},  L. Guillou (\cite{guil1})
utilisait d\'ej\`a des triangulations libres.
D'autre part, M. Flucher avait combin\'e des triangulations libres avec le lemme de
Franks pour obtenir des points fixes d'hom\'eo\-mor\-phismes du tore, autour
de la conjecture d'Arnol'd (\cite{fluc1}). 
P. Le Calvez et A. Sauzet ont alors exploit\'e ces id\'ees pour obtenir ce
qui est sans doute la preuve la plus claire du th\'eor\`eme de \hbox{Brouwer} (\cite{leca3}).
Enfin, la propri\'et\'e de minimalit\'e est introduite dans la th\`ese
d'A. Sauzet (\cite{sauz1}), 
 qui y \'etudie les propri\'et\'es
combinatoires de ces triangulations de mani\`ere intensive,
afin de construire des courbes disjointes de leur
image pour les hom\'eo\-mor\-phismes du tore, de l'anneau ou de la sph\`ere.
On peut citer \'egalement l'utilisation de discr\'etisation de la
dynamique dans les preuves combinatoires du th\'eor\`eme de Conley-Zehnder
donn\'ees par S.~Alpern et V.~S.~Prasad (\cite{alpe1})~; dans ce
contexte, la tr\`es jolie utilisation
du r\'esultat combinatoire appel\'e ``th\'eor\`eme des mariages'' remonte \`a P. Lax (\cite{lax1}).

\subsection{D\'efinition}
\label{ss.definition-decomposition}
On pose $U=\S^2 \setminus \fixe(h)$. 
On dira qu'un sous-ensemble de $U$
est \emph{born\'e} si son adh\'erence dans $\S^2$ est disjointe de $\fixe(h)$. 
               
\begin{figure}[h!]
\par
\centerline{\hbox{\input{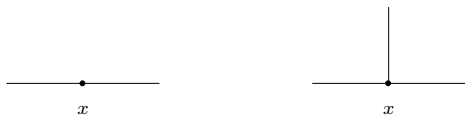}}}
\par
\caption{\label{fig31}Aspect local d'un graphe triadique}
\end{figure}

\begin{defi}
  \index{triadique}
Un \res{graphe triadique} est un sous-ensemble ferm\'e $F$ de $U$, tel
que, en tout point $x$ de $F$, $F$ est localement hom\'eomorphe \`a l'un
des deux dessins de la figure~\ref{fig31}. Si $F$ est localement
hom\'eomorphe au dessin de droite, on dit que
le point $x$ est un \res{sommet} de $F$.
\end{defi}

Soit $F$ un graphe triadique.
\begin{defi}
\label{def.briq}
\index{brique}
Une \res{brique} de $F$ 
est l'adh\'e\-ren\-ce (dans $U$) d'une composante connexe du
compl\'ementaire de $F$ dans $U$.
\end{defi}
\begin{defi}
\index{brique!adjacente}
Deux briques de $F$, distinctes, sont \res{adjacentes} si leur intersection n'est pas vide.
\end{defi}
\begin{defi}
\index{d\'ecom\-po\-si\-tion!ar\^ete}
Une \res{ar\^ete} de $F$ est l'adh\'erence (dans $U$) d'une
composante connexe du com\-pl\'e\-men\-taire dans $F$ des
sommets.
\end{defi}

\begin{defi}
Un graphe triadique $F$ est dit \res{compact} si aucun point fixe de
$h$ n'est dans l'adh\'erence d'une brique (autrement dit, toute brique est born\'ee).
\end{defi}
\begin{defi}
Le graphe triadique $F$ est dit \res{libre} (pour $h$) si
 toute brique de $F$ est libre.
\end{defi}
\begin{defi}\label{d.minimal}
Le graphe triadique $F$, libre, est dit \res{minimal} si pour tout
couple de briques adjacentes $B$ et $B'$, la r\'eunion $B \cup B'$
 n'est pas libre.
\end{defi}

\begin{defi}
Une \res{d\'ecom\-po\-si\-tion en briques} pour $h$ est un
graphe triadique, compact, libre, minimal.\footnote{Le vocabulaire
adopt\'e ici diff\`ere l\'eg\`erement de celui
de \cite{sauz1}~: notamment, on a remplac\'e ``brique ferm\'ee'' par
``brique'', ``born\'e'' par ``compact'', et ``maximal'' par ``minimal''.}
\end{defi}

Nous avons expliqu\'e plus haut les motivations des deux propri\'et\'es
principales, la liber\-t\'e et la minimalit\'e. Prendre un graphe triadique au
lieu d'une triangulation quelconque est une astuce topologique pour
que le bord de n'importe quelle r\'eunion de briques soit toujours une
sous-vari\'et\'e (\cf affirmation \ref{aff.cedr} ci-dessous).
 Les briques born\'ees sont une
commodit\'e, sans doute pas essentielle, qui \'evitera les \'etudes de
cas. 

\begin{rema}
Une d\'ecom\-po\-si\-tion en briques pour $h$ est aussi une d\'ecom\-po\-si\-tion
en briques pour $h^{-1}$.
\end{rema}

\begin{rema}
Dans la d\'efinition d'une d\'ecom\-po\-sition en briques, on ne suppose pas
que $F$ est connexe, ni qu'une ar\^ete contient deux sommets distincts,
et borde deux briques adjacentes distinctes. N\'eanmoins, en l'absence
de courbe d'indice~$1$, toutes ces propri\'et\'es d\'ecouleront des
propri\'et\'es de $h$ (voir le corollaire~\ref{cor.arete} ci-dessous).
\end{rema}

\begin{rema}\label{r.minimal}
 Soit $F$  un graphe triadique
libre. Si $F$ est minimal pour l'inclusion
(\emph{i.e.} s'il n'existe aucun graphe triadique
libre strictement inclus dans $F$),
 alors $F$  est minimal au sens
de la d\'efinition  \ref{d.minimal} ci-dessus.
\footnote{La r\'eciproque est \'egalement vraie si l'on suppose que toute
ar\^ete borde deux briques distinctes, ce qui est v\'erifi\'e en
l'absence de courbe d'indice~1 (voir le corollaire~\ref{cor.arete})~; nous
n'en aurons pas besoin.}
\end{rema}
\begin{demo}[de la remarque]
Sinon, il existe deux briques $B_1, B_2$ distinctes, adjacentes, dont la r\'eunion $B$
est libre. 
On v\'erifie facilement que $F'=F \setminus \inte(B)$ est un graphe
triadique libre, strictement inclus dans $F$.

\end{demo}

\subsection{Exemples}
\begin{figure}[htpb]
\par
\centerline{\hbox{\input{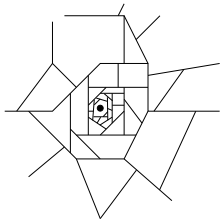}}}
\par
\caption{\label{fig36}Allure d'une d\'ecom\-po\-si\-tion au voisinage d'un
point fixe}
\end{figure}

La figure \ref{fig36} montre l'allure au voisinage d'un point fixe~:
remarquons que puisque les briques sont born\'ees, tout voisinage du
point fixe doit contenir une infinit\'e de briques.  Par contre, comme
$F$ est un graphe triadique, la famille des briques de $F$ est
localement finie au voisinage de tout point de $U$ (et la r\'eunion d'un
nombre quelconque de brique est toujours un ferm\'e de $U$).

\begin{figure}[hbtp]
\par
\centerline{\hbox{\input{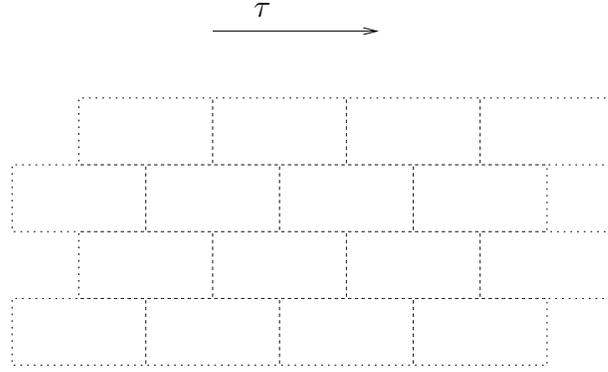}}}
\par
\caption{\label{fig10}Une d\'ecom\-po\-si\-tion pour la translation}
\end{figure}

\begin{figure}[hbtp]
\par
\centerline{\hbox{\input{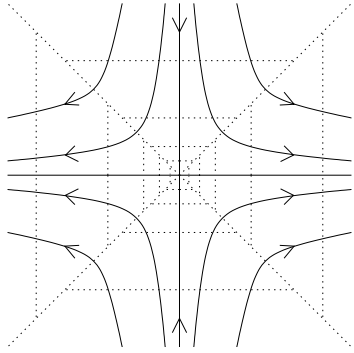}}}
\par
\caption{\label{fig11}Une d\'ecom\-po\-si\-tion pour le lin\'eaire hyperbolique}
\end{figure}

Les figures \ref{fig10} et \ref{fig11} montrent des d\'ecom\-po\-si\-tions pour
la translation $(x,y) \mapsto (x+1,y)$ et l'hom\'eo\-mor\-phisme
 lin\'eaire hyperbolique selle
$(x,y) \mapsto (2x,y/2)$ dans le mod\`ele du plan.


\subsection{Construction d'une d\'ecom\-po\-si\-tion}\label{sec.code}
\index{d\'ecom\-po\-si\-tion!construction}
\emph{La preuve du th\'eo\-r\`eme de cette section n'est pas essentielle \`a la compr\'ehension du
texte, et peut \^etre omise en premi\`ere lecture.}

\begin{theo}[A. Sauzet]
\label{the.exde}
Il existe une d\'ecom\-po\-si\-tion triadique, compacte, libre et  minimale pour $h$.
\end{theo}

Nous donnons une d\'emons\-tra\-tion fortement inspir\'ee de celle
d'A. Sauzet ; voici la principale modification : plut\^ot que de
construire d'abord une d\'ecom\-po\-sition dont les briques ne sont pas
compactes et qu'il faut re-d\'ecouper \emph{a posteriori}, on impose dans la
d\'ecom\-po\-sition la pr\'esence d'un certain nombre de briques, construites
\emph{a priori}, qui vont obliger la
d\'ecom\-po\-sition \`a \^etre compacte. Remarquons que ceci permet d'\'eviter
le recours au  lemme de Franks  (voir \cite{sauz1}, section~3.6)~; la construction
pr\'esent\'ee ici a donc l'avantage d'\^etre valable pour tout
hom\'eo\-morphisme de la sph\`ere (ou d'une autre surface compacte) 
avec un nombre fini de points fixes, sans
hypoth\`ese d'indice. 
    
\begin{demo}             
Tout d'abord, il est tr\`es facile de trouver un graphe triadique
 compact $F$ (par exemple, en utilisant le dual d'une triangulation
 localement finie de $U$). On peut tout
aussi facilement gagner la libert\'e : si une brique $B$ de $F$ n'est pas
libre, on la subdivise en briques assez petites (de diam\`etres inf\'erieurs
\`a $\epsilon=\mathrm{Inf}\{d(x,h(x) \mid x \in B \}$).
On voudrait ensuite  obtenir la minimalit\'e  en
choisissant  deux briques adjacentes dont la r\'eunion est libre, en
enlevant une ar\^ete dans leur intersection, et en recommen\c{c}ant tant que $F$
n'est pas minimale. En g\'en\'eral, malheureusement, ce proc\'ed\'e fait
perdre la compacit\'e en donnant
naissance \`a des briques non compactes (autrement dit, la famille des
d\'ecom\-po\-si\-tions triadiques libres $F'$ incluses dans $F$ est inductive
pour l'inclusion et on peut lui appliquer le lemme de Zorn, mais ceci est
faux en g\'en\'eral pour les d\'ecom\-po\-si\-tions triadiques, libres et compactes).

Le lemme suivant a pour but d'obtenir une propri\'et\'e suppl\'ementaire
qui emp\^eche l'apparition des briques non compactes~:
\begin{defi}[ (figure \ref{fig12})]
\label{def.anpd}
\index{anneau pr\'e-d\'ecompos\'e}
On appelle \res{anneau pr\'e-d\'ecompos\'e} tout
ensemble $C$ inclus dans $U$,  hom\'eomorphe \`a
l'anneau $\S^1 \times [0,1]$, muni d'un ensemble  $F_C$ appel\'e
\res{pr\'e-d\'ecom\-po\-sition de $C$} et v\'erifiant :
\begin{enumerate}
\item $F_C$ est un graphe triadique, inclus dans $C$ ;
\item $F_C$ contient le bord de $C$ ;
\item les adh\'erences des composantes connexes de $C \setminus F_C$
sont libres (on les appellera  \emph{briques de $F_C$})~;
\item la r\'eunion de deux briques de $F_C$ adjacentes n'est pas libre ;
\item{\bf propri\'et\'e suppl\'ementaire :} aucune brique ne rencontre
simultan\'ement les deux composantes connexes du bord de $C$.
\end{enumerate}
\end{defi}

\begin{figure}[hbtp]
  \par
\centerline{\hbox{\input{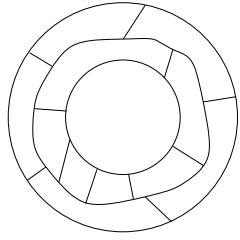}}}
\par
  \caption{\label{fig12}Topologie d'un anneau pr\'e-d\'ecompos\'e}
\end{figure}
\begin{lemm}\label{lem.pred}
Soit $J$ un cercle topologique dans $U$ tel que $h(J) \cap J \neq
\emptyset$. Alors tout voisinage de $J$ contient un anneau pr\'e-d\'ecompos\'e.
\end{lemm}

\paragraph{Utilisation du lemme}
A l'aide de ce lemme, on reprend toute la preuve du th\'eo\-r\`eme~\ref{the.exde}.
 Soit $x$ un point fixe de $h$ ; on trouve facilement
un cercle topologique qui entoure $x$, arbitrairement petit, et qui
rencontre son image.
Le lemme permet donc de trouver une suite $(C_i(x))_{i
\geq 0}$ de petits anneaux pr\'e-d\'ecompos\'es, entourant $x$, convergeant vers $x$
(pour la topologie de Hausdorff) et disjoints deux \`a deux (figure \ref{fig13}).
On effectue cette construction pour chacun des points fixes de $h$.
\begin{figure}[hbtp]
  \par
\centerline{\hbox{\input{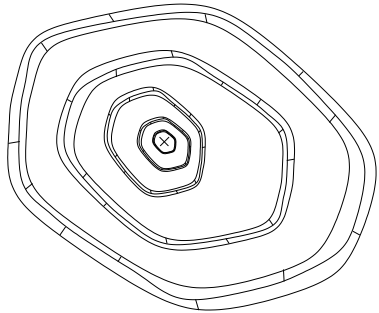}}}
\par
  \caption{\label{fig13}Suite d'anneaux pr\'e-d\'ecompos\'es autour d'un des
points fixes de $h$}
\end{figure}
On recommence alors la preuve expliqu\'ee plus haut : l'ensemble $\cup_{x,i}
F_{C_i(x)}$ est un graphe triadique, et il est facile de
l'\'etendre, \`a l'ext\'erieur de la r\'eunion des anneaux $C_i(x)$,
pour obtenir un graphe triadique et libre que l'on note $F$.
 On utilise le lemme de
Zorn pour choisir un graphe triadique libre 
$F'$, inclus dans $F$, qui est minimal pour l'inclusion.
D'apr\`es la remarque \ref{r.minimal}, $F'$ est \'egalement minimal au sens de la
d\'efinition \ref{d.minimal}.
On montre enfin que $F'$ est une d\'ecom\-po\-sition en briques
gr\^ace \`a l'affirmation suivante~:
\begin{affi}
Tout graphe triadique libre $F'$ inclus dans $F$ est compact.
\end{affi}
En effet, si $F'$ est inclus dans $F$, toute brique de $F'$ est
r\'eunion de briques de $F$. Une brique  de $F'$ non born\'ee devrait traverser de part
en part au moins un des
anneaux pr\'e-d\'ecompos\'es $C$ ; elle rencontrerait alors les deux bords de
cet anneau, donc contiendrait deux briques de $F_C$ adjacentes
(propri\'et\'e suppl\'ementaire de la d\'efinition \ref{def.anpd}),
 et ne pourrait pas \^etre libre
(propri\'et\'e 4 de la d\'efinition~\ref{def.anpd}). Ceci termine la
preuve du th\'eor\`eme~\ref{the.exde} \`a partir du lemme~\ref{lem.pred}.
\end{demo}

\begin{demo}[du lemme \ref{lem.pred}]
On part d'une ``d\'ecom\-po\-si\-tion libre minimale'' de  $J$, \ie une famille
finie $(J_0, \cdots, J_{k-1}), k \geq 2$ de sous-arcs de $J$  d'int\'erieurs disjoints
deux \`a deux, d'union $J$, libres, tels que la r\'eunion de deux arcs
adjacents ne soit pas libre  (son existence est imm\'ediate puisque $J
\cap h(J) \neq \emptyset$).
 \`a l'aide du th\'eor\`eme de Schoenflies, on \'epaissit ensuite
 $J$ en un anneau $C$ hom\'eomorphe \`a $J \times [0,1]$ suffisamment
fin pour que les ensembles $B_i \simeq J_i \times [0,1]$ soient encore libres
; l'ensemble $F_C$ r\'eunion des bords des $B_i$ v\'erifie alors les
points~1 \`a 4 de la d\'efinition \ref{def.anpd}.
 Par un argument classique de
transversalit\'e,
quitte \`a remplacer $C$ et les briques
$B_i$   par un anneau et des briques arbitrairement proches, on peut
renforcer le point~4 en supposant de plus qu'on a la propri\'et\'e~:

\emph{4'. pour tout couple $(B_i,B_{i+1})$ 
de briques adjacentes de $F_C$,
 l'image de {\bf l'int\'erieur} de l'une rencontre
l'autre.}

Nous allons maintenant  successivement d\'eformer et
re-d\'ecouper une \`a une les briques $B_i$ pour obtenir la propri\'et\'e
suppl\'ementaire de la d\'efinition \ref{def.anpd}.
Expliquons le traitement subi par $B_0$.

\paragraph{D\'eformation de la brique $B_0$}
On choisit une  brique adjacente \`a $B_0$, que l'on note  $B'$. 
D'apr\`es la propri\'et\'e~4', l'image par $h^\epsilon$ ($\epsilon=-1 \mbox{ ou } 1$)
de $B_0$ rencontre
l'int\'erieur de $B'$ ; on peut donc trouver un arc $\alpha_0$ inclus
dans $B'$,  d'int\'erieur inclus dans $B' \setminus h^\epsilon(B_0)$,
ayant une extr\'emit\'e sur l'ar\^ete entre $B_0$ et $B'$ et l'autre sur
$h^\epsilon(B_0) \cap \inte(B')$ (figure \ref{fig16}). On modifie
alors $B_0$ et $B'$ en pr\'elevant sur $B'$ un voisinage de $\alpha_0$
et en l'ajoutant \`a $B_0$ (par abus, on continuera \`a noter $B_0$ et
$B'$ les deux nouvelles briques). La nouvelle brique $B_0$  
cesse alors d'\^etre libre. Si  le voisinage
de $\alpha_0$ pr\'elev\'e sur $B'$ est suffisamment mince, on a encore 
la propri\'et\'e 4' entre $B'$ et les briques qui lui sont adjacentes (plus
pr\'ecis\'ement, pour chaque brique $B''$ adjacente \`a $B'$, 
on a $(h(B') \cup h^{-1}(B')) \cap
\inte(B'') \neq \emptyset$).
 On choisit aussi ce voisinage assez petit pour que
$B_0 \cap h^{-\epsilon}(B_0)$ ne rencontre pas les deux bords de l'anneau $C$.
\begin{figure}[hbtp]
  \par
\centerline{\hbox{\input{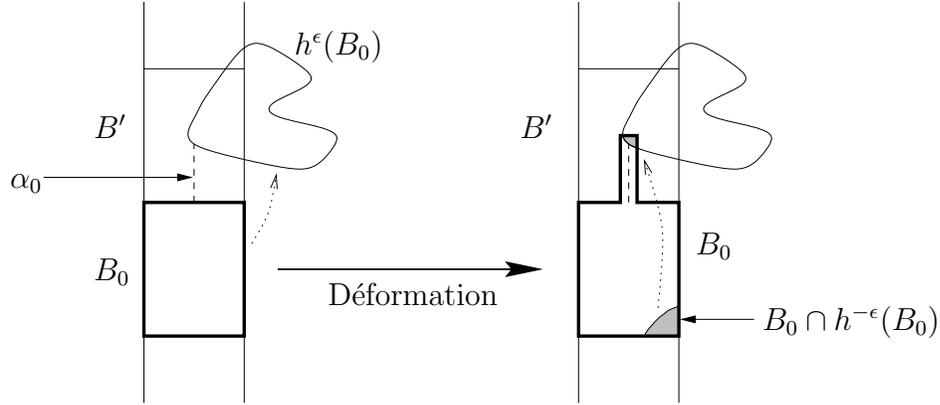}}}
\par
  \caption{\label{fig16}D\'eformation de la brique $B_0$}
\end{figure}

\paragraph{D\'ecoupage de la brique $B_0$}
Pour chaque brique $B$ adjacente \`a $B_0$, on choisit deux points
$x^1(B)$ et $x^2(B)$ dans $\inte(B_0) \cap h^{\epsilon(B)}(B)$ (o\`u
$\epsilon(B)=-1 \mbox{ ou } 1$). On d\'ecoupe alors $B_0$ en deux
briques $B_0^1$ et $B_0^2$ \`a
l'aide d'un arc $\beta_0$  de mani\`ere \`a ce que :
\begin{figure}[hbtp]
  \par
\centerline{\hbox{\input{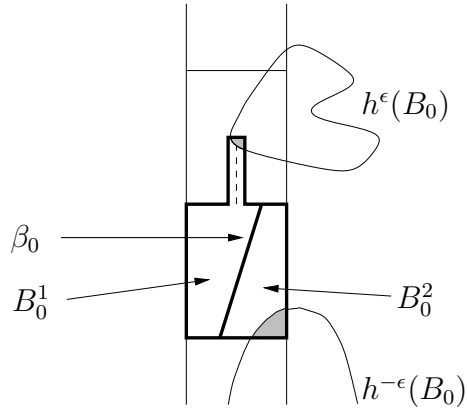}}}
\par
  \caption{\label{fig17}D\'ecoupage en deux de la brique $B_0$}
\end{figure}
\begin{enumerate}
\item ni  $B_0^1$ ni $B_0^2$ ne rencontre simultan\'ement les deux composantes
connexes de $\partial C$ ;
\item $h^\epsilon(B_0) \cap B_0 \subset B_0^1$ et $h^{-\epsilon}(B_0)
\cap B_0 \subset B_0^2$ ;  
\item pour chaque brique $B$ adjacente \`a $B_0$, $x^1(B) \in B_0^1$ et
$x^2(B) \in B_0^2$.
\end{enumerate}
Le deuxi\`eme point signifie que les deux nouvelles briques sont libres
; le dernier entra\^\i ne que la propri\'et\'e 4' est encore v\'erifi\'ee. On
obtient ainsi un nouveau graphe triadique $F_C^1$, r\'eunion des
fronti\`eres des briques apr\`es d\'eformation et d\'ecoupage.

\paragraph{D\'eformation et d\'ecoupage des autres briques}
A partir de $F_C^1$, on r\'eit\`ere la m\^eme op\'eration, successivement,  pour toutes les
autres briques $B_i$, en recopiant mot pour mot les \'etapes de
d\'eformation et de d\'ecoupage de $B_0$.

Le processus produit finalement un nouveau graphe triadique $F^k_C$ de $C$,
 qui d\'ecoupe $C$ en  $2k$
briques, dont aucune ne rencontre les deux bords de $C$ :  on a 
obtenu la propri\'et\'e suppl\'ementaire, ce qui montre que $F^k_C$ est
une pr\'e-d\'ecom\-po\-sition de $C$.
\end{demo}

\begin{defi}\label{def.pred}
On appelle \emph{pr\'e-d\'ecom\-po\-sition finie} un graphe triadique $F$ tel
que
\begin{enumerate}
\item $F$ est un ensemble compact de $U$ ;
\item si l'on appelle  \emph{pr\'e-briques}
les adh\'erences des composantes connexes du compl\'e\-men\-taire de $F$ ne
contenant pas de point fixe, les pr\'e-briques sont libres.
\end{enumerate}
\end{defi}

D'apr\`es la premi\`ere propri\'et\'e, les pr\'e-briques sont en nombre fini, et
leur union $K$ est compacte.
La construction des d\'ecom\-po\-sitions en briques  a pour corollaire
imm\'ediat~:
\begin{coro}\label{cor.pred}Soit $F$ une pr\'e-d\'ecom\-po\-si\-tion finie.
Il existe une d\'ecom\-po\-sition en briques $F'$ telle que $F' \cap K \subset
F \cap K$, autrement dit toute pr\'e-brique de $F$ est incluse dans une
brique de $F'$.
\end{coro}

\subsection{Premi\`eres propri\'et\'es}
\label{sub.trpl}
\textbf{On consid\`ere toujours un hom\'eo\-morphisme de la sph\`ere $\S^2$, pr\'eservant
l'orientation, n'ayant qu'un nombre fini de points fixes~; de plus, 
on suppose d\'esormais qu'il n'existe pas de courbe d'indice $1$
pour $h$.}
On se donne une d\'ecom\-po\-si\-tion en briques $F$ pour $h$ (th\'eo\-r\`eme~\ref{the.exde}).

\subsubsection*{$\bullet$ Lemme de Franks}
\begin{figure}[h!]
  \par
\centerline{\hbox{\input{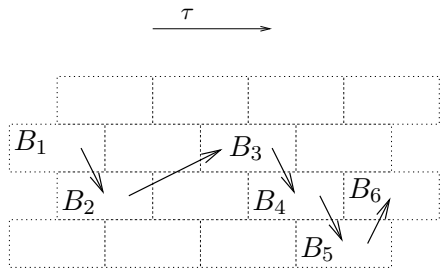}}}
\par
  \caption{\label{fig10b} Une cha\^ine de briques pour la translation}
\end{figure}

Nous commençons par appliquer le lemme de Franks aux briques de $F$.
\begin{defi}
Une \res{cha\^ine de briques} est une suite $B_1, \cdots, B_k$ de
briques de la d\'ecom\-po\-sition $F$ telles que pour tout $i=1, \cdots, k-1$,
$h(B_i)$ rencontre $B_{i+1}$ (figure \ref{fig10b}).
\end{defi}
Insistons sur le fait que, contrairement \`a la d\'efinition des
cha\^{\i}nes de disques, les briques d'une cha\^ine de briques ne sont
pas suppos\'ees \^etre deux \`a deux distinctes.

\begin{lemm}[discr\'etisation du lemme de Franks]
\label{lem.frdi}
Toutes les briques d'une cha\^ine de briques sont distinctes. En
particulier, toute cha\^\i ne de briques est une
cha\^\i ne de pseudo-disques, et il n'existe pas de cha\^ine de briques qui
soit une cha\^ine de pseudo-disques p\'eriodique.
\end{lemm}

\begin{demo}
Remarquons d'abord qu'une brique v\'erifie clairement la d\'efinition d'un
pseudo-disque (d\'efinition~\ref{def.psdi}).
Soit $(B_1, \cdots, B_k)$ une cha\^\i ne de briques, et supposons qu'il
existe deux entiers $i<j$ tels que $B_i=B_j$. Il est clair que, quitte
\`a extraire, on
peut supposer de plus que
 les briques $B_i, \cdots, B_{j-1}$ sont deux \`a
deux distinctes. Elles forment alors une cha\^ine de pseudo-disques
p\'eriodique, ce qui contredit le lemme de Franks (lemme~\ref{lem.fran}).
\end{demo}

\subsubsection*{$\bullet$ \'Enonc\'es des propri\'et\'es}
On note $\partial_U$ la fronti\`ere topologique dans $U$.
\begin{affi}
\label{aff.cedr} 
Si $V$ est une r\'eunion quelconque de briques,
alors $\partial_U V$ est une sous-vari\'et\'e sans bord de $U$,
ferm\'ee dans $U$, donc
une r\'eunion (finie ou non) de cercles topologiques et de droites topologiques.
\end{affi}

\begin{affi}
\label{aff.lune}
 Soit $B_1$ et $B_2$ deux briques adjacentes ; alors de deux
choses l'une~: 
\begin{enumerate}
\item soit $h(B_1)$ rencontre $B_2$,
\item soit $h(B_2)$ rencontre $B_1$.
\end{enumerate}
\end{affi}

\begin{affi}
\label{aff.brdi}
Les briques sont des disques topologiques ferm\'es.
\end{affi}

\begin{coro}~
\label{cor.arete}
\begin{enumerate}
\item Toute ar\^ete est incluse dans deux briques distinctes~;
\item toute ar\^ete  contient deux sommets distincts~;
\item $F$ est connexe.
\end{enumerate}
\end{coro}

\begin{affi}
\label{aff.brad}
  L'image par $h$  de toute brique $B$ rencontre au moins
  une brique adjacente \`a $B$. M\^eme chose pour l'image par $h^{-1}$.
\end{affi}

\subsubsection*{$\bullet$ Preuves des propri\'et\'es}
La preuve de l'affirmation \ref{aff.cedr} est imm\'ediate.

\begin{demo}[de l'affirmation \ref{aff.lune}]
Puisque $B_1 \cup B_2$ n'est pas libre, mais que $B_1$ et $B_2$ le
sont, c'est que  $h(B_1)$ rencontre $B_2$ ou $h(B_2)$ rencontre
$B_1$. D'autre part, si ces deux possibilit\'es \'etaient simultan\'ement
r\'ealis\'ees, alors $(B_1,B_2,B_1)$ serait une cha\^\i ne de briques qui
contredirait le lemme \ref{lem.frdi}.
\end{demo}

\begin{demo}[de l'affirmation \ref{aff.brdi}]
Soit $B$ une brique. D'apr\`es l'affirmation~\ref{aff.cedr}, son bord est
une r\'eunion de cercles et de droites topologiques ; comme $B$ est
compacte (par d\'efinition des d\'ecom\-po\-si\-tions en briques), c'est en fait une
r\'eunion de cercles. Il suffit donc de montrer que le bord de $B$ est
connexe,  ou encore que le compl\'ementaire de $B$ dans
la sph\`ere est connexe.

Comme $\partial B$ est une r\'eunion de cercles topologiques, chaque
 composante connexe de $\S^2 \setminus \inte(B)$ est un disque
 topologique ferm\'e. Soit $D_1$ l'un d'entre eux ; c'est aussi une
 r\'eunion de briques de la d\'ecom\-po\-sition.

Soit $B_1=B$, et $B_2$ une brique de $D_1$ adjacente \`a $B$.
En appliquant l'affirmation \ref{aff.lune},
 on voit que $D_1$ rencontre $h(B)$ ou
$h^{-1}(B)$. Supposons par exemple que $D_1$ rencontre $h(B)$ : comme
$B$ est libre, $h(B)$
 est  un ensemble connexe inclus dans le compl\'ementaire de
$B$, il est donc inclus dans $D_1$. Comme $\partial D_1 \subset B$, 
on a alors $h(\partial D_1) \subset D_1$.

Ceci nous donne deux possibilit\'es : ou bien $h(D_1) \subset D_1$, ou
bien $h(\adhe(\S^2 \setminus D_1)) \subset D_1$. Dans le premier cas, $D_1$
serait un disque attractif,\index{disque attractif}
 mais c'est exclu, car le lemme \ref{lem.diat} entra\^inerait
l'existence d'une courbe d'indice $1$, contrairement aux hypoth\`eses. On
est par cons\'equent dans le deuxi\`eme cas, et l'ensemble $\adhe(\S^2
\setminus D_1)$
est libre. Par minimalit\'e de la d\'ecom\-po\-si\-tion, cet ensemble \'etant une
r\'eunion de briques, il est en fait constitu\'e d'une seule brique, qui
est n\'ecessairement la brique $B$. Ceci prouve que le compl\'ementaire de
$B$ est le disque $D_1$, qui est connexe.
\end{demo}

\begin{demo}[du corollaire \ref{cor.arete}]

\paragraph{Point 1 (figure~\ref{fig46})} 
Soient $B$ une brique et $\alpha$ une ar\^ete incluse dans $B$, qui n'est
incluse dans
aucune autre brique. Nous allons montrer que sous ces hypoth\`eses, $B$
ne peut pas \^etre un disque topologique.

L'int\'erieur de $\alpha$ est inclus
dans l'int\'erieur de $B$ (sans quoi $\alpha$ borderait une autre brique). 
Par d\'efinition des briques, l'ensemble $O=B \setminus F$ est un ouvert connexe.
 Il existe donc  une courbe de Jordan $\gamma$ incluse dans
$O \cup \alpha$, rencontrant $\alpha$ en unique point, et de mani\`ere
 transverse (voir figure~\ref{fig46}). 

Soient 
$D_1$ et $D_2$ les deux disques topologiques d\'elimit\'es par $\gamma$.
Montrons que $D_1$ rencontre le compl\'ementaire de $B$. 
 Ou bien $D_1$ contient un point fixe de $h$, qui est dans le
compl\'ementaire de $B$ puisque $B$ est born\'ee (compacit\'e des
d\'ecom\-po\-si\-tions). Ou bien $D_1$ ne contient
aucun point fixe~: dans ce cas, l'ensemble $G=(F \cap D_1) \setminus
 \inte(\alpha)$ est un graphe triadique qui n'a qu'un nombre fini de
 sommets. Comme tout graphe fini sans feuille poss\`ede un cycle, 
on peut trouver une courbe de Jordan $\gamma'$
incluse dans $G$, donc dans $F \cap D_1$. L'un des deux disques
topologiques ouverts bord\'es par $\gamma'$  est alors inclus dans 
le compl\'ementaire de $B$.

En appliquant le m\^eme raisonnement \`a $D_2$, on voit que le
compl\'ementaire de $B$ a au moins deux composantes connexes, ce qui contredit
l'affirmation~\ref{aff.brdi}. 
\begin{figure}[htpb]
\par
\centerline{\hbox{\input{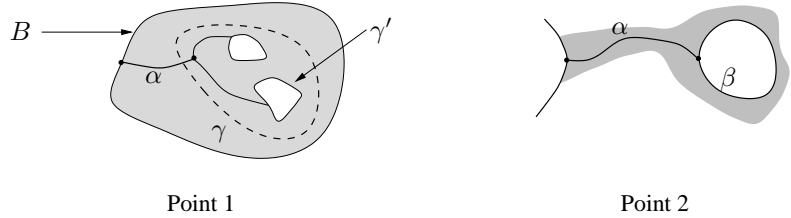}}}
\par
\caption{\label{fig46}Preuve du corollaire~\ref{cor.arete}}
\end{figure}

\paragraph{Point 2 (figure~\ref{fig46})} 
Si une ar\^ete $\beta$ contenait un unique sommet, alors l'autre
ar\^ete $\alpha$ adjacente \`a ce sommet serait contenue dans une unique
brique, contredisant le point~1.

\paragraph{Point 3}
Commençons par remarquer que tout point de $F$ appartient \`a la fronti\`ere
 d'au moins une brique~: ceci est une cons\'equence du point~1.

 Les briques sont en nombre d\'enombrable ; on peut les num\'eroter de
mani\`ere \`a ce que chaque brique $B_i$ soit adjacente \`a une  brique $B_j$ avec $j < i$.
Pour tout entier $n$, la r\'eunion $\partial B_0 \cup \cdots \cup
\partial B_n$ est alors
connexe puisque chaque bord $\partial B_i$ est connexe. Or d'apr\`es la
remarque initiale, $F=\cup
\partial B_i$, c'est donc une r\'eunion croissante d'ensembles connexes,
il est connexe.
\end{demo}

\begin{demo}[de l'affirmation \ref{aff.brad} (figure~\ref{fig19})]
  On raisonne par l'absurde. Soit $\cal E$ l'ensemble des briques
  adjacentes \`a $B$. Si $h(B)$ ne rencontre aucune des briques de $\cal
  E$, c'est que $h^{-1}(B)$ les rencontre toutes (d'apr\`es
  l'affirmation~\ref{aff.lune}).
 Soit $B_0 \in {\cal E}$ ; comme $B$ est libre, $h^{-1}(B)$ ne rencontre
  pas $B_0 \cap B$, en particulier $B_0$ n'est pas contenu dans
  $h^{-1}(B)$. La brique $B_0$ rencontre donc \`a la fois $h^{-1}(B)$ et
  son compl\'ementaire, elle rencontre donc la fronti\`ere $h^{-1}(\partial B)$.
 Il existe donc une  brique $B_1$, adjacente \`a $B$,
  telle que $h^{-1}(B_1) \cap B_0 \neq \emptyset$.
\begin{figure}[htpb]
\par
\centerline{\hbox{\input{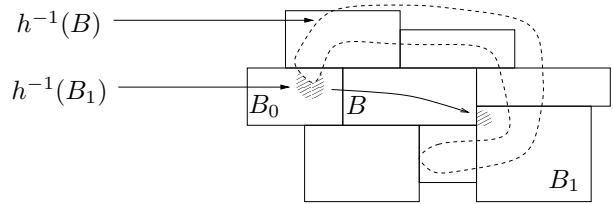}}}
\par
\caption{\label{fig19}Preuve de l'affirmation \ref{aff.brad}}
\end{figure}
  En r\'eit\'erant le processus, on trouve une suite $B_0, B_1, \cdots$
  de briques de $\cal E$ telles que $h(B_i) \cap B_{i+1} \neq
  \emptyset$. D'apr\`es la discr\'etisation du lemme de Franks (lemme \ref{lem.frdi}),
 les briques $(B_i)_{i \geq
  0}$ sont toutes distinctes ; mais c'est impossible puisqu'il n'y a qu'un
  nombre fini de briques adjacentes \`a $B$.
\end{demo}

\subsection{Construction des attracteurs $A^+(B)$}
\begin{figure}[p]
\par
\centerline{\hbox{\input{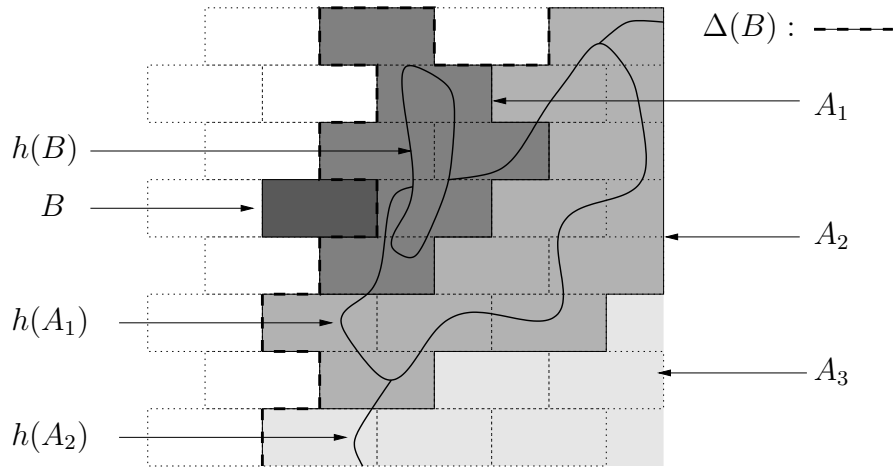}}}
\par
\caption{\label{fig32}Construction de l'attracteur $A^+(B)$ : les
premi\`eres \'etapes}
\end{figure}

Soit $B$ une brique de la d\'ecom\-po\-si\-tion.
On pose $A_0(B)=B$, et on d\'efinit par r\'ecurrence, pour tout entier
positif $i$,
$$
A_{i+1}(B)=\bigcup\{B' \mbox{ brique de } F \mid B' \cap h(A_i) \neq
\emptyset\}.
$$
On pose alors :
\begin{defi}[ (figure \ref{fig32})]
\index{$A^+(B)$}
\index{brique!attracteur}
$$
A^+(B)=\bigcup_{i \geq 1}A_i.
$$
\end{defi}

On montre facilement :
\begin{affi}
\label{aff.chai}
$A^+(B)$ est l'union des briques $B'$ telles qu'il existe une cha\^ine
de briques $(B_1=B, \cdots, B_k=B')$ avec $k\geq 2$.
C'est aussi le plus petit ensemble dont l'int\'erieur contient $h(B)$,
qui est une
r\'eunion de briques de $F$, et est un attracteur strict.
\end{affi}

\begin{figure}[p]
\par
\centerline{\hbox{\input{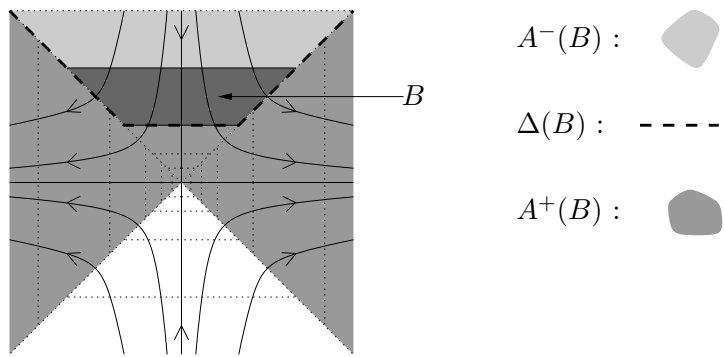}}}
\par
\caption{\label{fig18}Les ensembles $A^+(B)$ et $A^-(B)$ pour
l'hom\'eo\-mor\-phisme 
lin\'eaire hyperbolique selle (la d\'ecom\-po\-sition est celle de la figure~\ref{fig11})}
\end{figure}

On d\'efinit $A^-(B)$ de mani\`ere sym\'etrique~:
\begin{defi}\label{def.repu}\index{$A^-(B)$}
$A^-(B)$ est l'union des briques $B'$ telles qu'il existe une cha\^ine
de briques $(B_1=B', \cdots, B_k=B)$ avec $k\geq 2$.
\end{defi}
 La figure
\ref{fig32} illustre la construction de  l'ensemble $A^+(B)$, la
figure \ref{fig18} montre $A^+(B)$ et $A^-(B)$ pour une brique
d'une d\'ecom\-po\-si\-tion pour l'hom\'eo\-mor\-phisme
 lin\'eaire hyperbolique selle.

\begin{figure}[p]
\par
\centerline{\hbox{\input{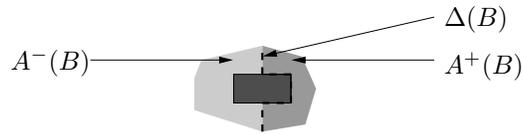}}}
\par
\caption{\label{fig38}Topologie des  ensembles $A^-(B)$,
 $A^+(B)$ et $\Delta(B)$ au voisinage d'une brique}
\end{figure}

\begin{affi}\label{aff.disj}
   L'attracteur $A^+(B)$ est connexe et rencontre $B$. Les int\'erieurs des ensembles $B$, $A^+(B)$ et $A^-(B)$ sont deux \`a deux
  disjoints. 
\end{affi}
\begin{affi}
\label{aff.allu}
 Toute brique adjacente \`a $B$ est soit dans $A^+(B)$, soit dans
$A^-(B)$, et $A^+(B) \cap B$ et $A^-(B) \cap B$ sont connexes et non vides.
\end{affi}
L'allure au voisinage d'une brique est donc celle de la figure \ref{fig38}.
\begin{affi}
\label{aff.drto}
Soit $A$ un attracteur strict, connexe,  qui est une r\'eunion de
briques. Alors  
les composantes connexes de la fronti\`ere de $A$ dans $U$ sont des droites topologiques.
\end{affi}

\begin{demo}[de l'affirmation \ref{aff.disj}]
Comme $B$ est connexe, la r\'eunion $A_1$ des briques qui rencontrent $h(B)$
est  connexe. De plus, elle rencontre $B$ car $h(B)$ rencontre au
moins une brique adjacente \`a $B$ (affirmation \ref{aff.brad}). 
La connexit\'e de $A^+(B)$ s'en suit. Comme $A^+(B)$ contient $A_1$, il
rencontre $B$.

La deuxi\`eme phrase de l'affirmation est une cons\'equence imm\'ediate de
la version discr\'etis\'ee du lemme de Franks (lemme \ref{lem.frdi}).
\end{demo}

\begin{demo}[de l'affirmation \ref{aff.allu}]
D'apr\`es l'affirmation~\ref{aff.lune}, toute brique adjacente \`a $B$ est
dans $A^+(B)$ ou dans $A^-(B)$ ; mais pas dans les deux d'apr\`es
l'affirmation pr\'ec\'edente. Celle-ci dit \'egalement que  les ensembles  $A^+(B) \cap B$
et $A^-(B) \cap B$ ne sont pas vides. Il reste \`a montrer leur connexit\'e.

Supposons que l'un de ces deux ensembles ne soit pas connexe. Alors on
peut trouver quatre briques $B_1$, $B_2$, $B_3$, $B_4$, adjacentes \`a $B$,
telle que $B_1$ et $B_3$ soient dans $A^+(B)$ et $B_2$ et $B_4$ dans
$A^-(B)$, et telles qu'en parcourant $\partial B$ on rencontre
successivement $B_1$, $B_2$, $B_3$ et $B_4$ (figure \ref{fig_allure_brique}).
\begin{figure}[hbtp]
  \par
\centerline{\hbox{\input{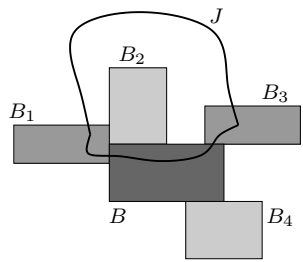}}}
\par \caption{\label{fig_allure_brique} Preuve de l'affirmation \ref{aff.allu}}
\end{figure}

Comme $A^+(B)$ est connexe et r\'eunion de briques, son int\'erieur est
connexe par arcs, et il existe un arc $\gamma \subset \inte(A^+(B)) $
joignant un point de $B_1$ \`a un point de $B_3$. A l'aide d'un autre
arc $\gamma' \subset B \cup B_1 \cup B_3$, on peut prolonger $\gamma$
en une courbe de Jordan $J$. Cette courbe de Jordan ne rencontre pas
$A^-(B)$ et s\'epare $B_2$ et $B_4$, ce qui contredit la connexit\'e de $A^-(B)$.
\end{demo}

\begin{demo}[de l'affirmation \ref{aff.drto}]
Soit $C$ l'une de ces  composantes connexes.
D'apr\`es l'affirmation \ref{aff.cedr}, il suffit de prouver que $C$
n'est pas un cercle topologique. Supposons le contraire ; comme $A$
est connexe, l'un  des deux disques topologiques ouverts bord\'es par $C$  est disjoint de
$A$, on le note $D$. Comme $A$ est un attracteur strict, on a $h(C) \cap
D=\emptyset$, donc de deux choses l'une : ou bien $D \subset h(D)$,
 mais le lemme \ref{lem.diat} interdit les disques attractifs~;\index{disque attractif}
 ou bien l'adh\'erence de $D$ est libre. Dans ce cas, 
l'adh\'erence de $D$ est exactement l'une des briques de la
d\'ecom\-po\-si\-tion (minimalit\'e de la d\'ecom\-po\-sition). Par ailleurs,
 les briques adjacentes \`a $D$ sont toutes
dans $A$ ; et l'une d'entre elles au moins, notons-la $B_1$, doit
rencontrer $h^{-1}(D)$ (affirmation \ref{aff.brad}). On a alors
$h(B_1) \cap D \neq \emptyset$, ceci contredit le fait que 
$A$ est un attracteur.
\end{demo}
\subsection{Construction des droites de \hbox{Brouwer} $\Delta(B)$}

D'apr\`es l'affirmation \ref{aff.allu}, il existe une unique composante
connexe de $\partial A^+(B)$ qui rencontre $B$.
\begin{defi}[ (figures \ref{fig32} et \ref{fig18})]
\index{brique!droite de \hbox{Brouwer}}
\index{$\Delta(B)$}
On note $\Delta(B)$ l'unique composante connexe de la fronti\`ere de
$A^+(B)$ dans $U$  qui rencontre $B$.
\end{defi}

\begin{prop}\label{pro.drbr}
L'ensemble $\Delta(B)$ est une droite de \hbox{Brouwer}.
L'en\-sem\-ble $B \cup \Delta(B)$ est libre, et le
domaine de \hbox{Brouwer} engendr\'e par $\Delta(B)$ contient $B$.
\end{prop}

\begin{demo}
\paragraph{$\Delta(B)$ est libre}
C'est clair puisque $A^+(B)$ est un attracteur strict.

\paragraph{$\Delta(B)$ est une droite topologique}
C'est une cons\'equence de  l'affirmation \ref{aff.drto}.

\paragraph{$\Delta(B)$ est une droite de \hbox{Brouwer}}
Si les deux extr\'emit\'es de $\Delta$ aboutissent au m\^eme point fixe
$x$ de $h$, il reste \`a montrer que $\Delta$ s\'epare 
$h(\Delta)$ et $h^{-1}(\Delta)$. Soit $D$ celui des deux disques
topologiques ferm\'es  bord\'es
par $\Delta(B) \cup \{x\}$ qui contient $A^+(B)$, et $D'$ son
compl\'ementaire.  Puisque $h(A^+(B))
\subset A^+(B)$, on a $h(\partial D) \subset D$. Si  $h(D')
\subset D$, alors $D'$ est libre, ce qui contredit la
minimalit\'e de la d\'ecom\-po\-si\-tion car $D'$
contient plusieurs briques (par exemple, $A^-(B) \subset D'$). On a
donc $h(D) \subset D$.

On en d\'eduit  $h(D) \subset D \subset h^{-1}(D)$, 
 et la propri\'et\'e de s\'eparation est v\'erifi\'ee.

\paragraph{Fin de la preuve}
Par d\'efinition de $A^+(B)$, $h(B) \subset \inte(A^+(B))$, donc $h(B)
\cap \Delta(B)=\emptyset$. D'autre part $h(\Delta(B)) \subset
\inte(A^+(B))$, donc $h(\Delta(B)) \cap B=\emptyset$ d'apr\`es
l'affirmation \ref{aff.disj}.  Comme $B$ et $\Delta(B)$ sont libres,
ceci montre que $B\cup \Delta(B)$ est libre.

Comme $B$ est connexe et rencontre $\Delta(B)$, on en d\'eduit que $B$
est dans le domaine de \hbox{Brouwer} engendr\'e par $\Delta(B)$ (plus
pr\'ecis\'ement, $h(B)$ est dans le domaine fondamental de $\Delta(B)$, 
$D(\Delta(B),h(\Delta(B)))$).
\end{demo}

\subsection{Preuve du th\'eor\`eme de \hbox{Brouwer}, corollaires}
\label{sec.brou}
\index{translation plane!preuve}
\begin{theo}[\hbox{Brouwer}]
Tout point de $\S^2 \setminus \fixe(h)$ est dans un domaine de \hbox{Brouwer}.
\end{theo}

\begin{demo}[(P. Le Calvez, A. Sauzet)]
Tout point de $U$ est dans une brique $B$, et toute brique $B$ est dans le
domaine de \hbox{Brouwer} engendr\'e par $\Delta(B)$ (proposition \ref{pro.drbr}).
\end{demo}

Le r\'esultat suivant nous servira \`a la section \ref{sec.reca} :
\begin{coro}[de la preuve du th\'eor\`eme de \hbox{Brouwer}]~
\index{arc libre!et domaine de \hbox{Brouwer}}\label{cor.precedent} \label{cor.arli}
\begin{enumerate}
\item Tout arc libre est contenu dans un
domaine de \hbox{Brouwer}~;
\item S'il existe un arc libre $\gamma$ de $x$ \`a $y$, alors il existe un arc
  libre de $x$ \`a $h^{n}(y)$ pour tout entier $n$.
\end{enumerate}
\end{coro}
\pagebreak
\begin{demo}
\paragraph{Premier point}
Un arc libre est inclus dans un disque topologique ferm\'e libre $D$~; le
bord d'un tel disque forme une pr\'e-d\'ecom\-po\-sition finie
(d\'efinition~\ref{def.pred}).
 D'apr\`es le
corollaire~\ref{cor.pred}, on peut choisir une d\'ecom\-po\-sition en briques
dont une brique contient $D$. Comme
  toute brique est incluse dans un domaine de \hbox{Brouwer}, on obtient
  le r\'esultat.
\paragraph{Deuxi\`eme point}

Remarquons tout d'abord que l'arc $\gamma$ \'etant libre, il est
disjoint de tous ses it\'er\'es par $h$ (lemme~\ref{lem.libr}). En particulier,
$x$ et $y$ ne sont pas dans la m\^eme orbite de $h$.

\emph{Examinons d'abord le cas o\`u $h$ est une translation affine du plan.} Soit
$p$ l'application quotient du plan dans l'anneau ouvert
$\R^2/h$. L'hypoth\`ese entra\^ine donc que $p(x) \neq p(y)$. On construit
alors facilement, dans l'anneau, un arc qui relie $p(x)$ \`a $p(y)$ et qui
fait ``$n$ tours de plus de $p(\gamma$)''. Cet arc se rel\`eve
en un arc $\gamma'$ joignant $x$ \`a $h^n(y)$~; comme $p(\gamma')$ est un
arc, $\gamma'$ est libre.

\emph{Revenons au cas g\'en\'eral.}  D'apr\`es le premier point,
 un arc libre est inclus dans un domaine de
  \hbox{Brouwer}. Mais puisque la restriction de $h$ \`a un tel domaine est
  conjugu\'e \`a une translation affine du plan (affirmation
 \ref{aff.dofo}, second point),
 le cas pr\'ec\'edent permet de conclure.
\end{demo}

\clearpage
\part{Preuve du th\'eo\-r\`eme principal de dynamique globale}\label{par.preu}
Dans cette partie, nous prouvons le th\'eo\-r\`eme principal de dynamique
global (\ref{the.prin}).
Rappelons que l'id\'ee de la preuve du th\'eo\-r\`eme~\ref{the.prin} est
d\'ecrite \`a la section~\ref{sub.idee}. La
section~\ref{sec.enon} contient les \'enonc\'es des lemmes et propositions
utilis\'es dans cette preuve~;
la preuve proprement dite occupe la section~\ref{sec.preu}.
 Les lemmes sont ensuite d\'emontr\'es dans les
sections~\ref{sec.cam} \`a \ref{sec.reca}.
\emph{Les sections~\ref{sec.preu} \`a \ref{sec.reca} sont largement
ind\'ependantes les unes des autres, et peuvent \^etre lues dans un ordre quelconque (la
seule exception notable concerne la section~\ref{sec.dbns}, qui utilise les
notions introduites au d\'ebut de la section~\ref{sec.cam}).}

\section{\'Enonc\'es des lemmes}\label{sec.enon}
Dans cette section, on \'enonce les lemmes et propositions qui seront
utilis\'es dans la preuve du th\'eo\-r\`eme principal~; leurs d\'emons\-tra\-tions
font l'objet des sections suivantes.
En~\ref{sub.cphb}, on transpose les d\'efinitions des croissants et des
p\'etales, introduites pr\'ec\'edemment sous l'hypoth\`ese (H2) du th\'eo\-r\`eme principal,
 au cadre des hom\'eo\-morphismes de \hbox{Brouwer} (hypoth\`ese (H1)). Ceci permet,
en~\ref{sub.tphb}, d'\'ecrire une version du th\'eo\-r\`eme principal dans ce
cadre (th\'eor\`eme~\ref{the.prinbrou}). Il reste alors \`a \'enoncer les r\'esultats permettant le
changement de cadre~: en~\ref{sub.dbns},
l'existence d'une droite de \hbox{Brouwer} reliant les deux points fixes~;
en~\ref{sub.inct} et \ref{sub.rele}, la possibilit\'e de relever canoniquement
la dynamique sur $\S^2 \setminus\{N,S\}$ en un hom\'eo\-morphisme de \hbox{Brouwer}.

\subsection{Croissants et p\'etales pour les hom\'eo\-morphismes de \hbox{Brouwer}}\label{sub.cphb}
\begin{figure}[htp]
\par
\centerline{\hbox{\input{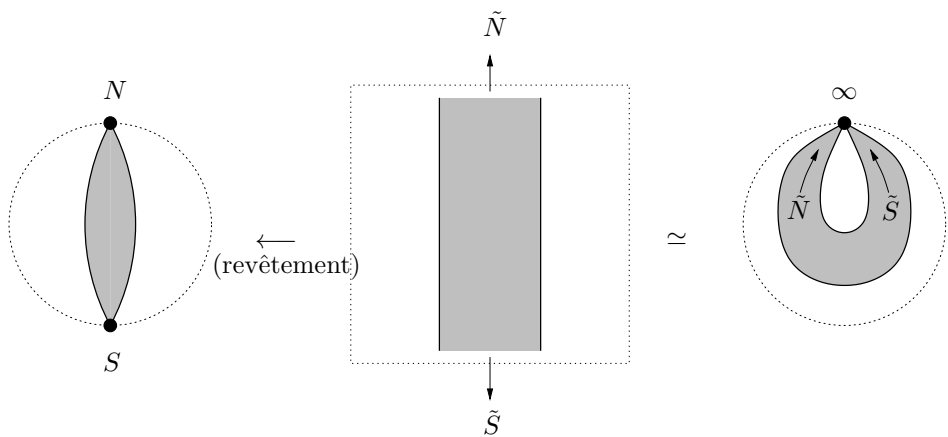}}}
\par
\caption{\label{fig47}Les deux notions de croissants~: pour un
hom\'eo\-morphisme de la sph\`ere fixant $N$ et $S$ (\`a gauche), ou pour un
hom\'eo\-morphisme de \hbox{Brouwer}, vu comme hom\'eo\-morphisme du  plan sans point
fixe (au milieu) ou comme hom\'eo\-morphisme de la sph\`ere fixant le point
\`a l'infini (\`a droite)}
\end{figure}

Nous commençons par d\'efinir p\'etales et croissants pour un
hom\'eo\-morphisme de \hbox{Brouwer} $H$, c'est-\`a-dire sous l'\textbf{hypoth\`ese (H1)}. 
La d\'efinition des p\'etales ne posera pas de probl\`eme si on assimile $H$
 \`a un hom\'eo\-morphisme de la sph\`ere fixant un
unique point. La d\'efinition des croissants est un peu plus
compliqu\'ee~: pour la comprendre, il peut \^etre utile de garder \`a
l'esprit que lorsqu'on appliquera la d\'efinition, $H$ sera le relev\'e
d'un hom\'eo\-morphisme $h$ de l'anneau $\S^2 \setminus\{N,S\}$, et que
les croissants de $H$ (pour la nouvelle d\'efinition) 
doivent correspondre \`a des croissants de $h$
(pour la d\'efinition donn\'ee \`a la section~\ref{sec.deen}). Il faudra
ainsi penser aux croissants de $H$ comme ayant deux extr\'emit\'es (bien
que, dans la sph\`ere $\R^2 \cup \{\infty\}$, ces extr\'emit\'es
aboutissent toutes les deux au point fixe $\infty$). Ceci est illustr\'e
sur la figure~\ref{fig47}.

\subsubsection*{$\bullet$ P\'etales}
Soit $H$ un hom\'eo\-morphisme de \hbox{Brouwer}.

\begin{defi}\index{p\'etale}\label{def.pehb}
Soit $\Delta$ une droite de \hbox{Brouwer} pour $H$. On appellera \res{p\'etale
attractif de $\Delta$} l'unique disque topologique ferm\'e $P^{+}(\Delta)$
de la sph\`ere $\R^2 \cup \{\infty\}$ dont la fronti\`ere est
$\Delta\cup\{\infty\}$ et qui contient $H(\Delta)$. De m\^eme, le  \res{p\'etale
r\'epulsif de $\Delta$} sera le p\'etale attractif de $\Delta$ pour
l'hom\'eo\-morphisme $H^{-1}$, not\'e $P^-(\Delta)$.
\end{defi}
Nous avons d\'ej\`a vu que $P^+(\Delta)$ est un attracteur strict
(remarque~\ref{rem.attr}). On confondra parfois le p\'etale et sa trace dans le
plan. 

\subsubsection*{$\bullet$ Croissants}
\begin{defi}\index{bande topologique}
Si $\Delta_0$ et $\Delta_1$ sont deux droites topologiques disjointes,
rappelons que $D(\Delta_0,\Delta_1)$  d\'esigne l'unique disque
topologique ouvert du plan
de fronti\`ere $\Delta_0 \cup \Delta_1$. 
On dira que l'adh\'erence, dans le
plan $\R^2$, de l'ensemble  $D(\Delta_0,\Delta_1)$ est une
\res{bande topologique}.
 Les deux droites $\Delta_0$ et $\Delta_1$ seront appel\'ees
\res{bords} de la bande.
\end{defi}

\begin{defi}\index{croissant}
\index{$D(\Delta_0,\Delta_1)$}
Si $(\Delta_0,\Delta_1)$ est un couple attractif de droites de \hbox{Brouwer}
disjointes, on appellera \res{croissant (attractif)  de $(\Delta_0,\Delta_1)$}
la bande topologique $\adhe(D(\Delta_0,\Delta_1))$.
 On d\'efinit de m\^eme le croissant (r\'epulsif) associ\'e \`a un couple
r\'epulsif.
\end{defi}
Nous avons d\'ej\`a vu que les croissants attractifs ainsi d\'efinis sont des
attracteurs stricts (affirmation~\ref{aff.topo}, figure~\ref{fig.couple-attractif}).

\subsubsection*{$\bullet$ Bouts d'une bande topologique}
\index{croissant!bouts}
Soit $A=\adhe(D(\Delta_0,\Delta_1))$ une bande topologique.
On notera $\hat A$ le compactifi\'e en bouts de $A$~; $\hat A$ est
hom\'eomorphe au disque ferm\'e. On le muni d'une orientation compatible avec
celle du plan. 
 On notera  $\t N$ et $\t S$ les deux bouts de
 $A$, de mani\`ere \`a ce qu'en parcourant la courbe bordant $\hat A$ dans
le sens positif, on rencontre $\Delta_0$, $\t S$, $\Delta_1$ et $\t N$
dans cet ordre (autrement dit, si l'on se place dans une carte o\`u
$\Delta_0$ et $\Delta_1$ sont des
 droites euclidiennes verticales avec $\Delta_0$ \`a gauche de
$\Delta_1$, $\t N$ est le bout ``en haut'' et $\t S$ est le bout ``en
bas''~: voir la figure~\ref{fig47}). Remarquons que ce choix de notation d\'epend de l'ordre des deux
droites de \hbox{Brouwer} dans le couple $(\Delta_0,\Delta_1)$
(si on permute ces deux droites, les bouts $\t N$ et $\t S$
sont \'egalement permut\'es).

Supposons de plus que $A$ est un croissant attractif.
Puisque $H(A) \subset A$, et que $H$ est propre, la restriction de $H$
\`a $A$ s'\'etend \`a $\hat A$. De plus, la configuration topologique des
deux droites et de leurs images (affirmation~\ref{aff.topo},
figure~\ref{fig.couple-attractif})
 montre que cette extension fixe les deux bouts $\t N$ et $\t S$.

\subsubsection*{$\bullet$ Croissants \`a dynamique  $\t N$-$\t S$}
\begin{defi}\index{croissant!\`a dynamique $\t N$-$\t S$}\label{def.cdns}
Un croissant attractif $A=\adhe(D(\Delta_0, \Delta_1))$  est dit \res{\`a dynamique $\t N$-$\t S$}
   si  pour tout voisinage $O_{\tilde N}$ de $\t N$ dans $\hat A$, il existe un p\'etale
   attractif $P$ pour $H$ tel que 
\begin{itemize}
\item $P \subset A$ ;
\item $A \setminus P \subset O_{\tilde N}$.
\end{itemize} 
\end{defi}
On d\'efinit de la m\^eme mani\`ere  les croissants attractifs \`a dynamique
 $\t S$-$\t N$, et les croissants r\'epulsifs \`a dynamique  $\t
N$-$\t S$ ou  $\t S$-$\t N$.

\subsubsection*{$\bullet$ P\'etales et croissants simpliciaux}

\begin{defi}\index{simplicial}
 Soit $F$ une d\'ecom\-po\-si\-tion
en briques pour $H$. On dira alors qu'une droite topologique est
\res{simpliciale} si elle est incluse dans $F$ (c'est-\`a-dire r\'eunion
d'ar\^etes de $F$). On dira qu'un p\'etale ou un croissant est simplicial
si sa fronti\`ere (dans $\S^2\setminus \fixe(h)$) est incluse dans $F$. On dira enfin qu'un croissant \`a
dynamique $\t N$-$\t S$ ou  $\t S$-$\t N$ est simplicial si c'est un
croissant simplicial, et si les p\'etales de la
d\'efinition~\ref{def.cdns} (ou \ref{def.cans}) peuvent \^etre choisis simpliciaux.
Ces d\'efinitions se transposent sans difficult\'e au cadre d'un
hom\'eo\-mor\-phisme de la sph\`ere ayant deux points fixes (hypoth\`ese (H2)).
\end{defi}

\begin{defi}\index{sous-croissant}
Soit $D$ une bande topologique~; un croissant  $A'$ est un
\res{sous-croissant} de $D$ si $A'$ est inclus dans $D$,  et si
 l'int\'erieur de $A'$ s\'epare les deux bords de $D$ (de mani\`ere \'equivalente,
si l'adh\'erence de $A'$  dans $\hat A$ contient les deux bouts de $\hat A$).
\end{defi}
Remarquons que si $A'$ est un sous-croissant d'une bande topologique $D$, alors
les bouts de $A'$ s'identifient aux bouts de $D$.

\begin{defi}\index{croissant!minimal}
On dira qu'un croissant attractif simplicial $A$ est \res{minimal} si
il ne poss\`ede pas de sous-croissant attractif simplicial diff\'erent de $A$.
\end{defi}

\subsection{Une version du th\'eo\-r\`eme principal pour les hom\'eo\-mor\-phismes
de \hbox{Brouwer}}\label{sub.tphb}
Dans cette partie, nous commençons par \'enoncer une version du th\'eor\`eme principal
pour les hom\'eo\-mor\-phismes de \hbox{Brouwer} (th\'eor\`eme~\ref{the.prinbrou}). Cet \'enonc\'e
 est ensuite d\'ecoup\'e en trois morceaux~: la proposition~\ref{pro.cam}
affirme qu'il existe des croissants simpliciaux minimaux, et qu'ils contiennent
toujours suffisamment de p\'etales~; le lemme~\ref{lem.codb1} permettra
de parler de famille maximale de croissants simpliciaux minimaux~;
enfin, la proposition~\ref{pro.codb2} calcule l'indice partiel
dans une bande en fonction du type des croissants minimaux
qu'elle contient.

\begin{theo}\label{the.prinbrou}
Soient $H$ un hom\'eo\-mor\-phisme de \hbox{Brouwer}, et $F$ une d\'ecom\-po\-si\-tion en
briques pour $H$. Soit $(\Delta_0,\Delta_1)$ un
couple de droites de \hbox{Brouwer} simpliciales disjointes, tel que
$\inpa(\Delta_0,\Delta_1)=-p<0$.

La bande $\adhe(\Delta_0,\Delta_1)$ contient alors $2p$
sous-croissants simpliciaux minimaux~:
 $p$ croissants attractifs \`a dynamique $\t N$-$\t S$,
et $p$ croissants r\'epulsifs \`a dynamique $\t S$-$\t N$,
d'int\'erieurs disjoints deux \`a
deux, les croissants attractifs et r\'epulsifs \'etant altern\'es
entre $\Delta_0$ et $\Delta_1$.
\end{theo}

Remarquons que dans le cas o\`u l'indice partiel est strictement
 positif,  on peut appliquer le th\'eor\`eme au couple
 $(\Delta_1,\Delta_0)$. On peut tester le th\'eor\`eme sur l'exemple
 ``multi-Reeb'' de la figure~\ref{fig35} (section~\ref{sub.exip}).

\begin{prop}
Soit $H$ un hom\'eo\-morphisme de \hbox{Brouwer}, et $F$ une d\'ecom\-po\-sition en
briques pour $H$.
\label{pro.cam}
  \begin{enumerate}
  \item Tout croissant attractif simplicial $A$ poss\`ede un sous-croissant
  attractif simplicial  minimal~;
  \item tout croissant attractif simplicial mini\-mal est \`a dynamique
  $\t S$-$\t N$ ou $\t N$-$\t S$.
  \end{enumerate}
L'\'enonc\'e est encore vrai en changeant ``attractif'' en ``r\'epulsif''.
 \end{prop}

\begin{figure}[htp]
\par
\centerline{\hbox{\input{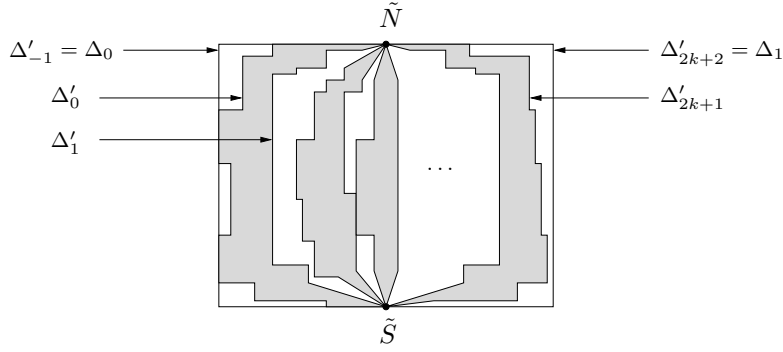}}}
\par
\caption{\label{famille_F}Famille maximale de croissants minimaux}
\end{figure}

\begin{lemm}[figure~\ref{famille_F}]
\label{lem.codb1}
Soit $H$ un hom\'eo\-morphisme de \hbox{Brouwer}, et $F$ une d\'ecom\-po\-sition en
briques pour $H$.
Soit $(\Delta_0,\Delta_1)$ un couple de droites de \hbox{Brouwer}
simpliciales disjointes.

Il existe alors un entier $k\geq -1$ et  une famille (\'eventuellement vide) 
$${\cal F}=\{\adhe(D(\Delta'_0,\Delta'_1)),\dots,
\adhe(D(\Delta'_{2k},\Delta'_{2k+1}))\}$$
 de sous-croissants simpliciaux de la bande topologique $\adhe(D(\Delta_0,\Delta_1))$, telle que
\begin{itemize}
\item les disques topologiques ouverts
$D(\Delta'_{2i},\Delta'_{2i+1})$ sont deux \`a deux disjoints~;
\item ces croissants sont minimaux~;
\item la famille $\cal F$ est
maximale pour l'inclusion~: autrement dit, il n'existe pas de famille ${\cal F}'$ de sous-croissants
 de $D(\Delta_0,\Delta_1)$, d'int\'erieurs deux \`a deux disjoints, simpliciaux, minimaux, qui 
contienne strictement $\cal F$. 
\end{itemize}
\end{lemm}

\begin{demo}[du lemme~\ref{lem.codb1}]
On se place sous les hypoth\`eses du lemme.
Soit $k_0$ la \emph{largeur simpliciale} de
la bande $\adhe(D(\Delta_0,\Delta_1))$, c'est-\`a-dire le plus petit entier tel
qu'il existe une famille de $k_0$ ar\^etes de la d\'ecom\-po\-sition $F$,
dont la r\'eunion est connexe et rencontre $\Delta_0$ et $\Delta_1$.

Si ${\cal F}=\{ C_1, \cdots, C_k \}$ est une famille de croissants
simpliciaux, attractifs ou
r\'epulsifs pour $H$, d'int\'erieurs deux \`a deux disjoints, et 
tels que l'int\'erieur de chaque
croissant s\'epare $\Delta_0$ et $\Delta_1$,
alors il est
clair que $k \leq k_0$.

Ceci entra\^{\i}ne imm\'ediatement l'existence d'une  famille $\cal F$ de
croissants minimaux qui
est maximale pour l'inclusion.
\end{demo}

Si $\cal F$ est une famille de croissants v\'erifiant les conclusions du
  lemme pr\'ec\'edent, on supposera
  toujours que les droites $\Delta'_i$ sont num\'erot\'ees
  dans l'ordre naturel de $\Delta_0$ vers $\Delta_1$\footnote{Autrement dit, que
  l'int\'erieur du croissant $D(\Delta'_{2i},\Delta'_{2i+1})$ s\'epare
  la droite $\Delta_0$, la droite $\Delta'_{2i}$ et les croissants
 $D(\Delta'_{2j},\Delta'_{2j+1})$ pour
  $j<i$, d'une part, de la droite $\Delta_1$, de la droite
  $\Delta'_{2i+1}$, et des croissants $D(\Delta'_{2j},\Delta'_{2j+1})$
  pour $j>i$, d'autre part.}. Dans ce cas, le bout $\t N$ de
  chacun des croissants de $\cal F$ s'identifie naturellement au bout $\t N$
   de la bande $\adhe(D(\Delta_0,\Delta_1))$ (m\^eme remarque pour les bouts
  $\tilde S$). On posera aussi
  $\Delta'_{-1}=\Delta_0$  et $\Delta'_{2k+2}=\Delta_1$.

 On note alors $\#{({\cal F},a,\t N\t S)}$, $\#{({\cal F},r,\t S\t
N)}$, $\#{({\cal F},a,\t S\t N)}$ et
$\#{({\cal F},r,\t N\t S)}$ le nombre de croissants de chaque type
  dans $\cal F$ (attractif $\t N$-$\t S$, \emph{etc.}).

\begin{rema}\label{rem.parite}
Dans une telle famille $\cal F$, deux croissants successifs
$\adhe(D(\Delta'_{2i},\Delta'_{2i+1}))$  et
$\adhe(D(\Delta'_{2i+2},\Delta'_{2i+3}))$ sont toujours l'un attractif et
l'autre r\'epulsif. En effet, si par exemple ils \'etaient tous deux
attractifs, alors la bande qui les s\'epare serait un croissant r\'epulsif
simplicial, donc contiendrait un croissant r\'epulsif simplicial minimal
d'apr\`es la proposition~\ref{pro.cam}. Ceci contredirait la maximalit\'e
de $\cal F$. Le m\^eme raisonnement montre que si le couple $(\Delta_0,\Delta_1)$ est
indiff\'erent, $\cal F$ contient autant de croissants attractifs que de
r\'epulsifs. 
\end{rema}

\index{cha\^{\i}ne!et indice partiel}
\begin{prop}\label{pro.codb2}\index{cha\^{\i}ne!entre deux droites de \hbox{Brouwer}}
Soit $H$ un hom\'eo\-morphisme de \hbox{Brouwer}, et $F$ une d\'ecom\-po\-sition en
briques pour $H$.
Soit $(\Delta_0,\Delta_1)$ un couple de droites de \hbox{Brouwer}
simpliciales disjointes, et $\cal F$ une famille maximale de
sous-croissants simpliciaux minimaux de la bande
$\adhe(D(\Delta_0,\Delta_1))$ (c'est-\`a-dire que $\cal F$ v\'erifie la
conclusion du lemme~\ref{lem.codb1}). Alors~:
\begin{enumerate}
\item pour tout $i=-1,\dots,k$,
l'indice partiel $\inpa(H,\Delta'_{2i+1},\Delta'_{2i+2})$ entre les bords
de deux croissants successifs est nul, et si
ces deux droites sont disjointes, il existe une cha\^ine de pseudo-disques de
l'une \`a l'autre, constitu\'ee de briques de la
d\'ecom\-po\-sition\footnote{D\'efinition~\ref{def.chtr}.}~; 
\item on a la formule suivante~:
$$
\inpa(H,\Delta_0,\Delta_1)= \frac{ \left( \#{({\cal F},a,\t S\t
N)}+\#{({\cal F},r,\t N\t S)} \right) -
\left( \#{({\cal F},a,\t N\t S)}+\#{({\cal F},r,\t S\t N)} \right)}{2}.
$$
\end{enumerate}
\end{prop}
Remarquons en particulier que si l'indice partiel entre $\Delta_0$ et $\Delta_1$ est
tr\`es n\'egatif, on obtient un grand nombre de croissants attractifs \`a dynamique
$\t N$-$\t S$ ou r\'epulsifs \`a dynamique  $\t S$-$\t
N$. Inversement, si la famille $\cal F$ est vide, on a $k=-1$,
et la proposition affirme l'existence d'une cha\^{\i}ne de briques de
$\Delta_0$ \`a $\Delta_1$.

\subsection{Existence d'une droite de \hbox{Brouwer} \`a extr\'emit\'es Nord-Sud}\label{sub.dbns}

\begin{prop}\label{pro.dbns}
Soit $h$ un hom\'eo\-morphisme de  la
sph\`ere $\S^2$, pr\'eservant l'orientation, avec $\fixe(h)=\{N,S\}$
(\textbf{hypoth\`ese (H2)}). On suppose  que $\indi(h,N) \neq 1$.
 Soit $F$ une d\'ecom\-po\-sition en briques pour $h$.

  Alors il existe au moins une droite de \hbox{Brouwer} $\Delta_0$ pour $h$,
  simpliciale (\emph{i.e.} incluse dans $F$), dont l'adh\'erence contient $N$ et $S$.
\end{prop}

\subsection{D\'efinition d'un indice pour les hom\'eo\-morphismes de \hbox{Brouwer}
commutant \`a la translation}\label{sub.inct}
On se place sous l'\textbf{hypoth\`ese (H2)}~: $h$ est un hom\'eo\-morphisme de  la
sph\`ere $\S^2$, pr\'eservant l'orientation, avec $\fixe(h)=\{N,S\}$.
\index{$\A$}
\index{$\t \A$}
Nous allons m\'elanger trois points de
vue, qu'on appelle ``mod\`ele de la sph\`ere, du plan, de
l'anneau''. Le passage d'un point de vue \`a l'autre se fait de la
mani\`ere suivante :

\begin{itemize}
\item on identifie $\S^{2} \setminus \{S\}$ au plan $\R^{2}$ \emph{via} un
\homeo\ pr\'eservant l'orientation qui envoie $N$ sur $(0,0)$ (par
exemple la projection st\'er\'eographique) ;
\item  soit $\tau$  la
translation $(\t \theta,r) \mapsto (\t \theta+1,r)$. On identifie l'anneau $\A=\R^{2}/\tau$
 \`a $\R^{2} \setminus \{(0,0)\}$ \emph{via}
l'\homeo\ $(\theta,r) \mapsto \exp(-r+2i\pi \theta)$ (``$\exp$'' est
l'exponentielle complexe, on a \'egalement identifi\'e $\R^2$ au plan complexe).
\end{itemize}
Il est clair que l'\homeo\ $h$  vit dans chacun des trois mod\`eles.  
Ces mod\`eles
sont orient\'es de mani\`eres compatibles, en choisissant $t \mapsto
\exp(it)$ comme orientation du cercle unit\'e du plan.

\begin{figure}[hbtp]
  \par
\centerline{\hbox{\input{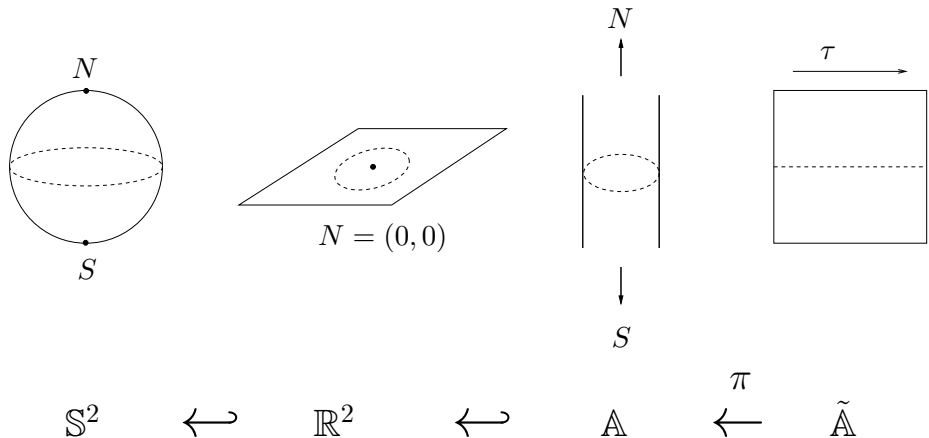}}}
\par
  \caption{\label{fig41}Les diff\'erents mod\`eles}
\end{figure}

On a un quatri\`eme point de vue avec le rev\^etement universel de l'anneau $\A$ :
$$
\pi : 
\left\{
\begin{array}{rcl}
\R^{2} &\longrightarrow & \A \\
(\t \theta,r) & \longmapsto & (\t \theta+\Z,r).
\end{array}
\right.
$$
 Pour ne pas confondre ce dernier point de
vue avec le mod\`ele du plan, on notera $\t \A$ le rev\^etement
universel de $\A$. La figure \ref{fig41} r\'esume les diff\'erents points de vue.

\begin{defi}
\index{relev\'e}
Un \res{relev\'e de $h$} par $\pi$ est un \homeo\ $H$ de $\t \A$ tel que 
 $\pi H=h \pi$. Si $C$ est un sous-ensemble connexe de $\A$, un
\res{relev\'e de $C$} est une composante connexe de $\pi^{-1}(C)$.
\end{defi}


On note maintenant $\com(\tau)$ l'ensemble des \homeo s du plan $\t \A$,
pr\'eser\-vant l'orientation,  qui commutent avec la
translation $\tau :(\t \theta,r) \mapsto (\t \theta+1,r)$ ; c'est aussi l'ensemble des
relev\'es d'\homeo s de l'anneau $\A$ isotopes \`a l'identit\'e (voir
par exemple \cite{lero12}, section 3). Comme
$h$ pr\'eserve l'orientation et fixe $N$ et $S$, il est isotope \`a
l'identit\'e dans le mod\`ele de l'anneau. Autrement dit~: 
\begin{affi}
Tout relev\'e $H$ de $h$ est un \'el\'ement de $\com(\tau)$, sans point fixe.
\end{affi}

Nous allons montrer~:
\begin{affi}
\label{aff.comm}
Soit $H \in \com(\tau)$ un hom\'eo\-morphisme de \hbox{Brouwer}
 et   $\gamma$ une courbe du plan v\'erifiant $\tau(\gamma(0))=\gamma(1)$.
 Le nombre $\indi(H, \gamma)$ est un entier qui ne d\'epend pas de la
 courbe $\gamma$. De plus,
 il est invariant par conjugaison~: si $g \in \com(\tau)$, on a $\indi(gHg^{-1},\gamma)
=\indi(H,\gamma)$.
\end{affi}

\begin{defi}
\index{$\com(\tau)$}
\index{indice!dans $\t \A$}
  Cet entier est appel\'e \res{indice de $H$} (par rapport \`a $\tau$), et not\'e $\indi(H,\tau)$.
\index{$\indi(H,\tau)$}
\end{defi}

\begin{demo}[de l'affirmation~\ref{aff.comm}] 
 Soit $\Gamma$ l'espace des courbes $\gamma$ du plan v\'erifiant
 $\gamma(1)=\tau(\gamma(0))$, muni de la topologie de la convergence
 uniforme. Cet espace est connexe par arcs (et m\^eme convexe).
  Le nombre $\indi(H,\gamma)$ d\'epend contin\^ument de
 $\gamma$~; comme $H$ commute avec la translation~$\tau$, il prend
 des valeurs enti\`eres sur $\Gamma$~; donc ce nombre ne d\'epend
 pas de $\gamma$.
\index{espace d'hom\'eo\-mor\-phismes}
D'autre part, l'espace $\com(\tau)$ est connexe (puisque constitu\'e de
 relev\'es d'hom\'eo\-morphismes isotope \`a l'identit\'e). Le nombre $\indi(gHg^{-1},\tau)$
d\'epend contin\^ument de $g$, donc $\indi(H,\tau)$ est constant sur
l'espace connexe des conjugu\'es \`a $H$ \emph{via} un \'el\'ement de $\com(\tau)$.
\end{demo}

\subsection{Existence d'un relev\'e canonique}\label{sub.rele}
\begin{prop}\label{pro.rele}
Soit $h$ un hom\'eo\-morphisme de  la
sph\`ere $\S^2$, pr\'eservant l'orientation, avec $\fixe(h)=\{N,S\}$
(\textbf{hypoth\`ese (H2)}). On suppose que $\indi(h,N) \neq 1$.

Il existe alors un  unique relev\'e $\t h$ de $h$ par $\pi$ v\'erifiant 
 $\indi(\t h,\tau)=\indi(h,N)-1$. On a de plus~:
\begin{enumerate}\index{arc libre}
\item \label{ite.redr} si $\Delta$ est une droite de
\hbox{Brouwer} pour $h$ dont l'adh\'erence contient $N$
 et $S$, alors tout relev\'e de $\Delta$  est une droite de \hbox{Brouwer} pour $\t h$~;
\item \label{ite.real} si $\gamma$ est un arc libre pour $\t h$, $\pi(\gamma)$ est un
arc libre pour $h$~;
\item \label{ite.rede} si $F$ est une d\'ecom\-po\-sition en briques pour $h$, alors $\t
F=\pi^{-1}(F)$ est une d\'ecom\-po\-sition en briques pour $\t h$.
\end{enumerate}
\end{prop}

\begin{defi}\label{d.reca}\index{relev\'e!canonique}
Le relev\'e $\t h$ de $h$ donn\'e par cette proposition est appel\'e
\res{relev\'e cano\-nique} de $h$.
\end{defi}

\section{Preuve du th\'eo\-r\`eme \`a partir des lemmes}
\label{sec.preu}
Nous allons prouver un \'enonc\'e qui pr\'ecise le th\'eo\-r\`eme~\ref{the.prin}~:
\begin{theo}[``version simpliciale'' du th\'eo\-r\`eme principal]
\label{the.prinbis}
Soit $h$ un hom\'eo\-mor\-phis\-me  de la sph\`ere, pr\'eservant l'orientation,
fixant uniquement les deux points $N$ et $S$ (\textbf{hypoth\`ese
(H2)}). On suppose que $\indi(h,N)=1-p
<1$. Soit $F$ une d\'ecom\-po\-sition en briques pour $h$.

Il existe alors $p$ croissants attractifs simpliciaux, minimaux, \`a dynamique Nord-Sud,
et $p$ croissants r\'epulsifs simpliciaux, minimaux, \`a dynamique Sud-Nord,
d'int\'erieurs disjoints deux \`a
deux, les croissants attractifs et r\'epulsifs \'etant
cycliquement altern\'es autour de $N$ et $S$.
\end{theo}
On d\'eduit facilement le th\'eo\-r\`eme \ref{the.prin} de cet
\'enonc\'e. En effet, on peut toujours trouver une d\'ecom\-po\-sition en
briques pour $h$ (th\'eo\-r\`eme~\ref{the.exde}). On applique alors le
th\'eo\-r\`eme \ref{the.prinbis}. Puis on remplace chaque croissant
attractif par son
image par $h$.
On obtient ainsi une nouvelle famille de croissants  deux \`a
deux d'intersection r\'eduite \`a $\{N,S\}$ (et non plus seulement
d'int\'erieurs disjoints), qui v\'erifie donc la conclusion du
th\'eo\-r\`eme~\ref{the.prin}.
 Remarquons par contre que les nouveaux
croissants ne sont plus simpliciaux.

\begin{demo}[du th\'eo\-r\`eme \ref{the.prinbis}]
Soit $h$ v\'erifiant les hypoth\`eses du th\'eo\-r\`eme, et $F$ une d\'ecom\-po\-sition
en briques pour $h$.

On applique d'abord la proposition~\ref{pro.dbns}
qui nous donne une droite de \hbox{Brouwer} $\Delta_0$, simpliciale, dont l'adh\'erence
contient les points $N$ et $S$.

On utilise maintenant la proposition~\ref{pro.rele}~: on note $H=\t h$
l'hom\'eo\-morphisme de \hbox{Brouwer}, relev\'e canonique de $h$. L'automorphisme du
rev\^etement universel $\pi$ est toujours not\'e $\tau$. On note aussi $\t F$ la
d\'ecom\-po\-sition en briques pour $H$ obtenue en relevant $F$
(point~\ref{ite.rede} de la proposition).
Soient $\t \Delta_0$ un relev\'e de $\Delta_0$, et $\t \Delta_1=\tau(\t
\Delta_0)$. Ces deux ensembles sont des droites de \hbox{Brouwer} pour $H$
(point~\ref{ite.redr} de la proposition).
Elles sont disjointes. 
On peut alors relier l'indice partiel entre ces deux droites \`a
l'indice de $H$ en tant que  commutant \`a la translation~: 
\begin{affi}\label{aff.padi}
$$\inpa(H,\t \Delta_0,\t \Delta_1)=\indi(H,\tau).$$
\end{affi}
En effet, le th\'eo\-r\`eme de Schoenflies permet de supposer que
$\Delta_0$ est une droite verticale de l'anneau $\A$, donc que $\t
\Delta_0$ et $\t \Delta_1$ sont deux droites euclidiennes verticales
dans le plan $\t \A$. Dans cette situation, d'apr\`es les d\'efinitions,
 chacun des deux indices est
\'egal \`a l'indice de $H$ le long de n'importe quelle courbe allant d'un point
$\t x$ de $\t \Delta_0$ \`a $\tau(\t x)$, ce qui prouve l'affirmation.

On applique maintenant le lemme~\ref{lem.codb1}
au couple $(\t \Delta_0,\t\Delta_1)$ de droites de \hbox{Brouwer}
simpliciales pour $H$ (pour la d\'ecom\-po\-sition $\tilde F$).
On obtient ainsi une famille maximale $\t {\cal F}$ de croissants
simpliciaux minimaux. D'apr\`es la proposition~\ref{pro.cam}, ces
croissants sont tous \`a dynamique $\t N$-$\t S$ ou $\t S$-$\t N$.
Et on a  la formule  de la proposition~\ref{pro.codb2} reliant
l'indice partiel
$\inpa(H, \t \Delta_0,\t \Delta_1)$ 
au nombre de croissants de chaque type. Remarquons que d'apr\`es
l'affirmation~\ref{aff.padi} ci-dessus et les propri\'et\'es du relev\'e
canonique, on a
$$\inpa(H, \t \Delta_0,\t \Delta_1)= \indi(H,\tau) = \indi(h,N)-1=-p<0
\ ;$$
 la formule  de la proposition~\ref{pro.codb2}
montre alors que la famille $\t {\cal F}$ contient
au moins 2 croissants.

D'autre part, la restriction de la projection $\pi$ au disque
topologique ouvert
$D(\t\Delta_0,\t\Delta_1)$ se prolonge contin\^ument aux bouts $\t N$
et $\t S$ de la bande $\adhe(D(\t\Delta_0,\t\Delta_1))$, en une
application $\hat \pi$ qui envoie $\t N$ sur $N$ et $\t S$ sur $S$~;
cette application est un hom\'eo\-morphisme sur son image.

Ceci permet de voir facilement que les projections des croissants de
 $\t {\cal F}$ sont des croissants  pour $h$, et que ``les types
 dynamiques se correspondent'' (c'est-\`a-dire que l'image par $\pi$ d'un croissant attractif \`a
 dynamique $\t N$-$\t S$ pour $H$ est un croissant attractif \`a dynamique
 Nord-Sud pour $h$, \textit{etc.}).

\paragraph{Le compte y est}
Soit maintenant ${\cal F}$ l'image par $\pi$ de la famille $\t {\cal
F}$. C'est une famille de croissants attractifs ou r\'epulsifs, \`a
dynamique Nord-Sud ou Sud-Nord, d'int\'erieurs deux \`a deux disjoints.
Il reste \`a extraire de la famille $\cal F$  une sous-famille ${\cal F}'$
constitu\'ee de $p$ croissants attractifs Nord-Sud et $p$ croissants r\'epulsifs
Sud-Nord, la difficult\'e \'etant que les croissants des deux types
doivent \^etre cycliquement altern\'es autour de $N$.

Rappelons tout d'abord que $\cal F$ contient exactement le m\^eme
nombre de croissants attractifs que de croissants r\'epulsifs, et qu'ils
sont altern\'es autour de $N$. Ceci provient de la maximalit\'e de $\t
{\cal F}$ (remarque~\ref{rem.parite}.

D'autre part, la formule de la proposition~\ref{pro.codb2} donne alors (avec des
notations \'evidentes)~:
$$
(*) \ \ \ \ \underbrace{
\left( \#{({\cal F},a,NS)}+\#{({\cal F},r,SN)} \right)
}_{\mbox{types recherch\'es}}
-\left( \#{({\cal F},a,SN)}+\#{({\cal F},r,NS)} \right)
 =2p.
$$  
Il y a deux cas :
\begin{enumerate}
\item Si la famille $\cal F$ ne contient que des croissants des types recherch\'es
(attractifs Nord-Sud et r\'epulsifs Sud-Nord), comme il y a autant
d'attractifs que de r\'epulsifs, il y en a $p$ de chaque sorte, et $\cal
F$ v\'erifie la conclusion du th\'eor\`eme.

\item Dans le cas contraire, d'apr\`es la formule,  puisque $2p >0$,
 $\cal F$ contient n\'eanmoins des croissants des
types recherch\'es ; il y a donc dans $\cal F$ deux croissants
 successifs, l'un qui est d'un des deux types recherch\'es, l'autre non. On
 consid\`ere la sous-famille ${\cal F}_1$ obtenue \`a partir de $\cal F$
 en enlevant ces deux croissants. Comme on a enlev\'e un croissant
 non-recherch\'e et un croissant recherch\'e, la formule $(*)$ est encore
 valable en remplaçant $ \cal F$ par la sous-famille  ${\cal F}_1$ ; comme ils \'etaient
adjacents, ${\cal F}_1$ est encore constitu\'ee de croissants
alternativement attractifs et r\'epulsifs. 
\end{enumerate}

On it\`ere ce proc\'ed\'e d'extraction jusqu'\`a obtenir une sous-famille
${\cal F}_l$ qui ne contient plus que des croissants des types
recherch\'es. Cette sous-famille satisfait encore \`a la formule $(*)$, et
v\'erifie donc la conclusion du th\'eor\`eme d'apr\`es l'\'etude du premier cas.

Ceci termine la preuve du th\'eo\-r\`eme \ref{the.prinbis}.
\end{demo}

\section{\'Etude des croissants minimaux (Proposition~\ref{pro.cam})}
\label{sec.cam}
Pour cette section, on suppose que  $H$ est un hom\'eo\-morphisme de
\hbox{Brouwer}, et que $F$ est une d\'ecom\-po\-sition en briques pour $H$.
 On consid\`ere un couple attractif $(\partial^-A,\partial^+A)$ de
droites de \hbox{Brouwer} disjointes, le croissant $A$ associ\'e, et sa
fronti\`ere $\partial A=\partial^-A \cup \partial^+A$.
On note $N$ et
$S$ les deux bouts de $A$.
Tous les ensembles consid\'er\'es seront suppos\'es \^etre inclus dans $\hat
A=A\cup\{N,S\}$.\footnote{
De fait, nous n'utiliserons que la restriction de $H$ \`a $A$, et 
ce qui suit est valable pour tout plongement de $\hat A$ dans
 $\hat A$, v\'erifiant des hypoth\`eses \'evidentes.}

Remarquons notamment que, puisque $A$ est un attracteur strict simplicial, pour
toute brique $B$ incluse dans $A$, l'attracteur $A^+(B)$ et la droite
de \hbox{Brouwer} $\Delta(B)$ sont \'egalement inclus dans $A$.

\subsection{Types topologiques et types dynamiques des droites de \hbox{Brouwer}}
\label{ss.type}

\begin{figure}[h!tp]
  \par
\centerline{\hbox{\input{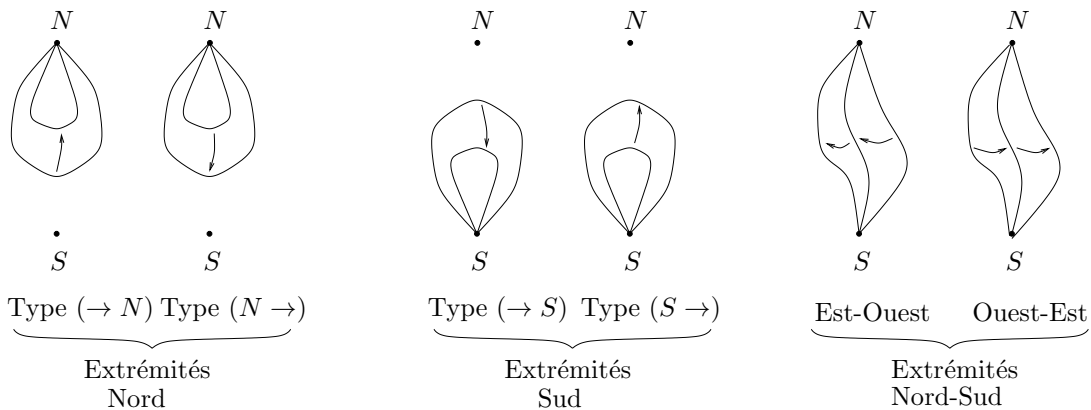}}}
\par
  \caption{\label{fig33}Extr\'emit\'es et types dynamiques des droites de \hbox{Brouwer}}
\end{figure}

\begin{defi}[ (figure \ref{fig33})]
\index{droite de \hbox{Brouwer}!type topologique}
\index{droite topologique!extr\'emit\'es}
Soit  $\Delta$ une droite topologique dans $A$~; on appelle
\res{extr\'emit\'es de $\Delta$} l'intersection de l'ensemble $\{N,S\}$
avec l'adh\'erence de $\Delta$ (dans $\hat A$). 
Il y a donc trois possibilit\'es, et selon les cas,
 la droite topologique $\Delta$ est dite \res{\`a
extr\'emit\'es Nord}, \res{\`a extr\'emit\'es Sud} ou \res{\`a extr\'emit\'es
 Nord-Sud}  (ou encore \res{de type
topologique Nord, Sud} ou \res{Nord-Sud}).
\end{defi}

\begin{defi}[ (figure \ref{fig33})]
\index{droite de \hbox{Brouwer}!type dynamique}
Une droite de \hbox{Brouwer} $\Delta$, incluse dans $A$, qui est  une droite
topologique \`a extr\'emit\'es Nord est dite \res{de type dynamique $\fn$}
 si $\Delta \cup \{N\}$ s\'epare  $H(\Delta)$ et le point
fixe $S$ (ou, de mani\`ere \'equivalente, si $H(\Delta)\cup \{N\}$ ne s\'epare pas
$\Delta$ et $S$). 

Dans le cas contraire, elle est dite \res{de type dynamique $\nf$}. On
d\'efinit de m\^eme les droites de \hbox{Brouwer} \res{de type dynamique
$\fs$} et \res{$\sf$}.
\end{defi}
Remarquons que si $\Delta$ est de type dynamique $\fn$, alors le
p\'etale attractif $P^+(\Delta)$ est inclus dans $\hat A$~; de m\^eme,
si  $\Delta$ est de type dynamique $\nf$, alors le
p\'etale r\'epulsif $P^-(\Delta)$ est inclus dans $\hat A$.

\begin{defi}[ (figure \ref{fig33})]
Une droite de \hbox{Brouwer} $\Delta$ qui est  une droite
topologique \`a extr\'emit\'es Nord-Sud est dite \res{de type
dynamique Est-Ouest} si $H(\Delta)$ s\'epare $\Delta$ de $\partial^-A$
dans $A$~; dans le cas contraire, elle est dite \res{de type
dynamique Ouest-Est}.
\end{defi}
Par exemple, les droites $\partial^-A$ et $\partial^+A$ sont
respectivement de types dynamiques Ouest-Est et Est-Ouest.

\begin{figure}[h!tp]
\par
\centerline{\hbox{\input{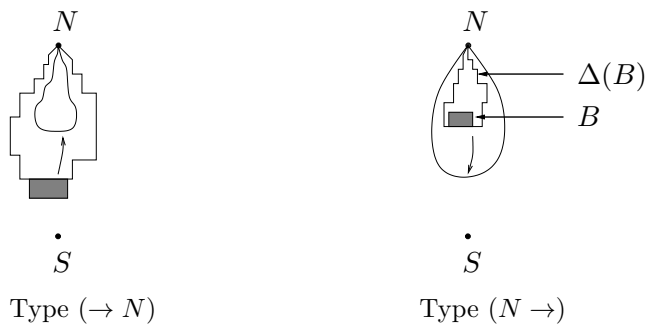}}}
\par
\caption{\label{fig39}Briques de type Nord}
\end{figure}

\begin{defi}[ (figure \ref{fig39})]
\index{brique!type topologique}
\index{brique!type dynamique}
Le \res{type topologique d'une brique} $B$ de la d\'ecom\-po\-sition $F$,
incluse dans $A$,  est le type topologique  de
sa droite de \hbox{Brouwer}   associ\'ee $\Delta(B)$.
  De m\^eme, le \res{type dynamique} de $B$  est le type dynamique de
$\Delta(B)$.
\end{defi}
 
 Si $B$ est une brique de type dynamique
    $\fn$, $A^{+}(B)$ est un attracteur dont l'adh\'erence (dans $\hat
    A$) ne contient pas
    $S$. Puisque les points errent (lemme~\ref{lem.libr} et
    corollaire~\ref{cor.erre}), les it\'er\'es positifs des 
    points de $A^+(B)$ tendent alors 
    vers $N$ (et aussi ceux des points de $B$ puisque $H(B) \subset
    A^+(B)$). De m\^eme, les it\'er\'es positifs d'une brique de type
    dynamique 
    $\fs$ tendent vers $S$. On en d\'eduit que les types dynamiques $\fn$ et $\fs$ 
sont  ``incompatibles''~:
\begin{affi}
\index{brique!incompatibilit\'e}
\label{aff.brfn}
   Si $B$ est une brique de type dynamique  $\fn$, et $B'$ une brique
 de type  dynamique $\fs$, alors $B$ et $B'$  ne peuvent pas \^etre
 adjacentes. Mieux, les ensembles  $B\cup A^+(B)$ et 
$B'\cup  A^+(B')$ sont disjoints.
\end{affi}

De m\^eme, si $B$ est une brique de type dynamique $\nf$, le p\'etale
r\'epulsif de $\Delta(B)$, $P^-(\Delta(B))$, est un attracteur pour
$H^{-1}$, inclus dans $\hat A$, dont l'adh\'erence  ne contient pas $S$. On en d\'eduit~:
\begin{affi}
\label{aff.brfn2}
   Si $B$ est une brique de type dynamique  $\nf$, et $B'$ une brique
 de type  dynamique $\sf$, alors $B$ et $B'$  ne peuvent pas \^etre
 adjacentes.
\end{affi}

La preuve de la derni\`ere affirmation est laiss\'ee au lecteur~:
\begin{affi}\label{aff.OE}
Soit $\Delta$ une droite de \hbox{Brouwer} \`a extr\'emit\'es Nord-Sud,
simpliciale, incluse dans $A$. Si $\Delta$ rencontre $\partial^+A$,
alors $\Delta$ est de type dynamique Est-Ouest.
\end{affi}

\subsection{Courbes simpliciales}
\label{sub.svs}
\index{${\cal V}_1$}
\index{sous-vari\'et\'e simpliciale}
Dans cette section, nous soulignons l'un des principaux int\'er\^ets des
d\'ecom\-po\-si\-tions en briques~: la recherche de droites de Brouwer
simpliciales est grandement facilit\'ee par des propri\'et\'es de passage \`a la limite.

Soit ${\cal V}_1$ l'ensemble des sous-vari\'et\'es du plan, incluses dans $A$,
 de dimension $1$, \'eventuellement \`a bord, qui sont des r\'eunions d'ar\^etes de la
d\'e\-com\-po\-si\-tion $F$.\footnote{Autrement dit, toute r\'eunion
d'une famille d'ar\^etes qui ne contient pas trois ar\^etes ayant un
sommet commun.}
 Les \'el\'ements de ${\cal V}_1$ sont appel\'es \emph{sous-vari\'et\'es simpliciales}. Cet
ensemble ${\cal V}_1$ est muni d'une topologie de Hausdorff (induite par
le plongement $C \mapsto \adhe(C)$ de $\cal V$ dans l'espace des
compacts de $\hat A$). On a un crit\`ere tr\`es simple de
convergence : une suite $(V_k) \subset {\cal V}_1$ converge vers
$V_\infty$ si et seulement si pour toute brique $B$ de la
d\'ecom\-po\-si\-tion, pour tout $k$ assez grand, $B \cap V_k = B
\cap V_\infty$ (autrement dit, une suite converge si et seulement si
sa trace sur tout compact de $A$ est stationnaire). Insistons sur le fait que
${\cal V}_1$ contient l'ensemble vide.

\begin{affi}
\label{aff.comp}\hspace{1cm}
\begin{enumerate}
\item L'espace topologique ${\cal V}_1$ est compact ;
\item l'ensemble des sous-vari\'et\'es simpliciales libres pour $H$ est ferm\'e
    dans ${\cal V}_1$.
\end{enumerate}
\end{affi}
Les preuves sont laiss\'ees au lecteur.

Voici un exemple d'utilisation de l'espace ${\cal V}_1$.
\begin{lemm}\label{lem.svs}
Soit $(\delta_k)$ une suite de  sous-vari\'et\'es  simpliciales connexes, telle
qu'il existe deux suites $(x_k)$ et $(y_k)$ de points de $(\delta_k)$
qui tendent respectivement vers $N$ et $S$. Alors toute valeur
d'adh\'erence $L$ de $(\delta_k)$ dans ${\cal V}_1$ contient une droite
topologique \`a extr\'emit\'es Nord-Sud.
\end{lemm}

\begin{demo}
Quitte \`a raccourcir $\delta_k$, on peut supposer que $\delta_k$ est un
arc dont $x_k$ et $y_k$ sont les extr\'emit\'es.
Quitte \`a extraire, on peut de plus supposer que la suite $(\delta_k)$ converge
vers $L$.
Comme  $(x_k)$ tend vers $N$ et $(y_k)$  tend vers $S$,
 la  sous-vari\'et\'e simpliciale $L$  est une
sous-vari\'et\'e sans bord. Si $L$ contenait un cercle topologique, ce
cercle serait aussi inclus dans $\delta_k$ pour $k$ assez grand, ce n'est
 pas le cas~: les composantes connexes
de $L$ sont donc des droites topologiques. Nous
allons montrer que l'une d'entre elles est \`a extr\'emit\'es Nord-Sud.

Il est facile de trouver un ensemble $R$, r\'eunion de briques
incluses dans $A$, qui est un disque topologique ferm\'e, et qui
s\'epare $N$ et $S$ dans $\hat A$ (on prend une courbe $\gamma$ s\'eparant $N$ et
$S$, on prend la r\'eunion des briques qui rencontrent $\gamma$, et on
``bouche les trous'').  Les composantes connexes de $L \cap R$ sont
des arcs inclus dans $F$, ils sont en nombre fini ; soit ${\cal L}$
l'ensemble de ceux qui ``traversent''
$R$ (figure~\ref{fig20}). Pour $k$ assez grand, on a
$\delta_k \cap R=L\cap R$, et on voit que $\cal L$ contient un nombre
impair d'\'el\'ements.
\begin{figure}[hbtp]
  \par
\centerline{\hbox{\input{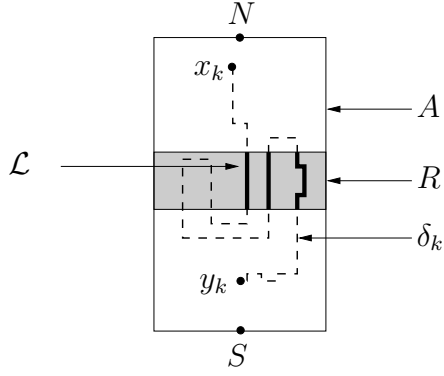}}}
\par
  \caption{\label{fig20}$\cal L$ contient un nombre impair d'\'el\'ements}
\end{figure}
On raisonne par l'absurde : si toutes les composantes connexes de $L$
 \'etaient des droites topologiques de type Nord ou de type Sud,
 $\cal L$ contiendrait un
nombre pair d'\'el\'ements. Par cons\'equent $\cal L$ contient au moins
 une droite topologique \`a extr\'emit\'es Nord-Sud.
\end{demo}

\subsection{\'Etude des croissants attractifs minimaux}

\begin{rema}[ (de lecture)]
La preuve de la proposition~\ref{pro.cam} est la plus ``technique'' de
ce texte.
Cependant, on peut se convaincre facilement que 
les preuves des th\'eo\-r\`emes de dynamique locale~\ref{the.fleu} et \ref{the.brow} 
(version topologique du ``th\'eor\`eme de la fleur'' de Leau-Fatou, et
indice des it\'er\'es) n'utilisent qu'une propri\'et\'e faible des croissants \`a
dynamique Nord-Sud, propri\'et\'e qui s'obtient assez rapidement (\'etape~\ref{ite.escr}
du deuxi\`eme point de la preuve ci-dessous, affirmation~\ref{aff.faib}).
Plus pr\'ecis\'ement, on peut d\'efinir :
\begin{defi}
Un croissant attractif $C$ est dit \res{faiblement \`a dynamique Nord-Sud}
   si  il existe un p\'etale attractif $P$ bas\'e en $S$, inclus dans
   $C$,  qui rencontre les deux bords de $C$.
\end{defi}
Dans tous les \'enonc\'es conduisant aux th\'eor\`emes~\ref{the.fleu} et
\ref{the.brow} (notamment dans le th\'eor\`eme~\ref{the.prinbis} et
l'affirmation~\ref{aff.ipcr}), on peut alors remplacer ``croissant attractif
\`a dynamique Nord-Sud'' par ``croissant attractif faiblement \`a
dynamique Nord-Sud''.

Par contre,  la version topologique du th\'eo\-r\`eme de
la vari\'et\'e stable (th\'eo\-r\`eme~\ref{the.stlo}) n\'ecessite la notion
forte de croissant \`a dynamique Nord-Sud, et donc la lecture de
toute la preuve de la proposition~\ref{pro.cam}.
\end{rema}

\begin{demo}[de la proposition \ref{pro.cam}]
\paragraph{Preuve du premier point}
L'existence d'un
sous-croissant attractif simplicial minimal dans $A$ va
essentiellement d\'ecouler du lemme de Zorn~: l'id\'ee est de consid\'erer
une droite de \hbox{Brouwer} $\Delta^-$, de type dynamique Ouest-Est, ``la
plus \`a droite possible'' dans $A$, puis une droite de \hbox{Brouwer}
$\Delta^+$, de type dynamique Est-Ouest, \`a droite de $\Delta^-$, ``la
plus \`a gauche possible'', pour obtenir un croissant minimal
$\adhe(\Delta^-,\Delta^+)$. D\'etaillons cette id\'ee.

Soit $\cal D$ la famille des droites de \hbox{Brouwer} incluses dans $A$,
 simpliciales, \`a extr\'emit\'es Nord-Sud, et de type dynamique Ouest-Est. Cette
 famille contient $\partial^-A$. Si $\Delta$ est un \'el\'ement de $\cal
 D$, alors $\Delta$ est disjointe de $\partial^+A$ d'apr\`es
 l'affirmation~\ref{aff.OE}. 
\`{A} chaque \'el\'ement $\Delta$ de $\cal D$ correspond donc un croissant
 attractif $\adhe(\Delta,\partial^+A)$ (affirmation~\ref{aff.topo}).
 On munit $\cal D$ de l'ordre
 induit, \emph{via} cette correspondance, par l'inclusion sur les
 croissants.

Montrons rapidement que $\cal D$ est inductive. Soit $(\Delta_n)_{n
\geq 0}$ une
suite d'\'el\'ements de $\cal D$ totalement ordonn\'ee, c'est-\`a-dire telle
que pour tout $n$, $\adhe(\Delta_{n+1},\partial^+A) \subset
\adhe(\Delta_n,\partial^+A)$. D'apr\`es le lemme~\ref{lem.svs}, il existe
une droite topologique $\Delta_\infty$, \`a extr\'emit\'es Nord-Sud, incluse dans une
valeur d'adh\'erence de la suite $(\Delta_n)$. Cette droite est libre
(affirmation~\ref{aff.comp}), c'est donc une droite de \hbox{Brouwer}~; on
montre facilement qu'elle est de type dynamique Ouest-Est~:
$\Delta_\infty$ est donc un \'el\'ement de $\cal D$.\footnote{Il est tout
\`a fait possible que la suite $(\Delta_n)$ ne converge pas vers
$\Delta_\infty$ (et que l'inclusion $\adhe(\Delta_\infty,\partial^+A) \subset
\cap_n \adhe(\Delta_n,\partial^+A)$ soit une inclusion stricte).}
De plus, il est clair que  $\adhe(\Delta_\infty,\partial^+A) \subset
\cap_n \adhe(\Delta_n,\partial^+A)$, ce qui prouve que $\cal D$ est
inductive.

On choisit alors un \'el\'ement $\Delta^-$ de $\cal D$ qui soit minimal
pour l'ordre.

De mani\`ere sym\'etrique, on consid\`ere la famille ${\cal D}'$ des droites
 de \hbox{Brouwer} incluses dans le croissant attractif
 $\adhe(\Delta^-,\partial^+A)$,  simpliciales, \`a extr\'emit\'es
 Nord-Sud, et de type dynamique Est-Ouest. On montre que cette famille est
 inductive. On trouve ainsi une seconde droite de \hbox{Brouwer} $\Delta^+$.
Le croissant attractif $\adhe(\Delta^-,\Delta^+)$ est alors un
 sous-croissant simplicial de $A$, attractif, minimal pour l'inclusion.

\paragraph{Preuve du deuxi\`eme point}
On suppose maintenant que le  croissant attractif simplicial $A$ est minimal.
Voici comment la preuve va se d\'erouler. Nous allons montrer que pour
toute brique $B$, incluse dans $A$ et rencontrant $\partial^- A$, la
droite de \hbox{Brouwer} $\Delta(B)$ est \`a extr\'emit\'es Nord ou \`a extr\'emit\'es
Sud. L'incompatibilit\'e entre ces deux types de briques oblige alors
toutes les briques (dans $A$) adjacentes \`a $\partial^-A$ \`a \^etre du
m\^eme type~; et ce type d\'eterminera le type dynamique du croissant
(Nord-Sud ou Sud-Nord). Si par exemple elles sont toutes \`a extr\'emit\'es
Nord, on \'etudie le p\'etale attractif $P^+(\Delta(B))$ pour $B$ adjacente \`a
$\partial^-A$ et tr\`es proche du point $S$. On montre d'abord que ce
p\'etale rencontre les deux bords de $A$ (ce qui prouve que $A$ est
faiblement \`a dynamique Sud-Nord). Finalement, on prouve que ce p\'etale
est arbitrairement grand dans $A$, ce qui conclut.

\begin{enumerate}
\item \label{ite.NS} 
\begin{affi}
\label{aff.brns}
Les seules droites de \hbox{Brouwer} simpliciales \`a extr\'emit\'es Nord-Sud incluses dans $A$
sont $\partial^- A$ et $\partial^+ A$.
\end{affi}
\begin{demo}
Soit $\Delta'$ une telle droite. Si $\Delta'$ est de type
dynamique Est-Ouest, alors elle ne peut pas rencontrer $\partial^-A$ (affirmation~\ref{aff.OE})~;
et la bande $\adhe(\partial^-A,\Delta')$ est un
 croissant attractif inclus dans $A$  (affirmation~\ref{aff.topo}).
 Par minimalit\'e, on en d\'eduit que
$\Delta'=\partial^+A$. Pour l'autre type dynamique,
on conclut de mani\`ere analogue.
\end{demo}

\item Soit $B$ une brique dans $A$ rencontrant la fronti\`ere $\partial A$.
\begin{affi}
\label{aff.long}
La brique $B$ est de type dynamique $\fn$ ou $\fs$.
\end{affi}

    \begin{demo}
     Comme $A$ est attractif, $\Delta(B)$ ne contient pas $B \cap
    \partial A$, et ne peut pas \^etre \'egale \`a  $\partial^- A$ ou
    \`a  $\partial^+ A$ ;
  elle n'est donc pas \`a extr\'emit\'es Nord-Sud d'apr\`es
    l'affirmation \ref{aff.brns}.

 Supposons par exemple que $B$ est de type topologique~$N$.
Comme $H(\Delta(B)) \subset A \setminus (\partial A \cup B)$, et que $B$
rencontre $\partial A$, la droite $H(\Delta(B))$
ne s\'epare pas $B$ et $S$. Or  $H(\Delta(B))$ ne s\'epare pas non plus $B$ et
$\Delta(B)$. Donc $H(\Delta(B))$ ne s\'epare pas $S$ et $\Delta(B)$, elle est
donc de type dynamique~$\fn$.

Pour les m\^emes raisons,  le type dynamique  $\sf$ est exclu.
    \end{demo}

\item  Puisque les types dynamiques $\fn$ et $\fs$ sont incompatibles
(affirmation \ref{aff.brfn}), toutes
les briques adjacentes \`a $\partial^- A$  sont du m\^eme type
dynamique.
Dans la suite, on se place dans le cas o\`u ce type dynamique commun est
$\fn$ (le cas $\fs$ s'en d\'eduit en \'echangeant les r\^oles de $N$ et $S$).
Soit alors  $(B_i)_{i \geq 0}$ n'importe quelle suite de briques adjacentes au
bord $\partial^{-}A$, telle que $\lim_{i
\rightarrow +\infty} B_i=S$. On note $P^+(B_i)$ le p\'etale attractif
associ\'e \`a la droite de \hbox{Brouwer} $\Delta(B_i)$.

\textbf{Nous allons montrer que pour tout $i$ assez grand, le p\'etale $P^+(B_i)$
est arbitrairement grand dans $A$, ce 
qui prouvera que le croissant attractif $A$ est \`a dynamique Sud-Nord}.  

\item \label{ite.escr} Soit $\Delta_\infty$ une valeur d'adh\'erence 
 de la suite  $(\Delta(B_i))_{i \geq 0}$ dans l'espace topologique
${\cal V}_1$ (section~\ref{sub.svs}).
\begin{affi}
 $\Delta_\infty$ contient $\partial ^-A \cup \partial ^+A$.
\end{affi}
\begin{demo}
Pour tout $i$, on note $\delta^-_i$ et
$\delta^+_i$ les deux composantes connexes de $\Delta(B_i) \setminus
B_i$ (le fait qu'il y en ait exactement deux provient de
l'affirmation~\ref{aff.allu}, figure~\ref{fig38}).
Par compacit\'e de ${\cal V}_1$, les suites de sous-vari\'et\'es  simpliciales
$(\delta^-_i)$ et   $(\delta^+_i)$ ont des valeurs d'adh\'erence. De
plus, elles v\'erifient chacune les hypoth\`eses du lemme~\ref{lem.svs}~: notons
alors $\delta_\infty^-$ et $\delta_\infty^+$ deux droites topologiques \`a
extr\'emit\'es Nord-Sud incluses dans leurs valeurs d'adh\'erence
 respectives. Puisque l'ensemble
des sous-vari\'et\'es  simpliciales libres est compact, $\delta_\infty^-$ et
$\delta_\infty^+$ sont libres~: ce sont deux droites de \hbox{Brouwer}. Puisque
pour tout $i$,  $\delta^-_i \cap \delta^+_i=\emptyset$, on a encore
$\delta_\infty^- \cap \delta_\infty^+=\emptyset$ (le fait d'\^etre
disjoint ``passe \`a la limite'' dans ${\cal V}_1$). Comme les seules
droites de \hbox{Brouwer} simpliciales  \`a extr\'emit\'es Nord-Sud sont  $\partial
^-A$ et $\partial ^+A$ (affirmation~\ref{aff.brns}), on en d\'eduit que 
$\delta_\infty^- \cup \delta_\infty^+=\partial ^-A \cup \partial ^+A$,
ce que l'on voulait.
\end{demo}

On en d\'eduit imm\'ediatement~:
\begin{affi}\label{aff.faib}
Pour $i$ assez grand, $P^+(B_i)$ est un p\'etale attractif bas\'e en $N$
qui rencontre les deux bords de $A$. Notamment, $A$ est faiblement \`a
dynamique Sud-Nord.
\end{affi}

\item \label{ite.ozero} 
Soient $x^{-}$ et $x^{+}$ deux points respectivement sur
$\partial^{-}(A)$ et $\partial^{+}(A)$. D'apr\`es le point
\ref{ite.escr}, quitte \`a r\'e-indexer les $(B_i)$, on peut supposer
que pour tout $i \geq 0$, $\Delta(B_i)$ contient $x^{-}$ et
$x^{+}$. On note (\cf figure \ref{fig3}, dessin de gauche)~: 
\begin{figure}[hbtp]
  \par \centerline{\hbox{\input{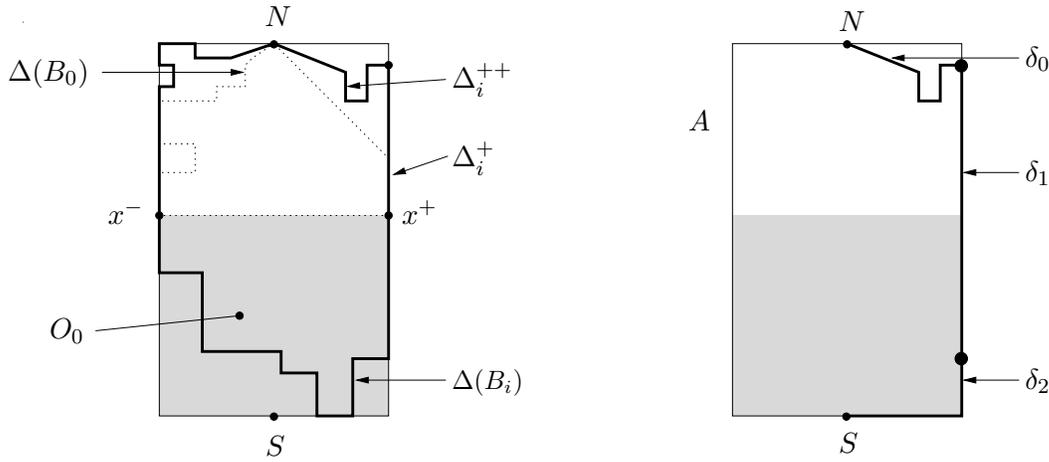}}} \par
\caption{\label{fig3}Anatomie de $\Delta(B_i)$}
\end{figure}

\begin{itemize}
  \item $\Delta^{+}_i$ la composante connexe de $\Delta(B_i) \cap
    \partial^{+}(A)$ contenant $x^{+}$ ; 
\item $\Delta^{++}_i$ la
    composante connexe de $\Delta(B_i) \setminus \Delta^{+}_i$ qui ne
    contient pas $x^{-}$ ; 
\item $O_i$ la composante connexe de $A \setminus
    \Delta(B_i)$ contenant $S$ ; 
\end{itemize} 
et on d\'efinit
    sym\'etriquement $\Delta^{-}_i$ et $\Delta^{--}_i$. Remarquons que
    d'apr\`es le point \ref{ite.escr}, $\Delta^{+}_i$ tend vers
    $\partial^{+}A$ quand $i$ tend vers $+\infty$.

\textbf{La fin de la preuve consiste \`a montrer que pour tout $i$ assez grand, $\Delta ^{++}_i$
est vide, et que $\Delta_\infty=\partial^-A \cup \partial ^+A$}.

\item \label{ite.escr2} La suite $(A^{+}(B_i))_{i\geq 0}$ est
essentiellement croissante ; plus pr\'ecis\'ement :
\begin{affi}~
\label{aff.escr}
\begin{itemize}
\item  Pour tout $j \geq 0$ il existe $i_j$ tel que pour tout $i \geq i_j$,
  $A^{+}(B_j) \subset A^{+}(B_i)$.
\item La suite $\Delta(B_i)$
  converge (vers $\Delta_\infty$) dans ${\cal V}_1$.
\end{itemize}
\end{affi}
 \begin{demo}
Le point \ref{ite.escr} dit que  $B_{j} \cap
 \partial^{-}(A)$ est incluse dans $\Delta(B_i)$ pour tout $i$ assez
 grand, ce qui n'est possible que si $B_{j} \subset A^{+}(B_i)$, ce
 qui entra\^\i ne la premi\`ere partie de l'affirmation puisque $A^{+}(B_j)$ est le plus petit
 attracteur strict,
 r\'eunion de briques, contenant $H(B_j)$.

On d\'eduit de cette premi\`ere affirmation que pour toute brique $B$ incluse dans $A$, de deux
choses l'une~: soit $B$ n'est inclus dans $A^+(B_i)$ pour aucun $i
\geq 0$, soit il existe $n_0$ tel que pour tout $i \geq n_0$,
$B$ est inclus dans $A^+(B_i)$. Par cons\'equent, la trace de la
suite $(A^+(B_i))$ sur tout compact de $A$ est constante \`a partir d'un
certain rang, donc la trace de $\Delta(B_i)$ \'egalement. La seconde
partie de l'affirmation en d\'ecoule.
 \end{demo}

\textbf{Quitte \`a extraire, on peut
supposer que la suite $(A^+(B_i))_{i \geq 0}$ est croissante, ce que nous ferons
dor\'enavant.}

\item \label{ite.N} Soit $i>0$ ; puisque $(A^+(B_i))$ contient $A^+(B_0)$,
il est clair que $\Delta^{++}_i \cap O_0=\emptyset$ (\cf
figure \ref{fig3}). Nous allons en d\'eduire :
\begin{affi}
\label{aff.de++}
Pour tout $i$ assez grand, $\Delta^{++}_i$ et $\Delta^{--}_i$ sont
vides.
\end{affi}
\begin{demo}
  Fixons $i$ ; soit $\delta_0=\Delta^{++}_i$, $\delta_1=\Delta^{+}_i$,
  et $\delta_2$ la composante connexe de $\partial^{+}A \setminus
  \delta_1$ dont l'adh\'erence contient $S$ (figure \ref{fig3},
  dessin de droite) ; on pose $\delta=\delta_0 \cup
  \delta_1 \cup \delta_2$.  Montrons alors que $\delta$ est libre. Les
  deux ensembles
  $\delta'_0=\delta_0 \cup \delta_1$ et $\delta'_2=\delta_1 \cup
  \delta_2$ sont libres, en tant que morceaux des droites de \hbox{Brouwer}
  $\Delta(B_i)$ et $\partial^+A$~;
  par cons\'equent, il y a
  juste \`a montrer que $h(\delta_0) \cap \delta_2=\emptyset$ et
  $h(\delta_2) \cap \delta_0=\emptyset$. La premi\`ere \'egalit\'e est
  toujours vraie car $A$ est un attracteur strict ; la deuxi\`eme est
  v\'erifi\'ee si on choisit $i$ est assez grand pour que $h(\delta_2)
  \subset O_0$ (ce qui est possible puisque, dans ${\cal V}_1$, $\lim_{i\rightarrow
  +\infty}\Delta^{+}_i=\partial^{+}A$, donc $\lim_{i\rightarrow
  +\infty}\delta_2=\emptyset$).  La droite $\delta$ est donc
   une droite de \hbox{Brouwer} simpliciale \`a extr\'emit\'es Nord-Sud ;  d'apr\`es le point
  \ref{ite.NS}, on a alors $\delta=\partial^{+}A$, c'est-\`a-dire
  $\Delta^{++}_i=\emptyset$.
\end{demo}

\item 
\begin{affi}
\label{aff.fina}
 On a $\Delta_\infty=\partial^-A \cup \partial ^+A$.
\end{affi}
\begin{demo}
Les composantes connexes de $\Delta_\infty$
sont des droites topologiques ; on raisonne par l'absurde en supposant
que l'une d'entre elles, not\'ee $\delta_\infty$, est distincte (et donc
disjointe) de
$\partial^-A$ et de $\partial ^+A$.

Gr\^ace aux deux affirmations pr\'ec\'edentes (\ref{aff.escr} et
\ref{aff.de++}), on voit que pour tout $i$ assez grand, $A^+(B_i) \cup
\{N\}$ est un voisinage de $N$ dans $\hat A$~; par cons\'equent, la seule
possibilit\'e est que $\delta_\infty$ 
soit \`a extr\'emit\'es Sud. De plus, elle est libre (affirmation \ref{aff.comp}).
On appelle $P_\infty$ le disque topologique ferm\'e de fronti\`ere
$\delta_\infty \cup \{S\}$ inclus dans $\hat A$.

Montrons que $P_\infty$ est un p\'etale r\'epulsif, autrement dit que
$\delta_\infty$ est  une droite de \hbox{Brouwer}.
Nous savons d\'ej\`a que c'est une droite topologique libre, il s'agit
donc de prouver que la condition de s\'eparation est respect\'ee. Soit $x$
un point arbitraire de $\delta_\infty$, et notons
$\delta_\infty^-$ et $\delta_\infty^+$ les deux arcs  d\'ecoup\'es par les
deux points $x$
et $S$ sur la courbe de Jordan $\delta_\infty \cup \{S\}$, de mani\`ere
\`a ce que les arcs $\partial^-A \cup \{N,S\}$, $\delta_\infty^-$, $\delta_\infty^+$
 et $\partial^+A \cup \{N,S\}$ soient dispos\'es dans cet ordre  autour
de $S$ (figure~\ref{fig2}, gauche).
\begin{figure}[htp]
  \par 
\centerline{\hbox{\input{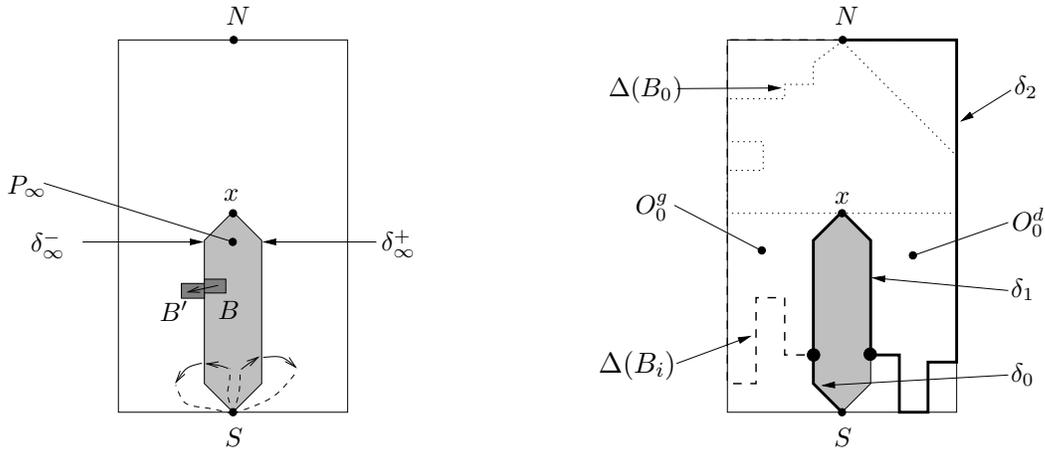}}} 
\par
\caption{\label{fig2}Construction de $\delta$}
\end{figure}

Consid\'erons deux briques adjacentes $B \subset P_\infty$
 et $B'\not \subset P_\infty$,  situ\'ees de part et d'autre de $\delta_\infty$.
Pour $i$ assez grand, l'attracteur $A^+(B_i)$ contient $B'$ mais ne contient pas
$B$. On ne peut donc pas avoir $H(B') \cap B
\neq \emptyset$ ; par minimalit\'e de la d\'ecom\-po\-sition en briques, on a  alors $H(B) \cap
B' \neq \emptyset$. 

En appliquant cet argument \`a des briques rencontrant
 $\delta_\infty^-$ ou $\delta_\infty^+$
 et arbitrairement proches de $S$, on obtient que les arcs
 $\partial^-A \cup \{N,S\}$, $H(\delta_\infty^-)$,
 $\delta_\infty^-$, $H^{-1}(\delta_\infty^-)$,$H^{-1}(\delta_\infty^+)$  
 $\delta_\infty^+$, $H(\delta_\infty^+)$ 
 et $\partial^+A \cup \{N,S\}$ (qui sont disjoints deux \`a deux) 
sont dispos\'es dans  cet ordre autour
de $S$.
Ceci prouve que $\delta_\infty$ s\'epare son image par $H$ de sa
 pr\'eimage : c'est bien une droite de \hbox{Brouwer} ; on voit \'egalement que
 son type dynamique est $\sf$.

Nous allons conclure la preuve de l'affirmation \ref{aff.fina} par un
argument similaire \`a celui utilis\'e dans la preuve de l'affirmation
\ref{aff.de++} (figure~\ref{fig2}, droite).
Comme $P_\infty$ est un p\'etale r\'epulsif simplicial,
  on a $\Delta(B_i) \cap
\inte(P_\infty)=\emptyset$ pour tout $i$.
Consid\'erons \`a nouveau un point $x$ de $\delta_\infty$ ; quitte \`a r\'e-indexer, on peut
supposer que $\Delta(B_i)$ contient $x$ pour tout $i \geq 0$. Soit $i$
assez grand (voir plus loin).  On note (figure \ref{fig2}) :
\begin{itemize}
\item $\delta_1$ la composante connexe de $\Delta(B_i) \cap
\delta_\infty$ contenant $x$ ;
\item $\delta_2$ la composante de $\Delta(B_i) \setminus \delta_1$
rencontrant $\partial^{+}A$ ;
\item $\delta_0=\delta_\infty^- \setminus \delta_1$.
\end{itemize}
Posons encore $\delta'_0=\delta_0 \cup \delta_1$, $\delta'_2=\delta_1
\cup \delta_2$, et $\delta=\delta_0 \cup \delta_1 \cup
\delta_2=\delta'_0 \cup \delta'_2$.

Montrons que $\delta$ est une droite de \hbox{Brouwer}. Les arcs $\delta'_0$
 et $\delta'_2$ sont libres, donc (comme au point \ref{ite.N}) il
 suffit de montrer que $H(\delta'_0) \cap \delta'_2=\emptyset$ et
 $H^{-1}(\delta'_0) \cap \delta'_2=\emptyset$. La deuxi\`eme \'egalit\'e
 est toujours vraie car $P_\infty$ est un r\'epulseur strict, donc
 $H^{-1}(\delta_0) \subset \inte(P_\infty)$. Montrons la premi\`ere~:
 l'ensemble $O_0$ a \'et\'e d\'efini au point~\ref{ite.ozero}~;
on \'ecrit $O_0 \setminus P_\infty$ comme l'union disjointe de
 $O_0^{g}$ et $O_0^{d}$ (points situ\'es ``\`a gauche'' et ``\`a
 droite'' de $P_\infty$, figure \ref{fig2}). D'apr\`es le point
 \ref{ite.escr2}, pour tout $i \geq i_0$, $\delta_2 \cap
 O_0^{g}=\emptyset$. Mais si $i$ tend vers $+\infty$, $\delta_1$ tend
 (dans ${\cal V}_1$)
 vers $\delta_\infty=\partial P_\infty$, donc $\delta_0$ tend vers $\emptyset$
 ; comme $P_\infty$ est un r\'epulseur, on peut choisir $i$ pour que
 $H(\delta_0) \subset O_0^{g}$, et dans ce cas l'\'egalit\'e est
 v\'erifi\'ee.

On a alors trouv\'e une droite de \hbox{Brouwer} \`a extr\'emit\'es  Nord-Sud
 $\delta$
 dans $A$ qui
n'est \'egale \`a aucun des deux bords de $A$, ce qui contredit
 l'affirmation \ref{aff.brns}.
\end{demo}

\item
 Les affirmations \ref{aff.de++} et \ref{aff.fina} entra\^inent
que pour tout voisinage $O_S$ de 
$S$, pour $i$ assez grand, $\Delta(B_i) \setminus \partial A \subset
O_S$. Ceci prouve que $A$ est un croissant attractif \`a
dynamique Sud-Nord.
\end{enumerate}
\end{demo}

\section{Indice partiel et construction de droites de
\hbox{Brouwer} (proposition \ref{pro.codb2})}
\label{sec.brsu}
Dans cette section, on suppose que  $H$ est un hom\'eo\-morphisme de
\hbox{Brouwer}, et que $F$ est une d\'ecom\-po\-sition en briques pour $H$.

\subsection{Indices partiels et d\'ecom\-po\-sition en briques}
\label{sse:ipdb}
\begin{figure}[htp]
\par
\centerline{\hbox{\input{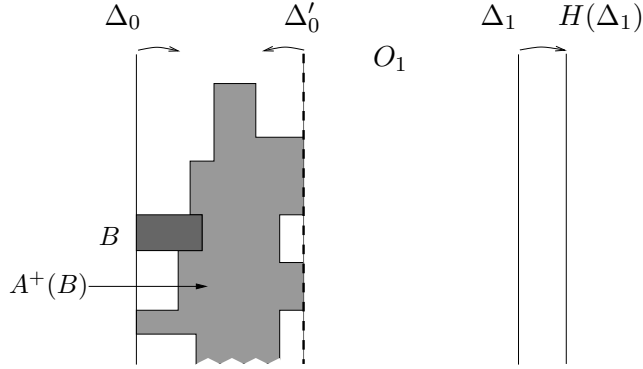}}}
\par
\caption{\label{fig-contre-courant}Obtention d'une droite de \hbox{Brouwer} \`a contre-courant}
\end{figure}

La formule contenue dans la proposition~\ref{pro.codb2} affirme
l'existence d'un certain nombre de croissants.
Le lemme suivant est un morceau de la proposition (obtenu
dans le cas particulier o\`u la famille $\cal F$ est vide).
 Il nous permettra de prouver cette
proposition en construisant, l'une apr\`es l'autre, les droites de
 \hbox{Brouwer} qui bordent les croissants.
\begin{lemm}[figure \ref{fig-contre-courant}]
\label{lem.cocr}
Soit $H$ un hom\'eo\-mor\-phisme de \hbox{Brouwer},
$F$ une d\'ecom\-po\-si\-tion en briques pour $H$, et
$(\Delta_0,\Delta_1)$ un couple indiff\'erent de  droites de
\hbox{Brouwer} simpliciales disjointes, $\Delta_0$ \'etant attractive.

Supposons que l'une des deux hypoth\`eses suivantes est v\'erifi\'ee~:
\begin{enumerate}
\item l'indice partiel $\inpa(H,\Delta_0,\Delta_1)$ est diff\'erent de
 $0$~;
\item il n'existe pas de
cha\^\i ne de pseudo-disques de $\Delta_0$ \`a $\Delta_1$, constitu\'ee
 de briques de la d\'ecom\-po\-si\-tion $F$.
\end{enumerate}
 Alors il existe une droite de \hbox{Brouwer} simpliciale
$\Delta'_0$ disjointe des droites topologiques $\Delta_0$ et $\Delta_1$, les s\'eparant, telle
que le couple $(\Delta_0,\Delta'_0)$ est attractif.
\end{lemm}

\begin{demo}[du lemme~\ref{lem.cocr}]
Rappelons que d'apr\`es la proposition \ref{pro.frip} (``lemme de
Franks'' pour l'indice partiel), l'hypoth\`ese~1 du lemme
 implique l'hypoth\`ese~2. Il suffit donc
d'\'ecrire la preuve sous  la deuxi\`eme hypoth\`ese du lemme. 

On note $O_1$ le disque ouvert $D(\Delta_0,\Delta_1)$ (figure~\ref{fig-contre-courant}).
 Soit $\cal B$ l'ensemble des briques de $F$  incluse dans $\adhe(O_1)$ et
 rencontrant $\Delta_0$. Soit  $B \in {\cal B}$.  Puisque $(B)$ est une
 cha\^\i ne de briques (de longueur 1), l'hypoth\`ese~2
 du lemme implique que $B \cap \Delta_1=\emptyset$.
D'apr\`es la
 caract\'erisation de $A^+(B)$ en termes de cha\^ine de briques
 (affirmation \ref{aff.chai}),
on a aussi   $A^{+}(B) \cap \Delta_1=\emptyset$. Puisque $\Delta_0$
 est attractive, on a $A^+(B) \subset P^+(\Delta_0)$~; 
on en d\'eduit que $A^+(B)$ est inclus dans $O_1\cup
 \Delta_0$.
D'autre part, remarquons que l'ensemble  $B \cup A^{+}(B)$ est encore
 un attracteur strict puisque, par d\'efinition de $A^+(B)$, on a $h(B)
 \subset \inte(A^+(B))$.   Posons alors
$$
A=\bigcup_{B \in {\cal B}} B \cup A^{+}(B).
$$
Cet ensemble est inclus dans $O_1 \cup \Delta_0$. C'est un attracteur
 strict,  puisqu'il est r\'eunion
d'attracteurs stricts, et il est connexe. D'apr\`es l'affirmation \ref{aff.drto},
les composantes connexes de $\partial A$ sont des  droites
topologiques. Soit $\Delta'_0$
l'unique   composante connexe de $\partial A$ qui s\'epare $\Delta_0$
et $\Delta_1$ (on peut par aussi d\'efinir $\Delta'_0$ comme la
 fronti\`ere de
la composante connexe de $\R^2 \setminus A$ qui contient $\Delta_1$).
Cette droite topologique est disjointe de $\Delta_0$ et de
$\Delta_1$, et elle est libre puisque $A$ est un attracteur strict.
 
V\'erifions qu'elle s\'epare son image et sa
pr\'eimage. Tout d'abord, remarquons que $\Delta_1$ s\'epare $\Delta_0$
et $H(\Delta_1)$ (car $\Delta_1$ est r\'epulsive pour
le couple $(\Delta_0,\Delta_1)$) ; comme $\Delta'_0$ s\'epare $\Delta_0$
et $\Delta_1$, on peut supposer que $\Delta_0$, $\Delta'_0$,
$\Delta_1$, $H(\Delta_1)$ sont quatre droites verticales dans le plan, dont les
abscisses sont rang\'ees dans cet ordre (\`a l'aide du th\'eo\-r\`eme de
Schoenflies-Homma, voir l'appendice). \index{Schoenflies-Homma!th\'eo\-r\`eme de}
 En particulier, $\Delta'_0$ ne
s\'epare pas $\Delta_1$ et $H(\Delta_1)$.
Soit maintenant  $U$ le demi-plan ouvert  de fronti\`ere $\Delta'_0$ qui
contient $\Delta_1$.  L'attracteur $A$ ne rencontre pas $U$, donc
$H(\Delta'_0)$, qui est inclus dans $A$, est disjoint de $U$.
Autrement dit, $U$ ne
rencontre pas la fronti\`ere de $H(U)$, et comme $U$ est connexe, $U$
est inclus dans $H(U)$ ou dans son compl\'ementaire. Mais d'autre part
$H(U)$ contient $H(\Delta_1)$, qui est aussi dans $U$, donc on a en
fait $U \subset H(U)$. 

L'ensemble $\adhe(U)$ est donc un p\'etale
r\'epulsif, et sa fronti\`ere $\Delta'_0$ est une droite de \hbox{Brouwer}
d'apr\`es la remarque~\ref{rem.attr}.
 Comme $A$ est un
attracteur, le couple $(\Delta_0,\Delta'_0)$ est attractif.
\end{demo}

On peut remarquer que dans cette construction, il existe une cha\^\i
ne de briques de $\Delta_0$ \`a $\Delta'_0$, et que par cons\'equent l'indice
partiel entre ces deux droites vaut $+1/2$ ou $-1/2$ (d'apr\`es le
``lemme de Franks pour l'indice partiel'', proposition \ref{pro.frip}).

\incn[existence : on peut partir de $\t \Delta_0'$, suivre un chemin
d'ar\^etes vers $\t \Delta_0$ jusqu'au premier moment o\`u on touche
$\partial A$. Soit $\t \Delta_1$ la composante de $\partial A$ qu'on a
rencontr\'e. Elle doit s\'eparer.
Propri\'et\'e de s\'eparation : ben, c'est pas compl\`etement \'evident. Utiliser
que l'on sait o\`u sont $H(\t \Delta_0)$ et $H(\t \Delta_0')$ (les 6
droites sont disjointes deux \`a deux, on peut
les supposer verticales), que $H(\t \Delta_1)$ s\'epare $H(\t \Delta_0)$
et $H(\t \Delta_0')$, donc est aussi verticale, $H$ ne peut pas
renverser $\t \Delta_1$. 
]

\subsection{Indice partiel et type des croissants}
\begin{demo}[de la proposition~\ref{pro.codb2}]
Soit $(\Delta_0,\Delta_1)$ un couple de droites de \hbox{Brouwer}
simpliciales disjointes, et
 $\cal F$ une famille maximale de sous-croissants simpliciaux minimaux
de la bande $\adhe(D(\Delta_0,\Delta_1))$.
\paragraph{Premier point}\index{droites topologiques!qui se traversent}
Soit $i$ un entier entre $-1$ et $k$. Puisque les int\'erieurs des
croissants sont disjoints deux \`a deux, les droites topologiques 
$\Delta'_{2i+1}$ et $\Delta'_{2i+2}$, appartenant \`a deux croissants
successifs, ne se traversent pas (d\'efinition
\ref{def.trav}), et l'indice partiel est bien d\'efini
(section~\ref{ss.inpat}). Dans le cas o\`u elles ne sont
pas disjointes, il vaut $0$ par d\'efinition.

Dans le cas o\`u elles sont disjointes, $(\Delta'_{2i+1},\Delta'_{2i+2})$ est un couple indiff\'erent
de droites de \hbox{Brouwer} simpliciales. On se place par exemple
 dans le cas o\`u $\Delta'_{2i+1}$ est attractive. Supposons par l'absurde que
 l'indice partiel
n'est pas nul, ou qu'il n'existe pas de cha\^{\i}ne de pseudo-disques de
$\Delta'_{2i+1}$ \`a $\Delta'_{2i+2}$ constitu\'ee de briques. On peut alors appliquer le
lemme~\ref{lem.cocr} au couple $(\Delta'_{2i+1},\Delta'_{2i+2})$, et
on obtient une droite de \hbox{Brouwer} $\Delta$, s\'eparant
$\Delta'_{2i+1}$ et $\Delta'_{2i+2}$, la bande
$\adhe(D(\Delta'_{2i+1},\Delta))$ \'etant un croissant attractif
simplicial dont l'int\'erieur est disjoint des autres croissants de la
famille $\cal F$. D'apr\`es la
proposition~\ref{pro.cam}, ce croissant contient un sous-croissant
attractif simplicial minimal $C$. L'existence de $C$ 
contredit la maximalit\'e de la famille $\cal F$.

\paragraph{Deuxi\`eme point}
Soit $C=\adhe(D(\Delta'_{2i},\Delta'_{2i+1}))$ un croissant de $\cal F$.
D'apr\`es la proposition~\ref{pro.cam}, le croissant $C$, \'etant minimal, est 
 \`a dynamique $\t S$-$\t N$ ou $\t N$-$\t S$. L'indice partiel \`a
travers le croissant est alors donn\'e par l'affirmation suivante~:
\begin{affi}\label{aff.ipcr}
Si $\adhe(D(\Delta,\Delta'))$ est un croissant attractif \`a dynamique  $\t
 N$-$\t S$, ou r\'epulsif \`a dynamique  $\t S$-$\t N$, alors l'indice partiel 
 $\inpa(\Delta,\Delta')$ vaut $-1/2$. Dans les deux autres cas,
 l'indice partiel vaut $+1/2$. 
\end{affi}
 \begin{demo}[de l'affirmation]
On suppose par exemple que le croissant est attractif. D'apr\`es la
d\'efinition~\ref{def.cdns} des croissants \`a dynamique  $\t N$-$\t S$ ou $\t S$-$\t N$,
il existe un p\'etale attractif, bas\'e en $\t S$ ou en $\t N$,
 dont le bord rencontre les deux bords du croissant
$D(\Delta,\Delta')$. Soit $\gamma$ 
un sous-arc du bord du p\'etale, rencontrant les deux bords du croissant, et
minimal pour l'inclusion. L'arc $\gamma$ est libre. Si le croissant  est 
 \`a dynamique $\t S$-$\t N$, $H(\gamma)$ est au-dessus de $\gamma$, et 
s'il  est  \`a dynamique $\t N$-$\t S$, $H(\gamma)$ est au-dessous de $\gamma$.
Le lemme~\ref{lem.arli} montre alors que l'indice partiel
 $\inpa(\Delta,\Delta')$ vaut $1/2$ dans le premier cas,
 $-1/2$ dans le second.
 \end{demo}

On obtient bien s\^ur l'indice partiel \`a travers les croissants
r\'epulsifs \`a dynamique  $\t S$-$\t N$ ou $\t N$-$\t S$ en appliquant
l'affirmation \`a $H^{-1}$.

 Terminons maintenant la preuve du deuxi\`eme point. Nous connaissons
l'indice partiel entre deux droites successives quelconques de la famille
$(\Delta'_{-1}, \dots, \Delta'_{2k+2})$.

Puisque les int\'erieurs des croissants sont deux \`a deux disjoints, pour
 tout $i$, les deux droites $\Delta'_i$ et $\Delta'_{i+1}$ ne se
 traversent pas. De plus, pour tous $i<j<k$, la droite $\Delta'_{j}$ s\'epare $\Delta'_{i}$ et
 $\Delta'_{k}$ (au sens de la d\'efinition~\ref{def.sepa2}).
\index{relation de Chasles!g\'en\'eralis\'ee}
 On peut donc appliquer $2k$ fois la relation de Chasles g\'en\'eralis\'ee
 (lemme~\ref{lem.chas2}), ce qui donne la formule voulue.
\end{demo}

\subsection{Preuve du th\'eor\`eme sur les hom\'eo\-mor\-phismes de \hbox{Brouwer}}

\begin{demo}[du th\'eor\`eme~\ref{the.prinbrou}] 
\`{A} partir de la
proposition~\ref{pro.codb2}, il suffit de recopier les arguments de la
fin de la preuve du th\'eor\`eme principal (section~\ref{sec.preu},
paragraphe intitul\'e ``le compte y est'').  Les d\'etails de la preuve sont laiss\'es
au lecteur. Ajoutons que cet \'enonc\'e n'est pas utilis\'e dans le reste du texte.
\end{demo}

\section{Existence d'une droite de \hbox{Brouwer} \`a extr\'emit\'es 
Nord-Sud (Proposition~\ref{pro.dbns})}
\label{sec.dbns}
Soit $h$ un hom\'eo\-morphisme de  la
sph\`ere $\S^2$, pr\'eservant l'orientation, avec $\fixe(h)=\{N,S\}$
(\textbf{hypoth\`ese (H2)}). On suppose  que $\indi(h,N) \neq 1$.
 Soit $F$ une d\'ecom\-po\-sition en briques pour $h$.

Les notions de \res{types topologiques}, de \res{types dynamiques},
 ainsi que
 \res{l'espace ${\cal V}_1$ des sous-vari\'et\'es  simpliciales} ont \'et\'e
d\'efinis dans le contexte d'un croissant attractif 
simplicial $A$ pour un hom\'eo\-morphisme de \hbox{Brouwer} $H$ (sections~\ref{ss.type} et
\ref{sub.svs}).
Ces d\'efinitions et les r\'esultats de ces sections se transposent
sans difficult\'e, et quasiment mot pour mot, au contexte de la pr\'esente
section (remplacer $H$ par $h$, le plan $\R^2$ et $A$ par $\S^2 \setminus \{N,S\}$, et le
compactifi\'e  $\hat A$ par $\S^2$).

\pagebreak
\begin{demo}[de la proposition~\ref{pro.dbns}]
\paragraph{Premier cas} 
Si il y a  une brique $B$ de type topologique Nord-Sud, on a gagn\'e.
\paragraph{Deuxi\`eme cas} Supposons maintenant qu'il existe
des droites de \hbox{Brouwer} $\Delta(B)$ de type topologique Sud arbitrairement proches de $N$.

Autrement dit, il existe 
 une suite  $(B_k)$ de briques de
   type topologique Sud et  une suite $(x_k)$ de sommets de $F$
   tendant vers $N$ telle que pour chaque $k$, $x_k \in \Delta(B_k)$.
La suite $(\Delta(B_k))$ v\'erifie alors les hypoth\`eses du
lemme~\ref{lem.svs}.\footnote{Ce lemme a \'et\'e \'enonc\'e dans le contexte
d'un croissant attractif pour un hom\'eo\-morphisme de \hbox{Brouwer}, mais
l'\'enonc\'e et la preuve se transposent quasiment mot pour mot au pr\'esent contexte
(le seul changement, dans la preuve, est que $R$ est maintenant un anneau et non un
disque topologique).}

Il existe donc une droite topologique $\Delta$ de type
Nord-Sud, incluse dans une valeur d'adh\'erence de la suite
$(\Delta(B_k))$. Cette droite est libre (affirmation~\ref{aff.comp}),
c'est donc une droite de \hbox{Brouwer} \`a extr\'emit\'es Nord-Sud~: on a encore gagn\'e.

\paragraph{Troisi\`eme cas} S'il existe
des droites de \hbox{Brouwer} $\Delta(B)$ \`a extr\'emit\'es Nord arbitrairement proches
de $S$, on conclut comme au deuxi\`eme cas en \'echangeant les r\^oles de
$N$ et $S$.

\paragraph{Et c'est tout} Nous allons montrer qu'on est forc\'ement
dans l'un des trois cas \'enum\'er\'es.

On raisonne par l'absurde, en supposant le contraire. Il existe alors
 un voisinage connexe
$O_N$ de $N$ qui ne rencontre pas de droite $\Delta(B)$ \`a extr\'emit\'es
 Sud, et  un voisinage connexe $O_S$ de $S$ qui ne rencontre pas de droite
$\Delta(B)$ \`a extr\'emit\'es Nord.
Remarquons que si $B$ est de type dynamique $\fn$, $\Delta(B)$ s\'epare
 $A^+(B)$ de  $S$, donc aussi de $O_S$~; en particulier, $A^+(B)$
 et $O_S$ sont disjoints.

\begin{affi}
Il existe  une
famille $\{B_1, \cdots, B_k\}$ de briques qui sont toutes du m\^eme
type dynamique ($\sf$ ou $\nf$), telles que $C=B_1 \cup
\cdots \cup B_k$ est connexe et s\'epare $N$ et $S$.
\end{affi}
\begin{demo}[de l'affirmation]
 On pose
$$
E_N=\bigcup \{B \cup A^{+}(B) \mid B \mbox{ brique de type dynamique } \fn \},
$$
$$
E_S=\bigcup \{B \cup A^{+}(B) \mid B \mbox{ brique de type dynamique } \fs \}.
$$

Puisque $E_N$ est une r\'eunion de briques et ne rencontre pas $O_S$,
$E_N \cup \{N\}$ est ferm\'e. Cet ensemble est aussi connexe (comme
r\'eunion de connexes contenant $N$). De plus, $E_N \cap
E_S=\emptyset$ (voir l'affirmation \ref{aff.brfn} sur
``l'incompatibilit\'e'' des types $\fn$ et $\fs$).
  Il n'est pas difficile de voir qu'il existe alors une
famille $\{B_1, \cdots, B_k\}$ de briques telles que $C=B_1 \cup
\cdots \cup B_k$ est connexe, ne rencontre pas $\inte(E_N \cup E_S)$
et s\'epare $N$ et $S$ (figure \ref{fig21} : on peut construire
$C$ en ``faisant le tour'' de $O_N \cup E_N$).
\begin{figure}[hbtp]
  \par
\centerline{\hbox{\input{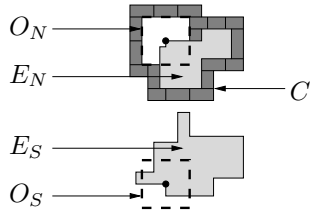}}}
\par
  \caption{\label{fig21}Construction de $C$}
\end{figure}
 Les briques $B_i$ ne peuvent pas \^etre de type Nord-Sud, sans quoi on
serait dans le premier cas ; ni de type dynamique $\fn$ ou $\fs$, sinon elles
seraient incluses dans $E_N$ ou $E_S$. Elles sont donc
chacune de type dynamique $\sf$ ou $\nf$ ; comme ces deux types dynamiques sont
``incompatibles'' (affirmation \ref{aff.brfn2}) et que $C$ est connexe,
 elles sont toutes du m\^eme
type dynamique.
\end{demo}

\begin{figure}[htp]
\par
\centerline{\hbox{\input{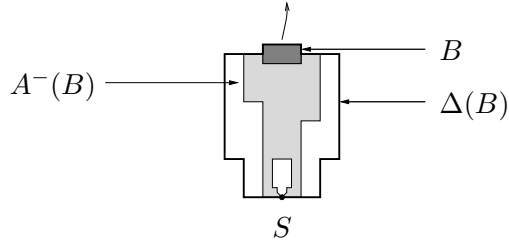}}}
\par
\caption{\label{fig40}L'ensemble $A^-(B)$ pour une brique de type
dynamique $\sf$}
\end{figure}
 On suppose par exemple que toutes les briques $B_i$ sont de type dynamique
$\sf$.\footnote{Ici, on conclurait facilement avec un lemme de Birkhoff, (voir
par exemple \cite{lero9}), mais un petit argument permet de l'\'eviter.}
Remarquons d'abord que dans ce cas, l'adh\'erence du r\'epulseur
$A^-(B_i)$~\footnote{D\'efinition~\ref{def.repu}.}
 contient $S$ et pas $N$ (figure \ref{fig40}).

Soit $C'=A^{-}(B_1) \cup \cdots \cup A^{-}(B_k)$. Cet ensemble est
 connexe,  c'est une r\'eunion
de r\'epulseurs, donc un r\'epulseur. De plus, $\adhe(C')$ contient
$S$ et pas $N$.

On consid\`ere la composante connexe de $\S^2 \setminus C'$ qui contient
$N$ ; soit $D$ son adh\'erence. Comme $C'$ est un r\'epulseur, on a
$h^{-1}(\partial D) \cap \inte(D)=\emptyset$ ; donc $D$ est un
attracteur. De plus, on voit facilement que c'est un disque
topologique. D'apr\`es le lemme \ref{lem.diat}, ceci contredit le fait
qu'il n'existe pas de courbe d'indice~$1$.
\index{disque attractif}
\end{demo}

Remarquons que le deuxi\`eme cas envisag\'e, dans lequel l'absence de
brique de type topologique Nord-Sud oblige
\`a effectuer un passage \`a la limite, survient pour la d\'ecom\-po\-si\-tion de
l'hom\'eo\-mor\-phisme lin\'eaire hyperbolique selle
dessin\'ee sur les figures \ref{fig11} et \ref{fig18}.

\section{Construction du relev\'e canonique (proposition \ref{pro.rele})}
\label{sec.reca}

\subsection{Preuve de la proposition~\ref{pro.rele}}

Soit $h$ un hom\'eo\-morphisme de  la
sph\`ere $\S^2$, pr\'eservant l'orientation, avec $\fixe(h)=\{N,S\}$
(\textbf{hypoth\`ese (H2)}). On suppose  que $\indi(h,N) \neq 1$.

\paragraph{Construction du relev\'e $\t h$}
  D'apr\`es la proposition~\ref{pro.dbns}, il existe au moins une droite de
\hbox{Brouwer} $\Delta_0$ dont l'adh\'erence contient $N$ et $S$.  

\begin{figure}[h!tpb]
\par
\centerline{\hbox{\input{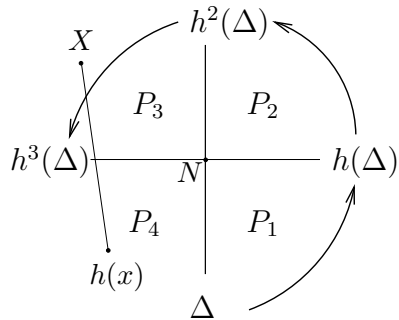}}}
\par
\caption{\label{fig42}L'image d'un point $X$ n'est jamais sur la
demi-droite issue de $N$, oppos\'ee \`a $X$}
\end{figure}

On se place dans le mod\`ele du plan.
Commençons par remarquer que les quatre droites
topologiques  $\Delta_0$, $h(\Delta_0)$, $h^{2}(\Delta_0)$ et
$h^{3}(\Delta_0)$ sont dispos\'e dans cet ordre cyclique autour du point
$N$. Ceci d\'ecoule facilement du fait que le domaine fondamental
$D(\Delta_0,h(\Delta_0))$ est libre (affirmation~\ref{aff.dofo}).
Quitte \`a conjuguer $h$
en utilisant le th\'eor\`eme de Schoenflies-Homma, (voir l'appendice) 
\index{Schoenflies-Homma!th\'eo\-r\`eme de}
on peut supposer que ces quatre droites topologiques
 sont des demi-droites euclidiennes issues de $N$
d\'ecoupant quatre quart-de-plans, dispos\'ees comme sur la figure
\ref{fig42}.

Notons $P_1=D(\Delta_0,h(\Delta_0)), \cdots,
 P_4=D(h^3(\Delta_0),\Delta_0)$ les quatre quart-de-plans comme sur la
 figure \ref{fig42}. L'ensemble $P_1$ \'etant un disque topologique ouvert libre, il est
 disjoint de tous ses it\'er\'es (lemme \ref{lem.libr})~:
 on a donc $h(P_1) = P_2$, $h(P_2)=P_3$, $h(P_3)
 \subset P_4$ et $h(P_4) \subset P_4 \cup P_1$.
 Notons $[AB]$ le segment euclidien, dans le mod\`ele du plan, joignant les deux points $A$
et $B$. Dans cette situation, pour tout $X \in \R^{2}$, on a la
propri\'et\'e  $N \not \in [X h(X)]$, autrement dit le
 segment  $[Xh(X)]$ est inclus dans l'anneau $\A$.

On consid\`ere~:
\begin{itemize}
\item $X_0$ un point quelconque de l'anneau $\A$~;
\item $\gamma=[X_0h(X_0)]$~;
\item $\t X_0$ un relev\'e de $X_0$ \`a $\t \A$~;
\item $\t \gamma$ l'arc relev\'e de $\gamma$ issu de $\t X_0$~;
\item $\t Y_0$ l'autre extr\'emit\'e de $\t \gamma$.
\end{itemize}
 D'apr\`es les propri\'et\'es
g\'en\'erales des rev\^etements, il existe un unique relev\'e $\t h$ de $h$
v\'erifiant $\t h(\t X_0)=\t Y_0$.

Comme $\t h$ commute avec l'automorphisme de rev\^etement $\tau$, il
ne d\'epend pas du choix du relev\'e $\t X_0$. Comme l'arc $[Xh(X)]$ varie
contin\^ument avec $X$, le relev\'e $\t h$ ne d\'epend pas non plus du
choix du point $X_0$.

\paragraph{Calcul des indices}
Notons $\t \theta$ l'application $(\t \theta,r) \mapsto \t \theta$ du
plan $\t \A$. Par construction, l'hom\'eo\-morphisme $\t h$ v\'erifie la
propri\'et\'e~:
$$
(*)  \mbox{ Pour tout } \t X \in \t \A, \quad \quad |\t \theta (\t X)-\t \theta
(\t h(\t X))| < \frac{1}{2}. 
$$
La formule reliant les indices d\'ecoule alors du lemme suivant~:
\begin{lemm}\label{l.cones}
Soit $g$ un hom\'eo\-morphisme de l'anneau $\A$, sans point fixe, isotope
\`a l'identit\'e, et $\t g$ un relev\'e de $g$ au plan $\t \A$. Supposons
que $\t g$ v\'erifie la propri\'et\'e $(*)$ ci-dessus. Alors $\indi(\t g,\tau)=\indi(g,N)-1$.
\end{lemm}
(La preuve du lemme~\ref{l.cones} est donn\'ee plus loin.)

\paragraph{Unicit\'e} 
Le lemme pr\'ec\'edent prouve notamment que l'indice de $\t h$ par rapport
\`a $\tau$ n'est pas nul.
Soit maintenant $\tilde h'=\tau^{k} \circ \tilde h$ un
autre relev\'e de $h$ (pour fixer les id\'ees, prenons $k > 0$, le cas
n\'egatif \'etant similaire), et $\t X$ un
point  de $\t \A$.  
De la propri\'et\'e $(*)$, on d\'eduit l'in\'egalit\'e
$$(**)  \quad  \t \theta (\tilde h' (\t X))= \t \theta (\tilde h (\t X))+k \geq \t
\theta (\tilde h(\t X))+1 > \t \theta(\t X) + 1/2$$
 ce qui entra\^ine facilement que l'indice de
$\tilde h'$ par rapport \`a $\tau$ est nul (``le vecteur
 $\t h'(\t X) - \t X$ pointe partout
vers la droite, il ne peut pas faire de tours''~: voir par exemple la
preuve du lemme \ref{lem.arli}).

\paragraph{Point 1~: relev\'e des droites de \hbox{Brouwer}}
Un relev\'e d'une droite topologique de $\S^2$ est clairement une droite
topologique de $\t \A$. Un relev\'e d'un ensemble libre est clairement
libre. Il reste \`a voir que tout relev\'e $\t \Delta$ d'une droite de \hbox{Brouwer} $\Delta$ dont
l'adh\'erence contient $N$ et $S$ v\'erifie la propri\'et\'e de s\'eparation~:
$\t \Delta$ s\'epare $\t h(\t \Delta)$ de $\t h^{-1}(\t \Delta)$.

Or cette propri\'et\'e est clairement v\'erifi\'ee si $\Delta$ est la droite
de \hbox{Brouwer} $\Delta_0$ utilis\'ee dans la construction de $\t h$ donn\'ee
ci-dessus. Par unicit\'e, l'hom\'eo\-morphisme $\t h$ ne d\'epend pas de la
droite utilis\'ee dans la construction, donc la propri\'et\'e est v\'erifi\'ee
pour toute droite $\Delta$.

\paragraph{Point 2~: projection des arcs libres}
Rappelons  que  le nombre $\indi(\t h, \tau)$ n'est pas nul.
Nous allons alors utiliser le lemme suivant (dont la preuve est donn\'ee
plus loin)~:
\begin{lemm}
\label{lem.indi}
  Soit $H \in \com(\tau)$, et supposons qu'il existe
un arc $\gamma$ dans $\t \A$, libre pour $H$, et un entier $k$ non nul 
tels que $\gamma(1)=\tau^k(\gamma(0))$.

Alors  $\indi(H,\tau)= 0$. 
\end{lemm}

Soit  $\gamma$ un arc de $\t \A$ libre pour $\t h$, montrons
 que $\pi(\gamma)$ est un arc libre pour $h$.
 La courbe $\pi(\gamma)$ est bien un arc : sinon l'arc $\gamma$ rencontre un de ses
 it\'er\'es  $\tau^k(\gamma)$, ce qui est exclu d'apr\`es le lemme
 \ref{lem.indi} ci-dessus (que l'on applique au sous-arc $\gamma'$ de
 $\gamma$ allant d'un point $x$ dans $\gamma \cap \tau^{-k}(\gamma)$
 au point $\tau^k(x)$).
 Si $\pi(\gamma)$ n'est pas libre,
 $\t h(\gamma)$ rencontre  $\tau^k(\gamma)$ pour un certain entier
 $k$ (avec $k \neq 0$ puisque $\gamma$ est libre pour $\t h$). 
Soit $x \in \gamma$ tel que $y=\t h ^{-1}(\tau^{k}(x))
 \in \gamma$. L'arc $\gamma$ est alors un arc libre pour $\t h$ contenant $x$ et
 $y$~: d'apr\`es le corollaire~\ref{cor.arli}, il existe aussi un
 arc libre pour $\t h$ qui contient $x$ et $\t h(y)=\tau^k(x)$. Ceci
 contredit encore le lemme \ref{lem.indi}.

\paragraph{Point~3~: relev\'e d'une d\'ecom\-po\-sition en briques}
Soit enfin $F$ une d\'ecom\-po\-sition en briques pour $h$, et $\t
F=\pi^{-1}(F)$~; il s'agit de montrer que $\t F$ est une d\'ecom\-po\-sition
en briques pour $\t h$.

Il est clair que $\t F$ est encore un graphe triadique.
Comme les briques de $F$ sont des disques topologiques, chaque brique
de $\t F$ est un disque topologique relev\'e d'une brique de $F$. En
particulier, le graphe triadique $\t F$ est libre et compact.
Il reste \`a prouver que $\t F$ est minimal, c'est-\`a-dire que la
r\'eunion de deux briques adjacentes n'est pas libre.

Soient $\t B_0$ et $\t B_1$ deux briques de $\t F$ adjacentes, et
appelons $B_0=\pi(\t B_0)$ et $B_1=\pi(\t B_1)$ les briques
correspondantes dans $F$. Les briques $B_0$ et $B_1$ sont adjacentes ;
donc l'image de l'une rencontre l'autre, par exemple il existe $x\in
B_0$ tel que $h(x) \in B_1$. Soit $\t x $ un relev\'e de $x$ dans $\t
B_0$, et $\t y$ un relev\'e de $h(x)$ dans $\t B_1$. Puisque $\t B_0$ et
$\t B_1$ sont adjacentes, il existe un arc $\gamma$ allant de $\t x$ \`a $\t y$
et inclus dans $\t B_0 \cup \t B_1$.

Or l'arc $\pi(\gamma)$ n'est pas libre, puisqu'il va de $x$ \`a $h(x)$~;
d'apr\`es le point pr\'ec\'edent de la proposition, l'arc $\gamma$ n'est pas
libre. \textit{A fortiori}, la r\'eunion $\t B_0 \cup \t B_1$ n'est pas
libre. Ceci prouve que $\t F$ est minimal, et ach\`eve la preuve de la
proposition~\ref{pro.rele}.

\paragraph{Transitions pour $\t F$}
Remarquons au passage que puisque $h(B_0) \cap B_1 \neq \emptyset$, 
le lemme de Franks (\ref{lem.frdi})
entra\^ine  $h(B_1) \cap B_0=\emptyset$,  donc $\t h(\t B_1)
\cap \t B_0=\emptyset$. La d\'ecom\-po\-sition $\t F$ \'etant minimale,
 c'est donc l'image de $\t B_0$ qui rencontre
$\t B_1$. Autrement dit, l'application $\pi$ envoie les ``transitions'' de $\t F$
sur les transitions de $F$ (ou, en terme de l'orientation des ar\^etes
d\'efinie dans \cite{sauz1}, est compatible avec les orientations).

\subsection{Preuve des deux lemmes}
\begin{figure}[htpb]
\par
\centerline{\hbox{\input{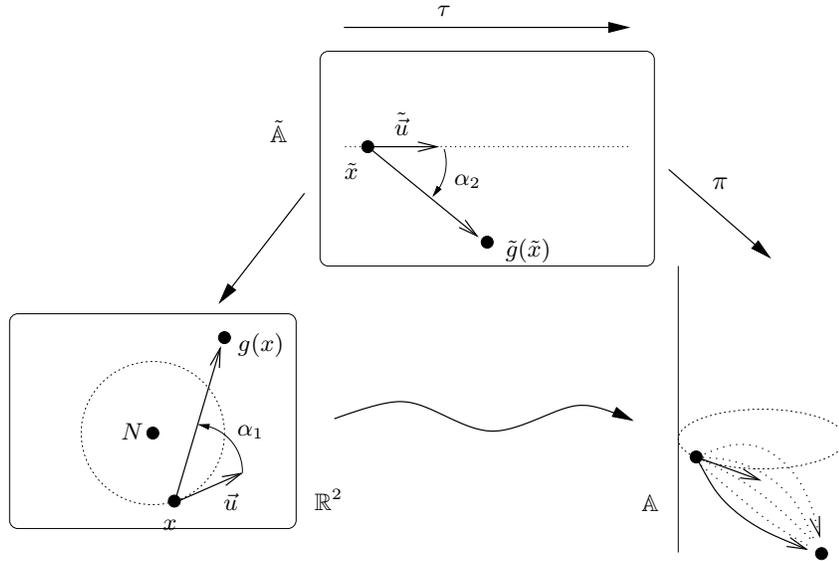}}}
\par
\caption{\label{fig48}Preuve du lemme~\ref{l.cones}}
\end{figure}
\begin{demo}[du lemme~\ref{l.cones}]
Voici une explication informelle, qui s'appuie sur la
figu\-re~\ref{fig48}.
On se place sous les hypoth\`eses du lemme. D'apr\`es la d\'efinition de
l'indice, le nombre $\indi(g,N)-1$
correspond alors, dans le mod\`ele du plan $\R^2$ (dessin en bas \`a
gauche), \`a la variation angulaire de l'angle $\alpha_1$ entre
le vecteur ``tournant'' $\vec u$ et la direction du segment euclidien de $x$ \`a $g(x)$,
lorsque le point $x$ parcourt 
un cercle qui entoure $N$ (le ``$-1$'' vient du fait que le vecteur $\vec u$ fait
un tour lorsque $x$ parcourt le cercle).

Le nombre $\indi(\t g,\tau)$, lui, correspond \`a la variation angulaire
de l'angle $\alpha_2$ repr\'esent\'e sur le dessin du haut, entre le
vecteur horizontal $\t{\vec u}$, relev\'e de $\vec u$,
 et la direction du segment de $\t x$ \`a $\t g (\t x)$,
lorsque le point $\t x$ parcourt un segment horizontal de longueur $1$
 (qui se projette sur le cercle).

Pour montrer que ces deux nombres n'en sont qu'un,
le probl\`eme vient du fait que les deux angles $\alpha_1$ (en un point
$x$) et $\alpha_2$ (en un relev\'e $\t x$) ne sont pas \'egaux (comme
on peut le voir  sur le dessin~: la projection de $\t \A$ sur $\A$ est une
isom\'etrie locale pour les m\'etriques naturelles sur les dessins,
 mais pas celle de $\t \A$ sur $\R^2$).
 
Cependant, gr\^ace \`a l'hypoth\`ese~(*),
ces deux angles ne prennent jamais des valeurs oppos\'ees. 
On en d\'eduit facilement que les deux
variations angulaires co\"{\i}ncident (une mani\`ere g\'eom\'etrique de voir
que l'application $\alpha_1$ est homotope \`a l'application $\alpha_2$
consiste \`a d\'eformer le plan $\R^2$ en l'anneau $\A$, en passant par
des c\^ones de plus en plus pointus~; \`a tout moment de la d\'eformation,
l'hypoth\`ese~(*) assure l'existence d'une unique g\'eod\'esique minimisante
de $x$ \`a $g(x)$~; au d\'ebut de la d\'eformation, cette g\'eod\'esique co\"{\i}ncide avec le
segment euclidien du plan $\R^2$~; \`a la fin, elle  co\"{\i}ncide avec le
segment g\'eod\'esique dans l'anneau $\A$, qui se rel\`eve en un segment
euclidien dans le plan $\t \A$).
\end{demo}

\begin{demo}[du lemme \ref{lem.indi}]
On se place sous les hypoth\`eses du lemme.
On consid\`ere, dans l'anneau $\A$, la courbe $\pi(\gamma)$ ; c'est
une courbe ferm\'ee, qui fait
$k$ fois le tour de $\A$. On voit facilement qu'une telle courbe
contient une courbe ferm\'ee simple, et que cette courbe ferm\'ee simple
fait une fois le tour de $\A$ (pas plus, sinon elle aurait un point
double ; et elle ne peut pas \^etre nulhomotope, sans quoi $\gamma$
aurait un point double). Cette courbe ferm\'ee simple poss\`ede un relev\'e
$\gamma'$ inclus dans $\gamma$, et $\gamma'$  est encore libre pour $H$~:
quitte \`a remplacer $\gamma$ par $\gamma'$,
 on peut supposer que $k=1$ et que $\pi(\gamma)$ est une courbe
ferm\'ee simple de $\A$.

On peut maintenant se ramener au cas o\`u  $\gamma=[0,1]\times\{0\}$ :
        en effet, le th\'eor\`eme de Schoenflies
        fournit un \homeo\ $g$ de $\A$, isotope \`a l'identit\'e,
        v\'erifiant $g(\pi(\gamma))=\pi([0,1]\times\{0\})$, et on
        remplace $H$ par $\t g H \t g^{-1}$ o\`u $\t g $ est un
        relev\'e de $g$ ; ceci ne change pas la valeur de
        $\indi(H,\tau)$ d'apr\`es l'affirmation~\ref{aff.comm}.

Rappelons que $\inte(\alpha)$ d\'esigne l'\emph{int\'erieur} d'un arc $\alpha$,
c'est-\`a-dire l'ensemble $\alpha(]0,1[)$.
 Soit, dans $\t \A$, les deux demi-droites $D_1=]-\infty,0] \times
        \{0\}$ et $D_2=[1,+\infty[ \times \{0\}$. Supposons que
        $H(\gamma)$ rencontre $D_1$ et $D_2$. On voit alors facilement
        que $\inte(\tau(H(\gamma)))$ doit rencontrer
        $\inte(H(\gamma))$ ; mais alors  $\inte(\tau(\gamma))$
        rencontrerait $\inte(\gamma)$, or ce n'est pas le cas.

 L'une au moins des deux demi-droites $D_1$ et $D_2$ est donc
        disjointe de $H(\gamma)$ ; supposons par exemple qu'il
        s'agisse de $D_1$. Alors, comme $\gamma$ est libre,
        $H(\gamma)$ est disjoint de $D_1 \cup \gamma$.
	Dans ce cas, pour tout $X\in \gamma$, le
        vecteur $H(X)-X$ ne pointe jamais dans la direction
        horizontale gauche. Ceci entra\^ine que $\indi(H,\tau)=0$
        (l'argument est le m\^eme que celui donn\'e dans la preuve du
        lemme \ref{lem.arli} et  repr\'esent\'e sur la figure \ref{fig27}).
\end{demo}

\subsection{Compl\'ements sur $\t h$ et sur les autres relev\'es de $h$} 
L'\'enonc\'e suivant contient quelques propri\'et\'es suppl\'ementaires des
relev\'es de $h$, qui ne seront pas utilis\'ees dans le texte.
Il n\'ecessite les d\'efinitions suivantes (et les d\'efinitions similaires
pour $\t h$)~:
\begin{defi}
\index{arc de translation}
Un \res{arc de translation} pour $h$ est un arc $\gamma$ qui relie un
point $x$ \`a son image $h(x)$, et tel que $\gamma \cap h(\gamma)=\{h(x)\}$.
\end{defi}
\begin{defi}
\index{singulier!couple}
Un couple $(x,y)$ de
points de $\A$  n'appartenant pas \`a la m\^eme orbite de $h$
est \res{singulier} si pour tout $\epsilon>0$, il
existe un point $z$ et un entier positif $n$, tel que $z$ soit
$\epsilon$-proche de $x$, et $h^n(z)$ $\epsilon$-proche de $y$.
\end{defi}
\begin{prop}[compl\'ement \`a la proposition~\ref{pro.rele}]\label{pro.comp}
Sous les hypoth\`eses de la proposition~\ref{pro.rele}, on a:
\begin{enumerate}
\item un couple $(x,y)$ de points de $\A$ est
singulier si et seulement si il existe des relev\'es $\t x$ et $\t y$ de
$x$ et $y$ tels que $(\t x,\t y)$ soit singulier pour $\t h$~;
\item le relev\'e canonique $\t h$ n'est pas conjugu\'e \`a une translation affine. En
particulier, il existe des couples singuliers pour $\t h$~;
\item les autres relev\'es $h'$ de $h$ sont tous conjugu\'es \`a une
translation affine~;
\item si $\gamma$ est un arc de translation pour $h$, alors tout
relev\'e $\t \gamma$ est un arc de translation pour le relev\'e canonique
$\t h$, et est libre pour les autres relev\'es de $h$.  
\end{enumerate}
\end{prop}
\begin{demo}
Pour le premier point, la propri\'et\'e essentielle est que dans les coordonn\'ees de
la figure~\ref{fig42}, toute orbite pour $\t h$ est enti\`erement
incluse dans une bande verticale de largeur $2$. 

Le deuxi\`eme point peut se d\'eduire du lemme~\ref{lem.indi}  (l'existence de couples singuliers
provient alors d'un r\'esultat de B.~Ker\'ekj\'art\'o, voir par exemple~\cite{lero4},
proposition~8).

Le troisi\`eme  point d\'ecoule imm\'ediatement de l'in\'egalit\'e $(**)$ (preuve
de l'unicit\'e dans la proposition~\ref{pro.rele}).

Le dernier point  d\'ecoule du point~2 de la proposition
\ref{pro.rele} concernant les arcs libres.

Les d\'etails sont laiss\'es au lecteur.
\end{demo}

Une cons\'equence de ce compl\'ement est qu'il existe aussi des couples
de points de $\A$ singuliers pour $h$. Ceci peut \'egalement se prouver directement, selon
la strat\'egie suivante~: dans le cas contraire, tous les points de
l'anneau $\A$ ont le m\^eme ensemble $\omega$-limite, qui est l'un des
deux points fixes (par exemple $N$). Dans ce cas, un argument
c\'el\`ebre de Birkhoff montre qu'il existe un disque attractif contenant
le point $N$ (voir par exemple \cite{lero9}).
 D'apr\`es le lemme~\ref{lem.diat}, l'indice de $N$ est alors 
\'egal \`a $1$.

\clearpage
\part{Applications en dynamique locale}

Nous d\'eduisons du
th\'eo\-r\`eme d'extension~\ref{the.exte}, du th\'eo\-r\`eme d'existence de
d\'ecom\-po\-sition en briques~\ref{the.exde} et du th\'eo\-r\`eme principal~\ref{the.prinbis}
les th\'eor\`emes de dynamique locale \ref{the.stlo},  \ref{the.fleu}
et \ref{the.brow}. 
La preuve du th\'eo\-r\`eme~\ref{the.fleu} utilisera aussi les r\'esultats
\'el\'ementaires sur les d\'ecom\-po\-sitions en briques~; celle du
th\'eo\-r\`eme~\ref{the.brow} n\'ecessite de plus un certain nombre des
r\'esultats interm\'ediaires de ce texte. \emph{Les sections~\ref{sec.fleu},
\ref{sec.bran} et \ref{sec.iter} sont compl\`etement ind\'ependantes}.

\section{Th\'eor\`eme de la fleur topologique}\label{sec.fleu}

\begin{demo}[du th\'eor\`eme \ref{the.fleu}]
Soit $h:U \rightarrow V$ un hom\'eo\-mor\-phisme local  dont le point
fixe est d'indice strictement sup\'erieur \`a $1$.
On applique le th\'eor\`eme d'extension \ref{the.exte},  et on se place
dans le mod\`ele de la sph\`ere, mais en identifiant $(0,0)$ \`a $S$ et l'infini \`a
$N$. On obtient ainsi un \homeo\ de la sph\`ere fixant uniquement $N$
et $S$, que l'on note encore $h$,  
avec $\indi(h,S)=1+p>1$. L'hom\'eo\-mor\-phisme $h$ v\'erifie alors les
hypoth\`eses du th\'eor\`eme principal, et il suffit maintenant de montrer que tout
voisinage de $S$ contient $p$ p\'etales attractifs et $p$ p\'etales
r\'epulsifs, d'intersections deux \`a deux \'egales \`a $\{S\}$, cycliquement altern\'es.

Soit $F$ une d\'ecom\-po\-sition en briques pour $h$ (th\'eo\-r\`eme~\ref{the.exde}).
D'apr\`es le th\'eor\`eme \ref{the.prinbis} (``version simpliciale'' du
th\'eo\-r\`eme principal), il existe $p$ croissants
attractifs \`a dynamique Nord-Sud, simpliciaux, et $p$ croissants r\'epulsifs \`a
dynamique Sud-Nord, \'egalement simpliciaux, deux \`a deux d'int\'erieurs
disjoints, et cycliquement altern\'es autour de
$S$. Par d\'efinition, un croissant attractif \`a dynamique Nord-Sud
contient un p\'etale simplicial
attractif bas\'e en $S$ ; de m\^eme un croissant r\'epulsif \`a dynamique Sud-Nord contient un
p\'etale simplicial r\'epulsif bas\'e en $S$.

Pour obtenir l'existence des $2p$ p\'etales
\emph{dans tout voisinage de $S$}, il suffit donc de montrer qu'un
p\'etale attractif simplicial  contient des p\'etales attractifs arbitrairement
proches de $S$~; plus pr\'ecis\'ement, il nous reste \`a prouver~:
\begin{lemm}
Soit $h$   un hom\'eo\-mor\-phisme  de la sph\`ere, pr\'eservant l'orientation,
fixant uniquement les deux points $N$ et $S$, tel que $\indi(h,N)\neq
1$. Soit $F$ une d\'ecom\-po\-si\-tion en briques pour $h$, et $O_S$ un
voisinage de $S$.
  Alors tout p\'etale attractif simplicial $P$ bas\'e en $S$  contient un p\'etale
  attractif simplicial bas\'e en $S$ et inclus dans $O_S$.
\end{lemm}
Soit $K=P \setminus O_S$~; il s'agit de montrer qu'il existe un p\'etale attractif
simplicial, inclus dans $P$, et disjoint de $K$.
On peut toujours supposer que $K$  est 
r\'eunion d'un nombre fini de briques de la
d\'ecom\-po\-si\-tion, et que chaque composante connexe de $K$ rencontre le
bord de $P$ (quitte \`a grossir $K$).

Puisque $K$ rencontre le bord de $P$, il existe  une brique $B_1$
 incluse dans $K$ et rencontrant le bord de
$P$. Comme $P$ est un attracteur strict r\'eunion de briques, $A^+(B_1)
\subset P$, et $\Delta(B_1)$ est \`a extr\'emit\'es Sud ; de plus,
on voit facilement que la brique $B_1$ ne peut qu'\^etre de type dynamique $\fs$.
 La droite $\Delta(B_1)$ d\'elimite donc
un p\'etale attractif $P_1=P^+(B_1)$, inclus dans $P$, et ne contenant pas
$B_1$. 

On applique ce proc\'ed\'e autant de fois que n\'ecessaire
jusqu'\`a obtenir un p\'etale attractif $P_k$ disjoint de $K$ (plus
pr\'ecis\'ement, on peut raisonner par r\'ecurrence sur le nombre
de briques contenues dans $K$, en remarquant que, \`a la $k$\`eme \'etape,
 chaque composante connexe de $K \cap P_k$ rencontre le bord de $P_k$).
\end{demo}

\section{Branches stables et instables locales}\label{sec.bran}
\index{branche stable}
\begin{demo}[du th\'eor\`eme~\ref{the.stlo}]
Soit $h:U \rightarrow V$ un hom\'eo\-mor\-phisme local dont le point
fixe est d'indice strictement inf\'erieur \`a $1$.
D'apr\`es le th\'eor\`eme d'extension \ref{the.exte}, il suffit de
traiter le cas o\`u l'hom\'eo\-mor\-phisme local $h$ est  la restriction
\`a un voisinage de $N$ (que l'on notera encore $U$)  d'un \homeo\ de
la sph\`ere fixant uniquement $N$ et $S$ (que l'on notera encore $h$), 
avec $\indi(h,N)=1-p<1$. 

D'apr\`es le th\'eor\`eme global \ref{the.prin}, il existe $p$ croissants
attractifs \`a dynamique Nord-Sud et $p$ croissants  r\'epulsifs \`a
dynamique  Sud-Nord, deux \`a deux d'intersection \'egale \`a $\{N,S\}$,
cycliquement altern\'es autour de $N$.

Soit $A$ l'un des croissants attractifs \`a dynamique Nord-Sud. 
Par d\'efinition, il existe un p\'etale attractif
 $P$ en $S$, inclus dans $A$,  tel que
 $A \setminus P \subset U$.  
On peut de plus supposer que la situation est hom\'eomorphe au dessin de
 la figure \ref{fig34}~(a), c'est-\`a-dire que l'ensemble
 $\alpha=\partial P \setminus \partial A$
 est connexe~: en effet, dans le cas contraire, on peut se ramener \`a
 cette situation en modifiant l\'eg\`erement le p\'etale $P$, comme
indiqu\'e sur la partie (b) de la figure.

\begin{figure}[htpb]
\par
\centerline{\hbox{\input{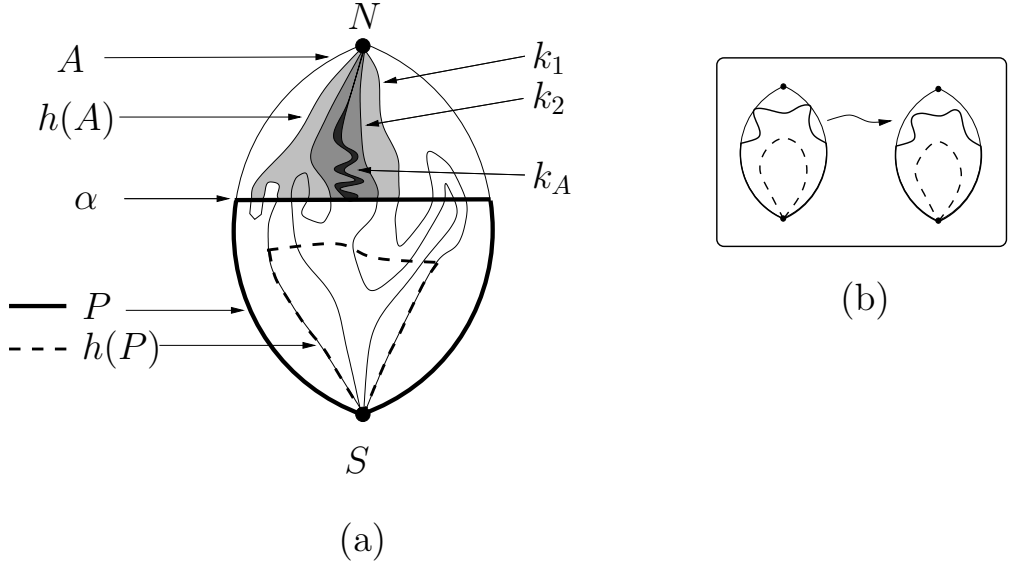}}}
\par
\caption{\label{fig34}Construction d'une branche instable locale}
\end{figure}

On pose alors  $W_{k_A}=\adhe(A \setminus P)$ (c'est un disque
topologique) ; et
pour tout $n\geq 0$, on appelle $k_n$ la composante connexe de $h^n(A)
 \cap W_{k_A}$ qui contient $N$. On pose $k_A=\cap_{n\geq 1}k_n$.
\begin{affi}
 Soit $U'$ un disque topologique ferm\'e, voisinage de $N$,
 tel que $U' \cap A=W_{k_A}$~; alors
l'ensemble $k_A$ est une branche instable $U'$-locale.
\end{affi}
\begin{demo}
Nous prouvons successivement que $k_A$ v\'erifie les quatre points de la
d\'efinition des branches instables $U'$-locales
(d\'efinition~\ref{def.bran}).

L'ensemble $k_n$  contient $N$, et on voit facilement qu'il
 rencontre $\alpha$. Il est  compact, connexe,  plein  (une composante
 connexe d'un ensemble plein est pleine).
 On a aussi  $h^{n+1}(A) \subset h^n(A)$, donc $k_{n+1}\subset k_n$.

 Par cons\'equent l'ensemble $k_A$ contient $N$, rencontre $\alpha$ et
est inclus dans $U'$ (ceci prouve le premier point).
Il est  compact et connexe comme intersection d\'ecroissante d'ensembles
 compacts et connexes~; il est plein comme intersection d'ensembles
 pleins. Ceci prouve le deuxi\`eme point.

Montrons que $k_A \subset h(k_A)$. Tout d'abord, $k_1 \subset
h(W_{k_A})$ ; en effet, on a 
$$
\begin{array}{rcl}
k_1 \subset h(A) \cap (W_{k_A}) & =  & h(W_{k_A} \cup P ) \cap W_{k_A}           \\
			& = & (h(W_{k_A}) \cap W_{k_A} )\ \cup \ (h(P) \cap
			W_{k_A}) \\
			& = &  	h(W_{k_A}) \cap W_{k_A}
\end{array}
$$
car $h(P) \cap W_{k_A}=\emptyset$ ($P$ est un attracteur strict).
  On a donc $k_{n+1}\subset k_1 \subset
h(W_{k_A})$ ; $k_{n+1}$ est un ensemble connexe inclus dans $h^{n+1}(A)
\cap h(W_{k_A})$ et qui contient $N$. Or $h(k_n)$ est la composante
connexe de $h^{n+1}(A) \cap h(W_{k_A})$ contenant $N$ ; on a donc $k_{n+1}
\subset h(k_n)$. On en d\'eduit $k_A \subset h(k_A)$.

Soit $x$ un point de $k_A$. Ses it\'er\'es n\'egatifs  sont tous dans $k_A$,
donc ne peuvent que tendre  vers $N$ d'apr\`es le corollaire
~\ref{cor.erre}. Ceci prouve le troisi\`eme point.

L'ensemble $W_{k_A}$ est un voisinage de $k_A \setminus \{N\}$ dans $U'$.
Soit $y$ un point de $W_{k_A}$.
D'apr\`es la d\'efinition des croissants
attractifs \`a dynamique Nord-Sud,  il existe un p\'etale attractif bas\'e en $S$ contenant
$y$, par cons\'equent les it\'er\'es positifs de $y$ tendent vers $S$, et
$y$ poss\`ede un it\'er\'e positif hors de $U'$.
Ceci prouve le quatri\`eme point.  
\end{demo}

On effectue la m\^eme construction pour chaque croissant attractif \`a
dynamique Nord-Sud
$A_1, \cdots A_p$,  et la construction similaire pour chaque croissant
r\'epulsif \`a dynamique Sud-Nord $A'_1, \cdots A'_p$ (en rempla\c{c}ant $h$ par $h^{-1}$).
Ceci  produit des ensembles $W_{k_{A_i}}$, $W_{k_{A'_i}}$ et des compacts
$k_{A_i}$, $k_{A'_i}$.  

On trouve facilement  un disque topologique ferm\'e $U'$, voisinage de $N$,
inclus dans $U\cap h(U)$, tel que pour tout $i$, $U' \cap A_i=W_{k_{A_i}}$ et $U' \cap
A'_i=W_{k_{A'_i}}$. Les compacts $k_{A_i}$ et $k_{A'_i}$ sont des branches
instables et stables $U'$-locales v\'erifiant les propri\'et\'es voulues.
\end{demo}

\section{Indices des it\'er\'es}\label{sec.iter}
Dans cette section, nous montrons que pour un
hom\'eo\-morphisme local d'indice diff\'erent de $1$, les indices de $h^n$
et de $h$ sont \'egaux pour tout $n$ non nul. Le point~1 du th\'eor\`eme,
qui revient \`a montrer que l'indice de $h^n$ est bien
d\'efini pour tout entier $n$ non nul,  d\'ecoule du lemme (bien
connu) suivant~:
\begin{lemm}[absence de r\'ecurrence locale]
Soit $h:U \rightarrow V$ un hom\'eo\-mor\-phisme local fixant uniquement
$0$, tel que l'indice du point fixe est diff\'erent de $1$. Il existe
alors un ensemble $U' \subset U$, voisinage de $0$,  ayant la propri\'et\'e suivante. Pour tout
point $x$ de $U'$, si la suite $(h^n(x))_{n \geq 0}$ est d\'efinie pour
tout $n$ positif et
enti\`erement incluse dans $U'$, alors cette suite tend vers $0$. 

En particulier, pour tout entier $n >0$, il existe un voisinage $U''$
de $0$ tel que l'hom\'eo\-mor\-phisme local $h^n : U'' \rightarrow h^n(U'')$
n'a pas d'autre point fixe que $0$.
\end{lemm}
Remarquons que l'exemple de la figure~\ref{fig49} montre qu'on ne peut pas
toujours choisir $U'=U$.

\begin{demo}[du lemme]
Le th\'eor\`eme d'extension~\ref{the.exte} permet de supposer que $h$ est
un hom\'eo\-mor\-phisme d\'efini sur le plan $\R^2$, et fixant uniquement $0$.
Dans ce cadre, l'hypoth\`ese d'indice permet d'appliquer la th\'eorie de \hbox{Brouwer}~: 
le lemme d\'ecoule alors du corollaire~\ref{cor.erre}.
\end{demo}

Passons maintenant au calcul de l'indice des it\'er\'es de $h$.  
L'id\'ee de la d\'emons\-tra\-tion est simple. Comme toujours, il suffit de
faire la preuve dans le cadre global. A l'aide des r\'esultats de ce
texte, nous pouvons d\'ecouper la sph\`ere en un certain nombre de
croissants attractifs et r\'epulsifs, \`a dynamique Nord-Sud ou Sud-Nord,
d'int\'erieurs disjoints deux \`a deux. On voit assez facilement que les
indices partiels pour $h$ et pour $h^n$ \`a travers chaque croissant
co\"{\i}ncident. D'autre part, entre deux croissants adjacents,
l'indice partiel pour $h$ est nul, et il
existe une cha\^{\i}ne de disques pour $h$, gr\^ace \`a laquelle on peut montrer
que l'indice partiel pour $h^n$ est encore nul. La relation de Chasles
 permet de conclure.

En pratique, on va se ramener, \emph{via} le relev\'e canonique, au
cadre des hom\'eo\-morphismes de \hbox{Brouwer}.
On commence donc par prouver une version du th\'eo\-r\`eme dans ce
cadre. Remarquons que l'une des difficult\'es de la preuve concerne le
changement de cadre~: en effet, si $H$ d\'esigne le relev\'e canonique de $h$, 
il n'est pas clair \emph{a priori} que $H^n$ soit le relev\'e canonique
de $h^n$ (tant qu'on ne sait pas quel est l'indice
de $h^n$, le relev\'e cnonique n'est pas d\'efini...).
 Ceci sera par contre  un corollaire imm\'ediat du th\'eor\`eme et
de la proposition~\ref{cor.iphn}.

Rappelons qu'une droite de \hbox{Brouwer} pour $H$ est encore une
droite de \hbox{Brouwer} pour toute les puissances non nulles de $H$ (remarque~\ref{rem.dbli}).

\index{indice partiel}
\begin{prop}[Indices partiels des it\'er\'es]
\label{cor.iphn}
Soit $H$  un hom\'eo\-morphisme de \hbox{Brouwer}, et $F$ une d\'ecom\-po\-sition en
briques pour $H$. Soit $(\Delta_0,\Delta_1)$ un couple de droites de \hbox{Brouwer}
simpliciales disjointes. Alors pour tout entier $n$ non nul, les
indices partiels  $\inpa(H,\Delta_0,\Delta_{1})$ et
$\inpa(H^n,\Delta_0,\Delta_{1})$ sont \'egaux.
\end{prop}

On se place sous les hypoth\`eses de la proposition~\ref{cor.iphn}, et
on consid\`ere  une famille v\'erifiant la
conclusion du lemme~\ref{lem.codb1}, autrement dit une famille maximale de
sous-croissants simpliciaux minimaux de la bande
$\adhe(D(\Delta_0,\Delta_1))$,
$${\cal F}=\{\adhe(D(\Delta'_0,\Delta'_1)),\dots,
\adhe(D(\Delta'_{2k},\Delta'_{2k+1}))\}.
$$
On a alors~:
\begin{lemm}[compl\'ement \`a la proposition~\ref{pro.codb2}]
\label{lem.iphn}
Pour tout $i=-1, \dots, 2k+1$,  pour tout entier $n$ non nul, on a
$$\inpa(H,\Delta'_i,\Delta'_{i+1})=\inpa(H^n,\Delta'_i,\Delta'_{i+1}).
$$
\end{lemm}

\begin{demo}[du lemme \ref{lem.iphn}]
Il suffit de traiter le cas $n>0$, puisque l'indice partiel pour $H$
est toujours \'egal \`a l'indice partiel pour $H^{-1}$
(affirmation~\ref{aff.ipfa}).

Soit $C$ un croissant de la famille $\cal F$. Le croissant $C$ \'etant minimal, d'apr\`es la
proposition~\ref{pro.cam}, il est \`a dynamique  $\t S$-$\t N$ ou $\t
N$-$\t S$. Il est clair que $C$ est encore un croissant pour $H^n$,
ayant le m\^eme type dynamique que pour $H$ (on peut par exemple
utiliser l'affirmation~\ref{aff.topo}). Puisque l'indice partiel
entre les deux bords de $C$ ne d\'epend que du type dynamique de $C$
(affirmation~\ref{aff.ipcr}), les indices partiels pour $H$ et pour
$H^n$ sont \'egaux. Ceci montre le lemme pour les entiers $i$ pairs.

\index{droites topologiques!qui se touchent}\index{cha\^{\i}ne!et indice partiel}
Soit maintenant $i$ un entier impair. Si les droites $\Delta'_{i}$ et
$\Delta'_{i+1}$ ne sont pas disjointes, elles se touchent
(d\'efinition~\ref{def.trav}), et l'indice partiel est nul pour $H$
comme pour $H^n$ (par d\'efinition). Il reste le cas o\`u ces deux droites
sont disjointes~; il s'agit de montrer que l'indice partiel pour $H^n$
est encore nul. 

Le couple $(\Delta'_i,\Delta'_{i+1})$ est indiff\'erent,
on suppose par exemple que $\Delta'_i$ est attractive.
D'apr\`es la proposition~\ref{pro.codb2}, l'indice partiel est nul
pour $H$, et il existe une cha\^{\i}ne
de pseudo-disques  $(D_1, \dots D_p)$, pour  $H$, de $\Delta'_{i}$ \`a
$\Delta'_{i+1}$ (constitu\'ee de briques).
Notons $(n_i)_{i=1,\dots,p-1}$ les temps de
transition de la cha\^{\i}ne.
Consid\'erons alors  la suite
$$(D_1, H^{(n-1) \times n_1}(D_2),
\dots,  H^{(n-1) \times (n_1+\cdots +n_{p-1})}(D_p)).$$
Les disques topologiques qui la composent sont deux \`a deux disjoints,
sans quoi il existerait une cha\^{\i}ne de pseudo-disques p\'eriodique pour
$H$ (ce qui contredirait le lemme de Franks~\ref{lem.fran}).
 De plus, cette suite est une
cha\^ine de pseudo-disques pour $H^n$, avec les m\^emes temps de transition
$(n_i)$. Cette cha\^ine de pseudo-disques va de la droite de \hbox{Brouwer} $\Delta'_i$ \`a la
droite $\Delta''_{i+1}=H^{(n-1) \times (n_1+\cdots +n_{p-1})} (\Delta'_{i+1})$.

  Le lemme de Franks pour l'indice partiel (proposition~\ref{pro.frip}) montre alors que
le nombre $\inpa(H^n,\Delta'_i,\Delta''_{i+1})$ est nul.

D'autre part, la droite $\Delta''_{i+1}$ est incluse dans le
domaine de \hbox{Brouwer} engendr\'e, pour $H^n$, par $\Delta'_{i+1}$.
 D'apr\`es le
deuxi\`eme point du corollaire \ref{cor.dotr}, on a 
$\inpa(H^n,\Delta'_{i+1},\Delta''_{i+1})=0$.

Les droites de \hbox{Brouwer} $\Delta'_i$, $\Delta'_{i+1}$ et
$\Delta''_{i+1}$ sont clairement disjointes, la deuxi\`eme s\'eparant les
deux autres. On peut donc appliquer la relation de Chasles
(lemme~\ref{lem.chas}), obtenant 
$\inpa(H^n,\Delta'_i,\Delta''_{i+1})=
\inpa(H^n,\Delta'_i,\Delta'_{i+1})+\inpa(H^n,\Delta'_{i+1},\Delta''_{i+1})$.
Deux des termes de la relation \'etant nuls, le troisi\`eme l'est aussi, ce qui
nous donne bien
$$
\inpa(H^n,\Delta'_i,\Delta'_{i+1})=0=\inpa(H,\Delta'_i,\Delta'_{i+1}).
$$
Ceci termine la preuve du lemme.
\end{demo}

\begin{demo}[de la proposition \ref{cor.iphn}]
\index{relation de Chasles!g\'en\'eralis\'ee}
 Comme dans la preuve de la formule de la proposition~\ref{pro.codb2},
on applique la relation de Chasles g\'en\'eralis\'ee
 (lemme~\ref{lem.chas2}), pour $H$ d'une part, et pour $H^n$ d'autre
 part, \`a la famille des droites $\Delta'_{-1}=\Delta_0,
 \dots, \Delta'_{2k+2}=\Delta_1$. D'apr\`es le lemme, les indices
 partiels entre deux droites successives sont les m\^emes pour
 $H$ et $H^n$, et on obtient le r\'esultat.
\end{demo}

\begin{demo}[du th\'eor\`eme \ref{the.brow}]
Gr\^ace au th\'eor\`eme d'extension, il suffit de prouver que si $h$ est
un hom\'eo\-mor\-phisme de la sph\`ere, pr\'eservant l'orientation,
fixant uniquement les deux points $N$ et $S$, tel que $\indi(N)=1-p
\neq 1$, alors pour
tout entier $n \neq 0$, $\indi(h^n,N)=\indi(h,N)$.

Comme dans la preuve du th\'eo\-r\`eme principal~:
\begin{itemize}
\item on se donne une d\'ecom\-po\-sition en briques $F$ pour $h$ (th\'eo\-r\`eme~\ref{the.exde})~;
\item on note $\Delta_0$ une droite de \hbox{Brouwer} simpliciale \`a
extr\'emit\'es Nord-Sud pour $h$ (son existence est assur\'ee par la
proposition~\ref{pro.dbns})~;
\item on note $H=\t h$ le relev\'e canonique de $h$ au plan $\t \A$, fourni par la
proposition~\ref{pro.rele}~;
\item on note  $\t \Delta_0$ un relev\'e de $\Delta_0$, et $\t \Delta_1=\tau(\t
\Delta_0)$ (o\`u $\tau$ d\'esigne l'automorphisme de rev\^etement).
\end{itemize}

 On a alors (comme dans l'affirmation~\ref{aff.padi})~:
$$\inpa(H,\t \Delta_0,\t \Delta_1)=\indi(H,\tau),$$
$$\inpa(H^n,\t \Delta_0,\t \Delta_1)=\indi(H^n,\tau).$$
D'autre part, on a aussi
$$\indi(h,N)=\indi(H,\tau)+1,$$
$$\indi(h^n,N)=\indi(H^n,\tau)+1.$$
En effet, la premi\`ere \'egalit\'e vient de la d\'efinition du relev\'e canonique. 
Montrons la deuxi\`eme. On reprend pour cela les id\'ees de la construction du relev\'e
canonique. Soit $n>0$ fix\'e~; puisque les indices sont des invariants
de conjugaison,
quitte \`a conjuguer $h$, on peut supposer que les $(2n+3)$ droites de
\hbox{Brouwer} $\Delta_0, \dots, h^{2(n+1)}(\Delta_0)$ sont des demi-droites
euclidiennes issues de $N$, d\'ecoupant le plan en $2n+3$ secteurs de
m\^eme taille (g\'en\'eralisant la figure~\ref{fig42}).
 Dans ce mod\`ele, on voit facilement que le relev\'e
canonique $H$ de $h$ v\'erifie l'in\'egalit\'e~:
$$
 \mbox{ Pour tout } \t X \in \t \A, \quad \quad |\t \theta (\t X)-\t \theta
(H^n(\t X))| < \frac{1}{2}.
$$
Comme $H^n$ est un relev\'e de $h^n$,
l'\'egalit\'e $\indi(h^n,N)=\indi(H^n,\tau)+1$  d\'ecoule alors du lemme~\ref{l.cones}.

Il ne nous reste plus qu'a montrer que 
$$\inpa(H,\t \Delta_0,\t \Delta_1)=\inpa(H^n,\t \Delta_0,\t \Delta_1).$$
Cette derni\`ere \'egalit\'e est une application de la
proposition~\ref{cor.iphn}.
\end{demo}

\pagebreak
\appendix
\addcontentsline{toc}{part}{Appendice}
\part*{Appendice}
\section{Th\'eor\`eme de Schoenflies-Homma et variantes}
\index{Schoenflies-Homma!th\'eo\-r\`eme de}
Le th\'eor\`eme de Schoenflies  est l'un des outils techniques les plus couramment utilis\'es dans les
probl\`emes de  topologie plane. Dans ce texte, nous appliquons un certain
nombre de variantes qui peuvent toutes \^etre vues comme des
corollaires d'un th\'eor\`eme de Homma. Dans cet appendice, nous \'enonçons
le th\'eor\`eme de Homma et expliquons quelques-unes de ses utilisations.

\subsection{Le th\'eor\`eme de Schoenflies}
\begin{theo}[Schoenflies]
Tout hom\'eo\-mor\-phisme entre deux courbes de Jordan incluses dans $\S^2$
s'\'etend en un hom\'eo\-mor\-phisme de $\S^2$.
\end{theo}

\subsection{Le th\'eor\`eme de Schoenflies-Homma}
 \begin{defi}[( \cite{homm2})] 
On appelle \emph{ensemble en forme de \og Y \fg}
 un ensemble qui est la r\'eunion de trois arcs ayant en commun un
 seul point, extr\'emit\'e de chacun d'eux.
 \end{defi}
 L'orientation de la
 sph\`ere induit un ordre cyclique sur l'ensemble des trois arcs qui
 constituent n'importe quel ensemble de ce type, au voisinage du point
d'intersection.

\begin{theo}[Homma]
\label{thhomm}
Soit  $F$ un ensemble compact, connexe et localement connexe de la
sph\`ere. Une application continue et injective de $F$ dans la sph\`ere
s'\'etend en un \homeo\ de la sph\`ere pr\'eservant l'orientation si
et seulement si elle pr\'eserve
 l'ordre cyclique de tous les ensembles en forme de \og Y \fg\
 contenus dans $F$.
\end{theo}

\subsection{Utilisations}
On utilise le th\'eor\`eme avec divers ensembles $F$, ces ensembles \'etant
\`a chaque fois des graphes simpliciaux finis. Les d\'etails sont laiss\'es
au lecteur. On a essentiellement
besoin de trois variantes~:
\begin{itemize}
\item Dans le cas o\`u $F$ est une courbe ferm\'ee simple ou bien un arc,
$F$ ne contient
aucun ensemble en forme de \og Y \fg, et on retrouve en particulier 
le th\'eor\`eme de Schoenflies.

\item Soient $\Delta_1, \dots , \Delta_k$ ($k \geq 2$), une famille de droites
topologiques du plan, deux \`a deux disjointes. Supposons que pour tout
 $0 \leq i_1 < i_2 < i_3 \leq k$, la droite $\Delta_{i_2}$ s\'epare
$\Delta_{i_1}$ de $\Delta_{i_3}$.
 Alors il existe un
hom\'eo\-mor\-phisme $h$ du plan, pr\'eservant l'orientation, qui envoie
chacune des droites $\Delta_i$ sur la droite verticale $i \times \R$.

\item  Soient $N$ et $S$ deux points antipodaux de la sph\`ere, et
$\Delta_1, \dots , \Delta_k$, $k \geq 2$ une famille d'arcs
d'extr\'emit\'es $N$ et $S$, deux \`a deux d'int\'erieurs disjoints. Alors il
existe un hom\'eo\-mor\-phisme de la sph\`ere, fixant $N$ et $S$, qui envoie
chacun des arcs sur un demi-grand cercle.
\end{itemize}

Remarquons que dans beaucoup de cas, on a besoin de sp\'ecifier l'image
d'un certain nombre de points par l'hom\'eo\-mor\-phisme recherch\'e~; ceci
arrive en particulier quand on veut obtenir un hom\'eo\-mor\-phisme du plan,
c'est-\`a-dire qui fixe le point \`a l'infini. Ceci
ne pose pas de probl\`eme (typiquement, on utilise que le groupe des
hom\'eo\-mor\-phismes du disque qui fixent chaque point du bord est transitif).

\pagebreak
\addcontentsline{toc}{part}{R\'ef\'erences}

\end{document}